\documentclass{article}[11pt]
\usepackage{graphicx} % Required for inserting images
\usepackage{latexsym, amssymb, amsmath,amsthm}
\usepackage[margin=1in]{geometry}
\usepackage{indentfirst}
\usepackage[font=footnotesize]{caption}
\usepackage{subcaption}
\usepackage[dvipsnames]{xcolor}
\usepackage{adjustbox}
\usepackage{lineno}
\usepackage[raggedrightboxes]{ragged2e}
\usepackage{comment}
\usepackage{soul}
\usepackage{url}
\usepackage{comment}
\usepackage[normalem]{ulem}
\usepackage{placeins}
\usepackage{multirow}
\usepackage{orcidlink}

\usepackage{tikz}
\usetikzlibrary{shapes,arrows}
\usetikzlibrary{matrix,chains,scopes,positioning,arrows,fit}
\tikzstyle{block1} = [rectangle, draw, fill=blue!20, text width=2em, text centered, rounded corners, minimum height=3em, minimum width=3em]
\tikzstyle{block2} = [rectangle, draw, fill=red!20, text width=2em, text centered, rounded corners, minimum height=3em, minimum width=3em]

\newtheorem{thm}{Theorem}

\usepackage{hyperref} 
\hypersetup{
    colorlinks=true,
    linkcolor=blue,
    filecolor=magenta,      
    urlcolor=cyan}

\usepackage{pdflscape}

\newcommand{\bigzero}{\mbox{\normalfont\Large\bfseries 0}}
\newcommand{\rvline}{\hspace*{-\arraycolsep}\vline\hspace*{-\arraycolsep}}

\newcommand{\OS}[1]{\textcolor{ForestGreen}{#1}}

\title{The behavioral spillover effect: Modeling behavioral interdependencies in multi-pathogen dynamics}
\author{Leah LeJeune$^{1,2,*}$\orcidlink{0000-0002-1150-7814}, Omar Saucedo$^{1,2,*, \diamond}$\orcidlink{0000-0003-3900-7712}, Lauren M. Childs$^{1,2,\sharp}$\orcidlink{0000-0003-3904-3895}, Navid Ghaffarzadegan$^{3,\sharp}$\orcidlink{0000-0003-3632-8588}}
\date{%
    $^1$Department of Mathematics, Virginia Tech\\
    225 Stanger Street Blacksburg, VA 24061-1026\\~\\
    $^2$Virginia Tech Center for the Mathematics of Biosystems, Virginia Tech\\~\\
    $^3$Department of Industrial and Systems Engineering, Virginia Tech\\
    7054 Haycock Rd, Falls Church, VA 22043\\~\\
    $^*$ These authors contributed equally.\\
    $^{\sharp}$ These authors also contributed equally. \\
    $\diamond$ Corresponding author: Omar Saucedo, osaucedo@vt.edu  \\%
}

\begin{document}
\maketitle

\begin{abstract}
    During the recent pandemic, a rise in COVID-19 cases was followed by a decline in influenza. In the absence of cross-immunity, a potential explanation for the observed pattern is behavioral: non-pharmaceutical interventions (NPIs) designed and promoted for one disease also reduce the spread of others. We study short-term and long-term dynamics of two pathogens where NPIs targeting one pathogen indirectly influence the spread of another – a phenomenon we term \textit{behavioral spillover}. We examine how perceived risk of and response to one disease substantially alters the spread of other pathogens, revealing how waves of different pathogens emerge over time as a result of behavioral interdependencies and human response. Our analysis identifies the parameter space where two diseases simultaneously co-exist, and where shifts in prevalence occur. Our findings are consistent with observations from the COVID-19 pandemic, where NPIs contributed to significant declines in infections such as influenza, pneumonia, and Lyme disease.
\end{abstract}


\section{Introduction}
The epidemic modeling literature has evolved over the past century, incorporating various mathematical techniques to capture the complex dynamics of infectious disease spread within social and environmental contexts \cite{kermack1927contribution,ross1917applicationIII,anderson1979population}. From early compartmental models based on ordinary differential equations to more sophisticated stochastic or network-based approaches, epidemiological models continuously advance to better represent real-world disease transmission patterns. The primary objective of these models are to enhance our understanding of disease dynamics \cite{blackwood2018introduction}, assess potential impacts of interventions to inform evidence-based policy-making \cite{mcbryde2020role,traulsen2023individual} and to offer short- and long-term projections \cite{cramer2022evaluation,rahmandad2022enhancing}. Recent advances include incorporation of change in human behavior and risk perception as well as mobility patterns in epidemiological models to further improve model-based insights \cite{lejeune2025formulating}.

Most studies focus on modeling the spread of a single virus, such as influenza or coronavirus \cite{coburn2009modeling,li2020substantial,he2020seir}; however, multiple pathogens often co-circulate, with the disease with the higher prevalence shifting over time \cite{neumann2022seasonality}. For example, influenza and respiratory syncytial virus frequently overlap during seasonal epidemics, with higher prevalence observed in colder months \cite{lofgren2007influenza}. In contrast, Lyme disease is more common during warmer seasons when outdoor activities increase \cite{moore2014meteorological}, and rhinoviruses peak in the spring and fall \cite{neumann2022seasonality}. This raises an important question: to what extent are the dynamics of multiple co-circulating pathogens independent of each other? This question became especially important during early months of the COVID-19 pandemic when a large fraction of public health capacity was devoted to treating COVID-19 patients \cite{rubin2020happens}.  

Despite the early concerns, an intriguing pattern emerged. As the coronavirus persisted over a relatively long period, a notable drop in influenza A and B cases was observed during the first two years of the pandemic \cite{jones2020covid,solomon2020influenza,perez2020dramatic,kiseleva2021covid}.  Interestingly, as the number of COVID-19 cases declined in late 2022, thanks to natural and vaccine-induced immunity, the number of influenza type-A cases skyrocketed (Figure \ref{fig:covidfludata}).
A potential explanation for the observed decline in influenza cases is behavioral: the non-pharmaceutical interventions (NPIs) implemented to control COVID-19, such as mask-wearing, social distancing, and reduced mobility, also limited the spread of other pathogens \cite{solomon2020influenza}. Since influenza transmission is partly aerosol-based, mask-wearing and other NPIs directly reduced influenza transmission as well. A second-order effect produced more complicated dynamics, related to the fact that during the two-year period of fewer influenza cases, little natural immunity was built against the virus, which possibly lead to a large wave of the disease once COVID-19 reached its endemic state.

This and other similar scenarios underscore the importance of analyzing interactions among multiple pathogens rather than relying solely on single-pathogen models. Single-pathogen approaches may overlook critical insights into competition, coexistence, and broader public health implications. Understanding the dynamics of multiple pathogens—including factors like cross-immunity and evolutionary processes—is crucial for capturing the complexity of real-world epidemics \cite{welsh2002no,bhattacharyya2015cross}. For example, some studies have explored how partial immunity from one infection can influence the spread of other pathogens (e.g., \cite{bhattacharyya2015cross,gog2002dynamics,lejeune2023effect}). Our focus, however, is distinct: we analyze how human behavior mediates the inter-dependencies among pathogens.

\subsection{Cross immunity}
Cross-immunity between different strains of the same pathogen or across different pathogen, referred to as heterologous immunity, has long been a topic of interest in immunology \cite{welsh2002no,welsh2010heterologous}. Pathogens that trigger an immune response after infection can potentially provide at least partial cross-protection against strains of the same or closely related pathogens. This cross-protection can influence the spread of disease, with an epidemic of one strain suppressing the transmission of another \cite{welsh2002no}. Early in the COVID-19 pandemic, it was hypothesized that some level of cross-immunity might exist among coronaviruses, potentially reducing the spread or severity of SARS-CoV-2 infection \cite{yaqinuddin2020cross}.

The phenomenon of cross-immunity and evolution of various strains and mutations have been explored through mathematical modeling, e.g. \cite{bhattacharyya2015cross,gog2002dynamics,kryazhimskiy2007state,eletreby2020effects}. Bhattacharyya et al. \cite{bhattacharyya2015cross} examine cross-immunity among respiratory viruses, focusing specifically on RSV, HPIV, and hMPV. Using compartmental models, they analyze patterns of laboratory-confirmed infections from the western United States between 2002 and 2014, arguing that closely related paramyxoviruses can suppress each other’s transmission through immune responses triggered by prior infections. In foundational work, Gog and Grenfell \cite{gog2002dynamics}, showed how cross-immunity among multiple pathogens can produce a range of complex, yet realistic behavior.

It is notable that model-based studies of cross-immunity often modify conventional compartmental models to represent complicated states of immunity or partial immunity to one strain but susceptibility to others \cite{lejeune2023effect,nikin2018unraveling}. However, based on recent pandemic experience, it is argued that such models often oversimplify or omit human behavior -- either by assuming no variation in contact rates over the course of an epidemic or by formulating them exogenously \cite{lejeune2025formulating}. In most contexts, fear of infection can lead to protective measures that suppress disease transmission, and as such  epidemic models should better incorporate human behavior change during a course of a pandemic \cite{funk2010modelling,funk2015nine,verelst2016behavioural}. For example, during the recent COVID-19 pandemic, as case numbers rose, public risk perception also increased, prompting individuals to adjust their daily routines to reduce the risk of infection. It also influenced government policies, leading to more stringent measures \cite{lim2023similar}. Altogether, these responses influenced behavior with actions such as quarantining, social distancing, and mask-wearing. In turn, this response reduced transmission, leading to fewer cases, which subsequently lowered risk perception and contributed to a resurgence of the disease. This phenomenon illustrates a negative (balancing) feedback loop—one of several mechanisms that can lead to infection waves \cite{rahmandad2022missing}.

\subsection{Modeling Human behavior}
There are various approaches to modeling the spread of infectious disease, ranging from mechanistic models—including deterministic and stochastic compartmental models \cite{blackwood2018introduction, allen2008introduction}—to curve-fitting approaches such as machine learning \cite{kraemer2025artificial, dandekar2020machine}. A comprehensive study of models used during the COVID-19 pandemic by Rahmandad et al. \cite{rahmandad2022enhancing} finds that mechanistic compartmental models outperform curve-fitting approaches, particularly in long-term projections. Importantly, their results show that incorporating behavioral feedback—where contact rates adjust in response to changes in prevalence or mortality—substantially improves projection accuracy.

To capture this effect, Rahmandad et al. \cite{rahmandad2022enhancing} propose the SEIRb model, which extends the standard Susceptible-Exposed-Infectious-Recovered (SEIR) framework by including a risk-responsive behavioral feedback mechanism. This allows infectivity to vary dynamically as human behavior adapts to the evolving state of the disease. As a result, even relatively simple mechanistic models can outperform more complex alternatives in long-term forecasting. In a related study, LeJeune et al. \cite{lejeune2024mathematical} develop a similar framework in which infectivity is explicitly modeled as a function of disease prevalence, with a delay reflecting lagged behavioral response. Their results demonstrate how such formulations improve the ability of models to capture realistic epidemic waves across data from 50 U.S. states and Washington, D.C.

These studies illustrate the growing emphasis on incorporating human behavior into epidemic models. More broadly, LeJeune et al. \cite{lejeune2025formulating} classify such approaches into exogenous and endogenous formulations. Exogenous approaches treat behavioral changes as externally imposed parameter shifts, whereas endogenous approaches model quantities such as contact rate or infectivity as functions of the disease state itself. Within endogenous frameworks, several extensions have been proposed. These include modeling heterogeneous compliance by dividing the susceptible population into behavioral subgroups \cite{abbas2022evolution,pant2024mathematical}, as well as representing individuals as payoff maximizers who balance economic costs against infection risks \cite{saad2023dynamics,espinoza2024adaptive}. Additional approaches incorporate the spread of fear, evolving social norms, and adherence fatigue \cite{epstein2021triple,qiu2022understanding,osi2025simultaneous}, or use game-theoretic formulations to capture strategic behavioral adaptation \cite{huang2022game, chang2020game}.

While these works highlight recent advancements in the consideration of cross-immunity and human behavior in epidemiological models, we emphasize here a critical gap which persists in the literature. On one hand, models of multi-pathogen dynamics often lack behavioral mechanisms; on the other hand, many behavioral epidemic models focus on single-pathogen contexts, overlooking how individual and collective responses influence outbreak dynamics. Particularly, compliance with NPIs can simultaneously influence the spread of multiple pathogens with similar transmission mechanisms. 

\subsection{Appearance in data}
\begin{figure}
    \centering
    \includegraphics[width=0.75\linewidth]{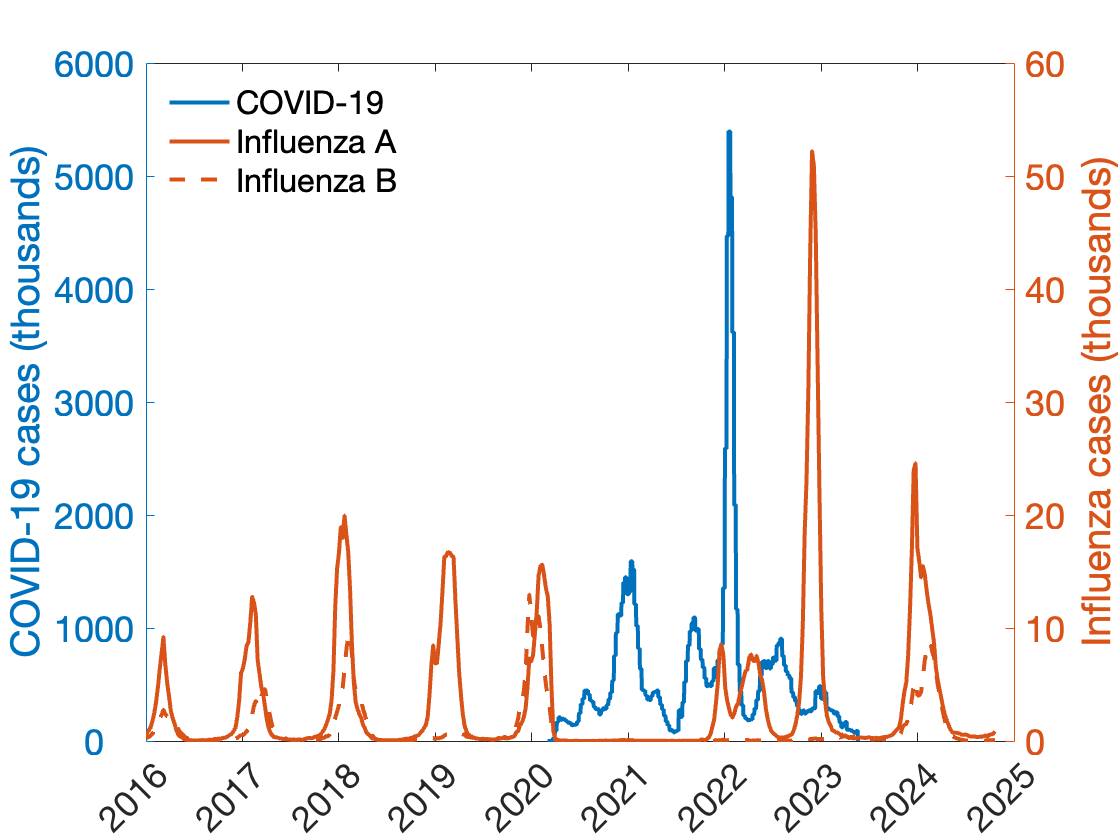}
    \caption{US Confirmed weekly cases of COVID-19 (blue, left y-axis) and Influenza A and B (red, right y-axis).  Source: based on CDC Influenza data \cite{CDC_FluView_YYYY} and Johns Hopkins Dashboard COVID-19 data \cite{dong2020interactive}.}
    \label{fig:covidfludata}
\end{figure}

In real-world settings, public behavior is rarely pathogen-specific. For instance, increasing mask-wearing or social distancing in response to one disease can inadvertently suppress the transmission of others. This was clearly observed during the COVID-19 pandemic, when widespread NPI substantially reduced flu cases \cite{jones2020covid,solomon2020influenza,perez2020dramatic}. As shown in Figure \ref{fig:covidfludata}, the United States has historically experienced annual flu seasons, typically peaking in late fall or early winter, until the onset of the recent COVID-19 pandemic.  During the surge in COVID-19 cases in 2020 and 2021, influenza cases dramatically declined, with very few cases of Influenza A or B. Interestingly, as COVID-19 was brought under control and transitioned into an endemic phase, seasonal flu patterns re-emerged. In fact, Flu A cases surged during the 2022–2023 season with a peak that was not only substantially larger but also arrived earlier than those observed in pre-pandemic years. This pattern is often attributed to behavioral factors: compliance with NPIs not only curbed COVID-19 transmission but also suppressed flu outbreaks. Once such measures were relaxed, flu incidence rebounded. This heightened rebound can be partly attributed to a decline in population immunity, as the prolonged suppression of influenza circulation during the pandemic left a larger-than-usual proportion of the population with reduced or waning protection against circulating strains \cite{chen2025immunity}.

Both COVID-19 and influenza have similar transmission mechanisms which can be targeted by comparable prevention measures, likely contributing to the reduction in influenza cases observed during the initial outbreak of COVID-19. However, examining the literature and data of other infectious diseases with distinct transmission routes, such as COVID-19 and Lyme disease, we see that changes in behavior in response to COVID-19 may have also caused changes in the spread of non-respiratory diseases. For example, the lifestyle disruption that COVID-19 caused across the world resulted in changes in time spent outdoors, affecting potential for human exposure to tick-borne diseases \cite{borșan2021recreational, jones2024lyme, mccormick2021effects}. Some regions reported an increase in tick-borne encephalitis, possibly due to more time spent out-of-doors \cite{jore2023outdoor}. Other regions, though, saw a reduction in reported cases of Lyme disease, potentially resulting from disruptions in the healthcare workflow and healthcare-seeking habits which led to under-reporting of cases \cite{jones2024lyme, mccormick2021effects, novak2021lyme}. In either case, these correlations further support the hypothesis that disruptions in habits and routines - i.e., changes in public behavior - in response to one disease can affect the spread of another disease with an independent transmission pathway. %Additionally, disruptions from COVID-19 on the healthcare workflow along with routines for seeking healthcare may explain under-reporting of Lyme disease cases .}

\subsection{Beyond specific diseases}
Similar dynamics are likely to arise in a variety of contexts and among pathogens with comparable transmission characteristics. For example, when individuals increase hand washing in response to a surge in one particular disease, they may also inadvertently reduce the transmission of other pathogens that spread through similar pathways. As such, modeling frameworks that incorporate behavioral coupling between co-circulating pathogens can offer unique insights into epidemic dynamics. These models can help answer the following important questions: How do general patterns of disease spread emerge when multiple pathogens are active and NPIs targeting one also reduce the transmission of another? How might these patterns change if such spillover effects are imperfect or only partially suppress the spread of a second disease? Can we identify ranges of reproductive numbers where one pathogen exceeds the other or where multiple pathogens can coexist?
Addressing these questions is essential for improving our understanding of—and response to—the complex interplay between human behavior and the concurrent spread of multiple infectious diseases.

In this paper, we study the short-term and long-term dynamics of two pathogens spreading simultaneously, where NPIs targeting one pathogen can indirectly influence the spread of others—a phenomenon we term \textit{behavioral spillover}. Generally, real-world infectious disease dynamics are influenced by many time-varying factors, including seasonality, policy changes, behavioral responses, emerging variants, environmental factors, and vaccination campaigns. The objective of this study, however, is not to reproduce the full complexity of epidemic dynamics, but rather a theoretical investigation: to isolate and analyze a specific mechanism—behaviorally mediated interactions between diseases—within a simple and analytically tractable framework. To that end, we present a mathematical model of two interacting pathogens and analyze their outbreak dynamics. The interaction is solely via human behavior; that is, no cross-immunity is considered, but we model how human risk response to one can fully or partially affect the spread of the other pathogen. We postulate that there are two regimes -- persistence by only one pathogen or coexistence of both -- and examine tipping points for change in patterns. This equilibrium is also influenced by the basic reproduction numbers of the pathogens. Over time, we observe a likelihood of periodic shift in which disease has higher prevalence between pathogens.

\section{Materials and Methods}
We develop a modified Susceptible-Infectious-Removed (SIR) model of the spread of two pathogens (denoted as $\{ A, B\}$) with no cross-immunity and include endogenous behavioral responses as shown in Figure \ref{fig:Eqn1FlowDiagram}. To begin, we formulate the spread of each pathogen independently (excluding the green curved dashed lines) and then formulate the behavioral interconnection between the spread of the pathogens, considering both perfect and imperfect interconnection.  

The core of our model structure is the simple SIR model with two modifications that include behavioral response and waning immunity. We refer to this family of models as SIRSb models to indicate the return of those in the recovered compartment to the susceptible compartment, and the inclusion of behavioral feedback that affects infectivity \cite{lejeune2024mathematical}.

The system of ordinary differential equations, System \eqref{eqn:SIRS}, represents an SIRS model of disease $i$ (here we focus on two diseases, $i \in\{A,B\})$ with waning immunity given by
\begin{equation}
\begin{aligned}\label{eqn:SIRS}
& \frac{d S_i}{d t}=-\beta_i(\cdot) S_i I_i+\frac{R_i}{\tau_R}, \\
& \frac{d I_i}{d t}=\beta_i(\cdot) S_i I_i-\frac{I_i}{\tau_I}, \\
& \frac{d R_i}{d t}=\frac{I_i}{\tau_I} -\frac{R_i}{\tau_R},\\
&\frac{d\widetilde{I}_i}{dt}=\frac{I_i-\widetilde{I}_i}{\tau_P}\\
& S_i+I_i+R_i=1,
\end{aligned} 
\end{equation}
where $S_i$, $I_i$, and  $R_i$ represent the proportion of individuals in the total population who are susceptible to disease $i$, infectious with disease $i$, and recovered from disease $i$, respectively. Notice that there is overlap among the two diseases in the population. The susceptible, infected, and recovered, populations with disease $A$ respectively, encompass the entire population, as do the susceptibles, infected, and recovered individuals with disease $B$. This stems from the assumption that all individuals have the potential to be infected with disease $A$, disease $B$, or both diseases simultaneously; in other words, individuals can become infected with both diseases without any biological cross-immunity. %For a given disease, an individual's infection status with respect to that disease has no influence on their infection status with respect to the other disease. 
For example, a disease $A$ infectious individual is also simultaneously one of susceptible, infectious or recovered from disease $B$.  Furthermore, infectious individuals do not change their behavior upon infection (e.g., no quarantine or isolation is included).
%, and individuals may be infected by both pathogens without any biological cross-immunity.}

Susceptible individuals become infected with disease $i$ through interactions with infectious individuals of disease $i$ with transmission rate $\beta_{i}(\cdot)$, which can depend on other variables. Infectious individuals recover after an average period of $\tau_I$ days. Once recovered, individuals lose immunity and are once again susceptible after on average $\tau_R$ days.

Here, $\widetilde{I}_i$ represents the perceived number of disease $i$ infectious individuals, which follows the true value $I_i$ with an average delay $\tau_P$. Specifically, $\widetilde{I}_i$ adjusts toward $I_i$ at rate $1/\tau_P$, so that discrepancies between perceived and actual infections decay over time \cite{nerlove1958, sterman2000}. This formulation corresponds to a standard exponential smoothing (information delay) process previously used in the literature to represent information flow (e.g., \cite{thompson2007polio, lim2023similar}).

\begin{figure}[h!]
\centering
\begin{tikzpicture}[node distance=2.0cm, scale=0.6]
\node[block1] (SA) at (0,0) [draw] {$S_{A}$};
\node[block1] (IA) at (4,0) [draw] {$I_{A}$};
\node[block1] (RA) at (8,0) [draw] {$R_{A}$};
\node[block1] (ItA) at (4, -3) [draw] {$\tilde{I}_{A}$};
\node[] (empty) at (2, 0) [] {};  
\node[] (empty2) at (2, 0) [] {};

\draw[thick, ->] (SA) -- (IA) node[midway, above] {$\beta_{A}$};
\draw[thick, ->] (IA) -- (RA) node[midway, above] {$\dfrac{1}{\tau_I}$};
\draw[thick, ->] (RA) -- ++(0,1.75) -| (SA) node[pos=0.25, above] {$\dfrac{1}{\tau_R}$};
\draw[blue, thick, dashed, ->] (IA.270) -- (ItA.90) node[midway, right] {};
\draw[blue, thick, dashed, ->] (ItA) --(empty) node [midway, left] {$m_A$};
%\end{tikzpicture}

%\begin{tikzpicture}[node distance=2.0cm, scale=0.6]
\node[block2] (SB) at (0,-8) [draw] {$S_{B}$};
\node[block2] (IB) at (4,-8) [draw] {$I_{B}$};
\node[block2] (RB) at (8,-8) [draw] {$R_{B}$};
\node[block2] (ItB) at (4, -5) [draw] {$\tilde{I}_{B}$};
\node[] (empty3) at (2, -8) [] {};

\draw[thick, ->] (SB) -- (IB) node[midway, below] {$\beta_{B}$};
\draw[thick, ->] (IB) -- (RB) node[midway, below] {$\dfrac{1}{\tau_I}$};
\draw[thick, ->] (RB) -- ++(0,-1.75) -| (SB) node[pos=0.25, below] {$\dfrac{1}{\tau_R}$};
\draw[red, thick, dashed, ->] (IB.90) -- (ItB.270) node[midway, right] {};
\draw[red, thick, dashed, ->] (ItB) --(empty3) node [midway, left] {$m_B$};
\draw[teal, thick, dashed, ->] (ItA) to[out=225, in=150] node[midway, left] {$m_A$} (empty3);
\draw[teal, thick, dashed, ->] (ItB) to[out=150, in=225] node[midway, left] {$m_B$} (empty);

\end{tikzpicture}
    \caption{SIRS model with waning immunity and behavioral response for disease $A$ and disease $B$. The system of differential equations for this flow-diagram are represented in System \ref{eqn:SIRS} with functions for $\beta_i$ in Table \ref{table:beta-v2} and $m_i$ in Equation \eqref{eq:e}. Solid arrows represent flow of individuals. Dashed lines represent flow of information. }
    \label{fig:Eqn1FlowDiagram}
\end{figure}

We incorporate the effect of risk response on transmission such that as prevalence increases (decreases) transmission  declines (rises), representing human response to risk. We model this behavioral response, $m_i$, as
\begin{align}
\label{eq:e}
& m_i=\exp \left(-k \widetilde{I}_i\right), %\label{eq:ei_orig}
% & \frac{d \widetilde{I}_i}{d t}=\frac{I_i-\widetilde{I}_i}{\tau_p}, \label{eq:tildeI}
\end{align}
which represents the effect of response (e.g., decline in contact rate) on the transmission rate of $\beta_i$. Here, $m_i$ is formulated such that increases in prevalence ($I$), tracked through a lagged compartment of the infectious population, $\widetilde{I}_i$, lead to a decline in $m_i$. As prevalence $I_i$  increases, perceived risk of  $\widetilde{I}_i$  (modeled by an exponential delay of prevalence) also increases. Higher perceived risk leads to stronger reaction, decreasing $m_{i}$, and consequently decreasing $\beta_{i}$.

\subsection{Transmission rate}
We consider three scenarios that differ with respect to the level of behavioral spillover. Mathematically, this leads to differences in the multipliers $m_A$ and $m_B$ in the formulation of the transmission rate, denoted by $\beta_i(\cdot)$, in the SIRS model given in System \eqref{eqn:SIRS}. The specific system of equations for each scenario, with the relevant formulation of transmission, is included in Appendix \ref{sec:app:eq} for the reader's reference.

\subsubsection{Two independent diseases}

\paragraph{Scenario 1: No spillover.}
For the case of two pathogens that are not behaviorally coupled and with no cross-immunity, we assume that the transmission rate of each pathogen is formulated independently as
\begin{align*}
& \beta_i=m_i \beta_{0, i}, %\label{eq:betai_orig}\
\end{align*}
where $\beta_{0,i}$ is the infectivity of disease $i$ and $m_i$ the multiplier of disease $i$. In this case, disease $A$ and disease $B$ are spreading independently.

\subsubsection{Two behaviorally-coupled diseases} 
Here, we retain the assumption of no cross-immunity between the diseases, but include behavioral interdependencies, where NPIs implemented to contain one disease can influence the spread of the other one. Specifically, we connect the two diseases $A$ and $B$ via their behavioral mechanism. We consider two scenarios: first, perfect spillover, where the impact of behavioral responses to one disease perfectly affects the spread of the other one, possibly due to very similar transmission mechanisms; and second, a more general context where the impact of behavioral response targeting one disease is only partially effective on the other one.  

\paragraph{Scenario 2: Perfect spillover.}
As a first scenario, we assume the simple case of perfect spillover. This represents a situation where a proportional decline in transmission of disease $A$, due to behavioral response, results in the same proportional decline in the transmission of disease $B$, and vice versa. This happens when the two diseases transmit in an identical manner. Behavioral response to one can alter the spread of the other, such that transmission of each disease is affected by both $m_A$ and $m_B$. Thus, we modify the transmission rate of disease $i$, $\beta_i$, to be   
\begin{align*}
\beta_i&=\beta_{0, i} m_A m_B  ,\\%\label{eq:betai_spillover}
&=\beta_{0, i} \exp \left(-k\left(\widetilde{I}_A+\widetilde{I}_B\right)\right).%\label{eq:betai_full}
\end{align*}

\paragraph{Scenario 3: Imperfect spillover.}
For a more general description, we consider imperfect spillover, where a decline in transmission of one disease results in an $s$ fractional decline in transmission of the other disease. For example, we express $\beta_A$ as  
\begin{align*}
\beta_A&=m_A \left(1-s\left(1-m_B\right)\right) \beta_{0, A}. %\label{eq:betaA_spillover}
\end{align*}
Formally, the effect of risk response to the prevalence of $B$ on the spread of disease $A$ is $\left(1-s\left(1-m_B\right)\right)$. We formulate $\beta_B$ in a similar manner (see Table \ref{table:beta-v2}).

It is important to note that this formulation encapsulates both no spillover (independent diseases) and perfect spillover. For perfect spillover ($s=1$), the overall effect of the response to $B$ on the spread of disease $A$ is $m_B$, and for no spillover ($s=0$), the overall effect of response to disease $B$ on the spread of $A$ is constant at one. Table \ref{table:beta-v2} summarizes all the formulations for $\beta_i(\cdot)$ with $i\in\{A,B\}$. We note that in model formulation and analytical results, we present no spillover ($s=0$) and perfect spillover ($s=1$) before discussing imperfect spillover ($0<s<1$), while for numerical results, we present results ordered according to increasing spillover.

\begin{table}[ht!]
\centering
\caption{Formulations of transmission rate, $\beta_i(\cdot)$ for two diseases, $A$ and $B$, under no spillover, perfect spillover and imperfect spillover.} 
\label{table:beta-v2}
\begin{tabular}{|c|l|l|}
\hline
& \textbf{Transmission rate of disease $A$}  & \textbf{Transmission rate of disease $B$} \\
\hline
\multirow{3}{*}{\shortstack[c]{\textbf{Independent diseases}\\\textbf{No spillover ($s=0$)}}} & 
\multirow{3}{*}{$\begin{aligned}
\beta_A&=m_A\beta_{0,A}\\
&=e^{-k\widetilde{I}_A}\beta_{0,A}
\end{aligned}$} &
\multirow{3}{*}{$\begin{aligned}
\beta_B&=m_B\beta_{0,B}\\
&=e^{-k\widetilde{I}_B}\beta_{0,B}
\end{aligned}$} \\ & & \\ & & \\
\hline
\multirow{3}{*}{\shortstack[c]{\textbf{Behaviorally-coupled diseases} \\
\textbf{Perfect spillover ($s=1$)}}} & 
\multirow{3}{*}{$\begin{aligned}
\beta_A&=m_A m_B\beta_{0,A}\\
&=e^{-k\widetilde{I}_A}e^{-k\widetilde{I}_B}\beta_{0,A}
\end{aligned}$} &
\multirow{3}{*}{$\begin{aligned}
\beta_B&=m_Am_B\beta_{0,B}\\
&=e^{-k\widetilde{I}_A}e^{-k\widetilde{I}_B}\beta_{0,B}
\end{aligned}$} \\
& & \\
& & \\
%$\beta_B=m_Am_B\beta_{0,B}$   \\
\hline
\multirow{3}{*}{\shortstack[c]{\textbf{Behaviorally-coupled diseases} \\
\textbf{Imperfect spillover ($0<s<1$)}}} &
\multirow{3}{*}{$\begin{aligned}
\beta_A&=m_A \big(1-s(1-m_B)\big)\beta_{0,A}\\
&=e^{-k\widetilde{I}_A}\big(1-s(1-e^{-k\widetilde{I}_B})\big)\beta_{0,B}
\end{aligned}$} &
\multirow{3}{*}{$\begin{aligned}
\beta_B &=\big(1-s(1-m_A)\big)m_B\beta_{0,A}\\
&=\big(1-s(1-e^{-k\widetilde{I}_A})\big)e^{-k\widetilde{I}_B}\beta_{0,B}
\end{aligned}$} \\
& & \\
& & \\
\hline
\end{tabular}
\end{table}
\section{Results}

We begin by numerically simulating the dynamics of our system with a range of spillover -- none, perfect and imperfect (Section \ref{sec:results:numerical}). We choose parameters from the literature on COVID-19 and influenza to highlight our application of this model to concrete scenarios. The parameter values used for analysis and simulation of the two diseases ($A$ and $B$) are given in Table \ref{table:params}.
Following our numerical examination, we determine the existence and stability of equilibria of the model under each spillover scenario (Section \ref{sec:results:equilibria}) and assess the structural and practical identifiability of the model (Section \ref{sec:results:identifiability}).

\subsection{Numerical results}
\label{sec:results:numerical}

\begin{figure}[h]
\centering

\begin{subfigure}[t]{.32\textwidth}
\centering    \includegraphics[width=\linewidth]{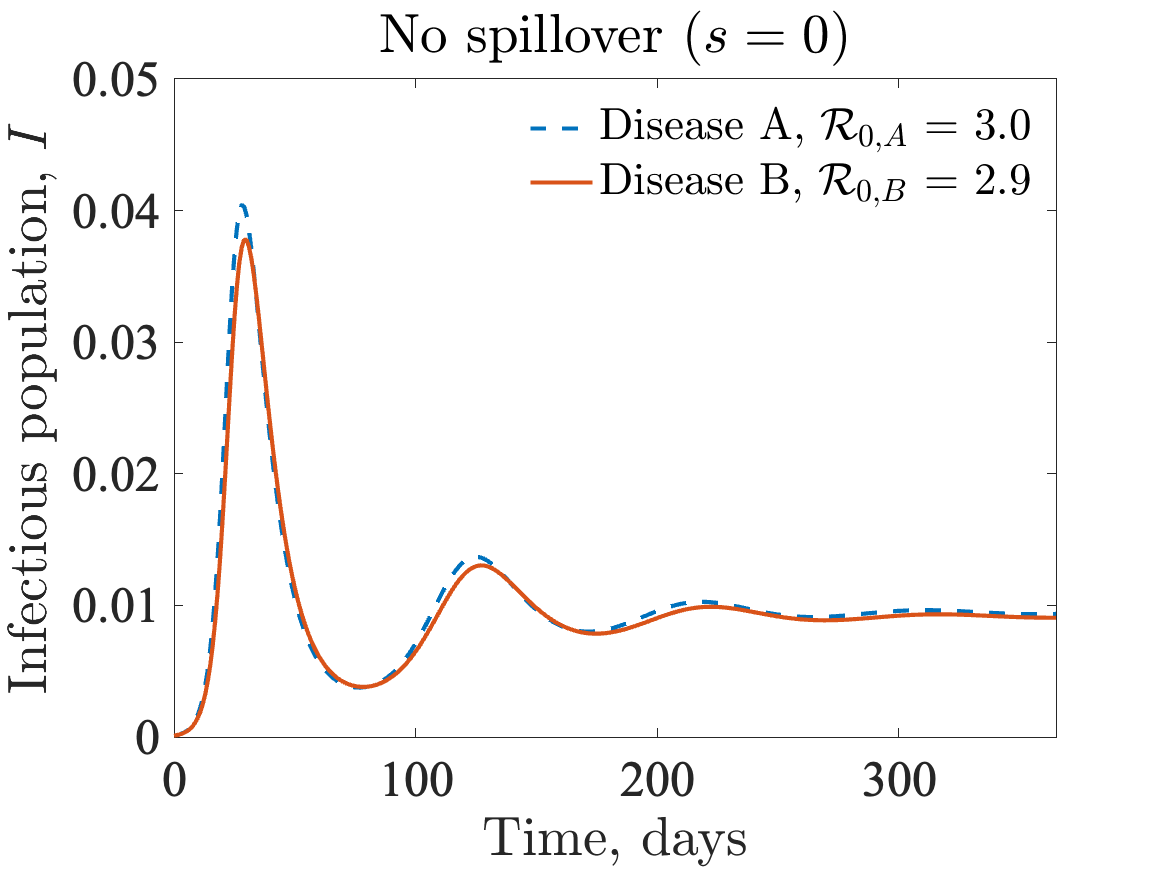}
\caption{}
\label{fig:s0b29}
\end{subfigure} 
\begin{subfigure}[t]{.32\textwidth}
\centering  
\includegraphics[width=\linewidth]{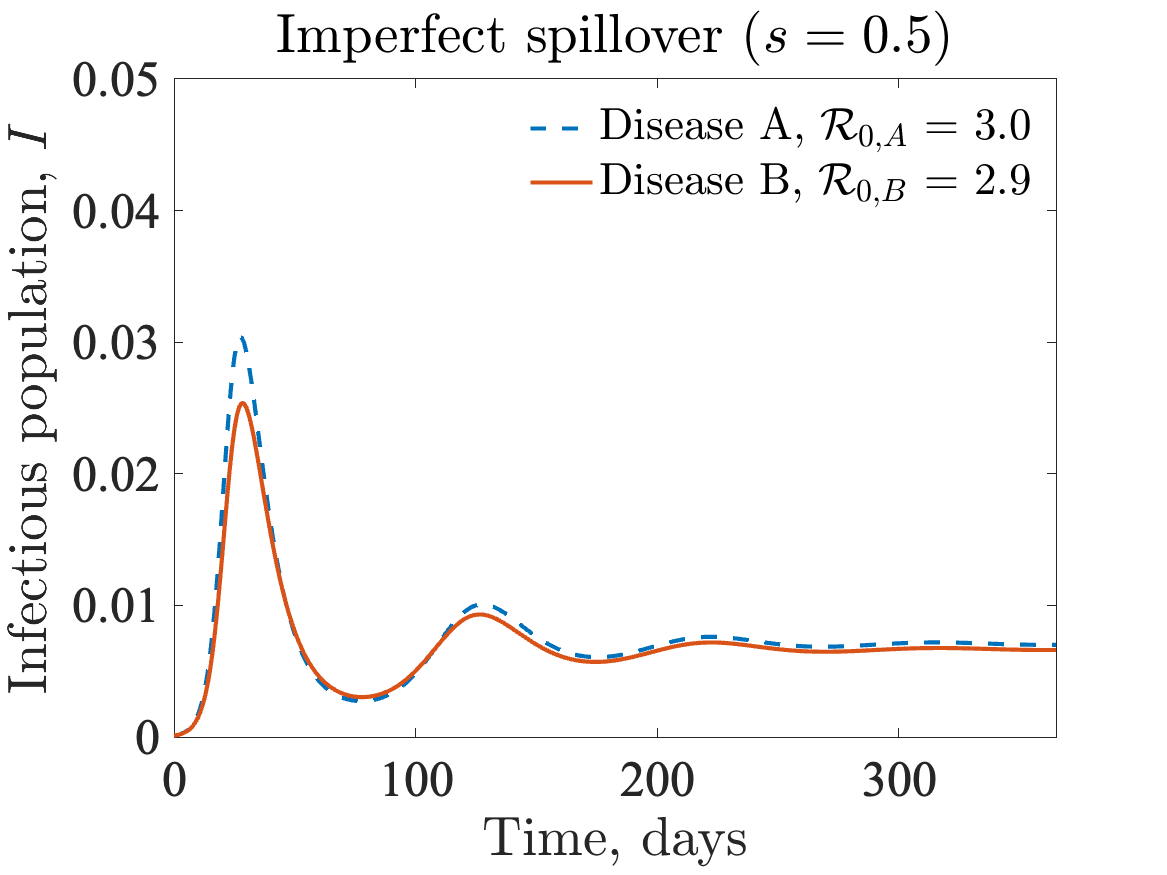}
\caption{}
\label{fig:s50b29}
\end{subfigure}
\begin{subfigure}[t]{.32\textwidth}
\centering \includegraphics[width=\linewidth]{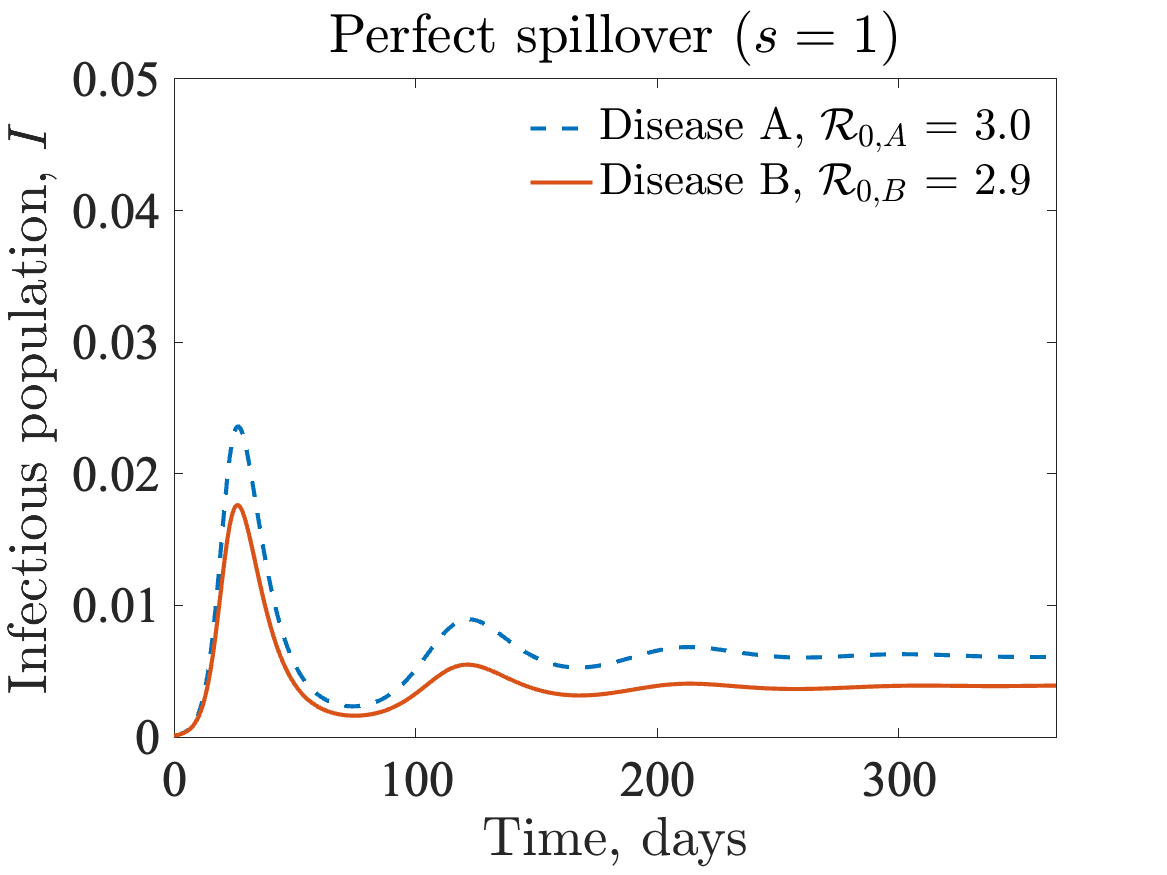}
\caption{}
\label{fig:s100b29}
\end{subfigure}\\

\begin{subfigure}[t]{.32\textwidth}
\centering \includegraphics[width=\linewidth]  {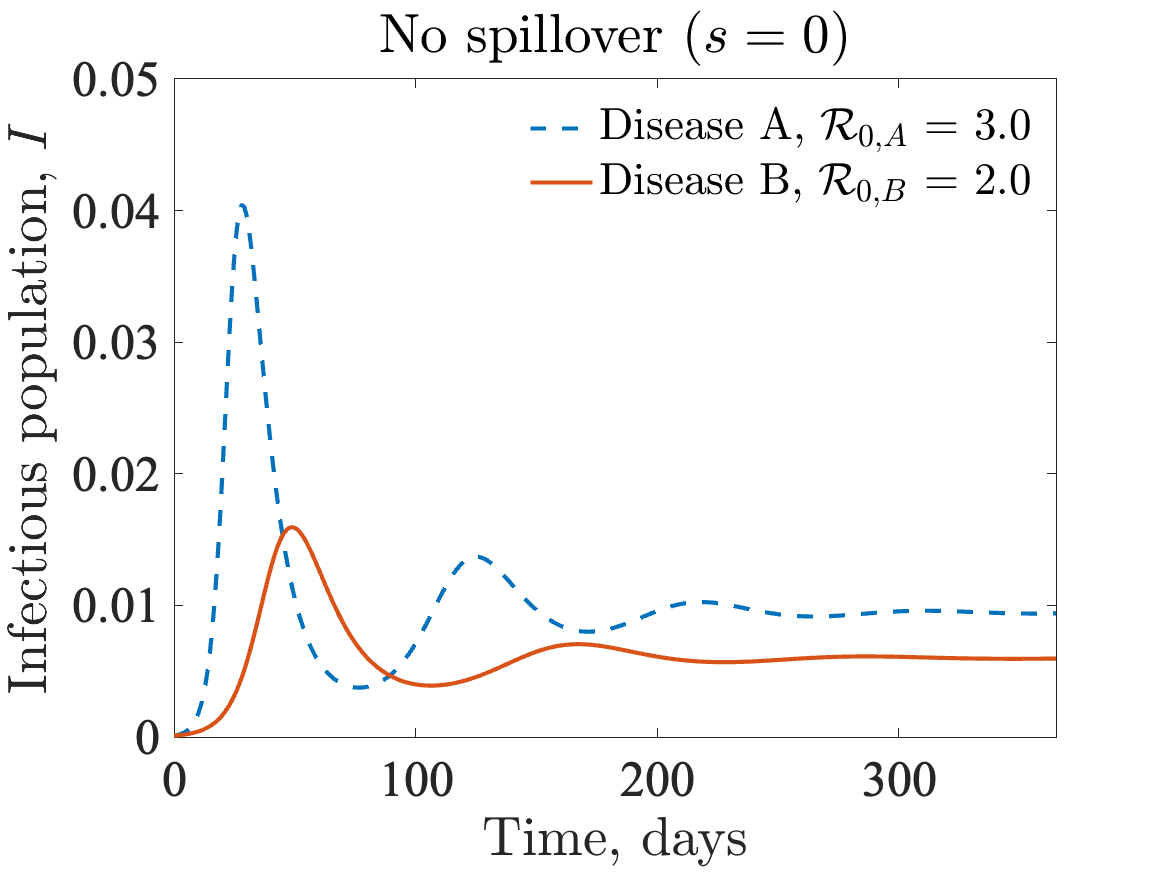}
\caption{}
\label{fig:s0b20}
\end{subfigure}
\begin{subfigure}[t]{.32\textwidth}
\centering \includegraphics[width=\linewidth]{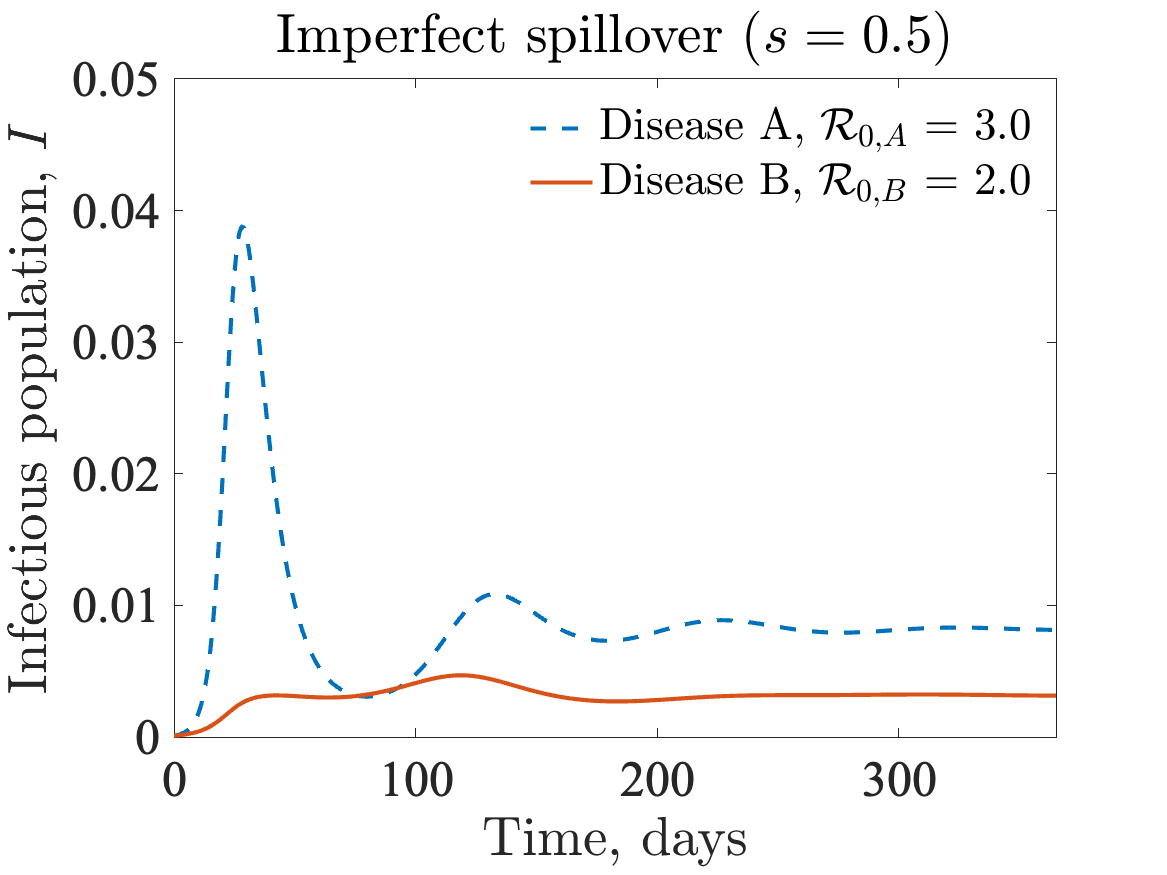}
\caption{}
\label{fig:s50b20}
\end{subfigure}
\begin{subfigure}[t]{.32\textwidth}
\centering \includegraphics[width=\linewidth]{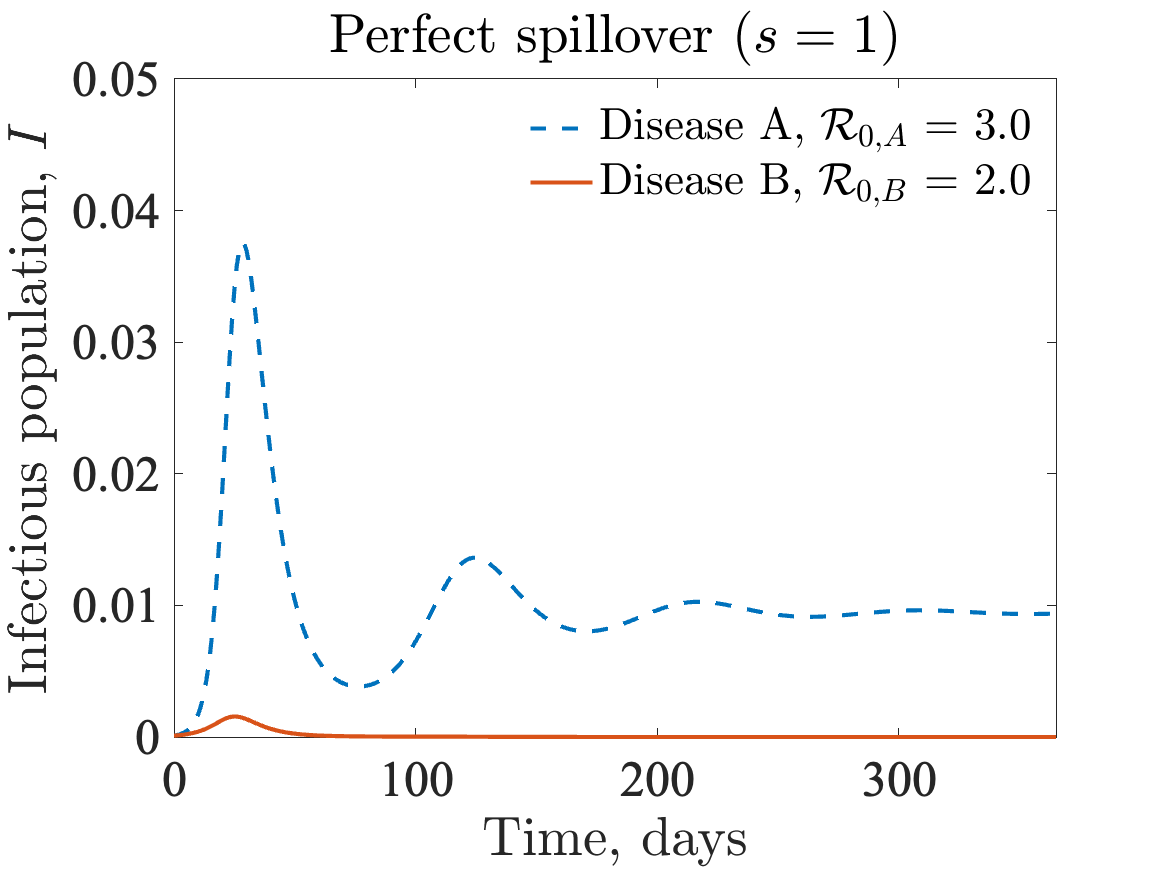}
\caption{}
\label{fig:s100b20}
\end{subfigure}\\

\begin{subfigure}[t]{.32\textwidth}
\centering \includegraphics[width=\linewidth]{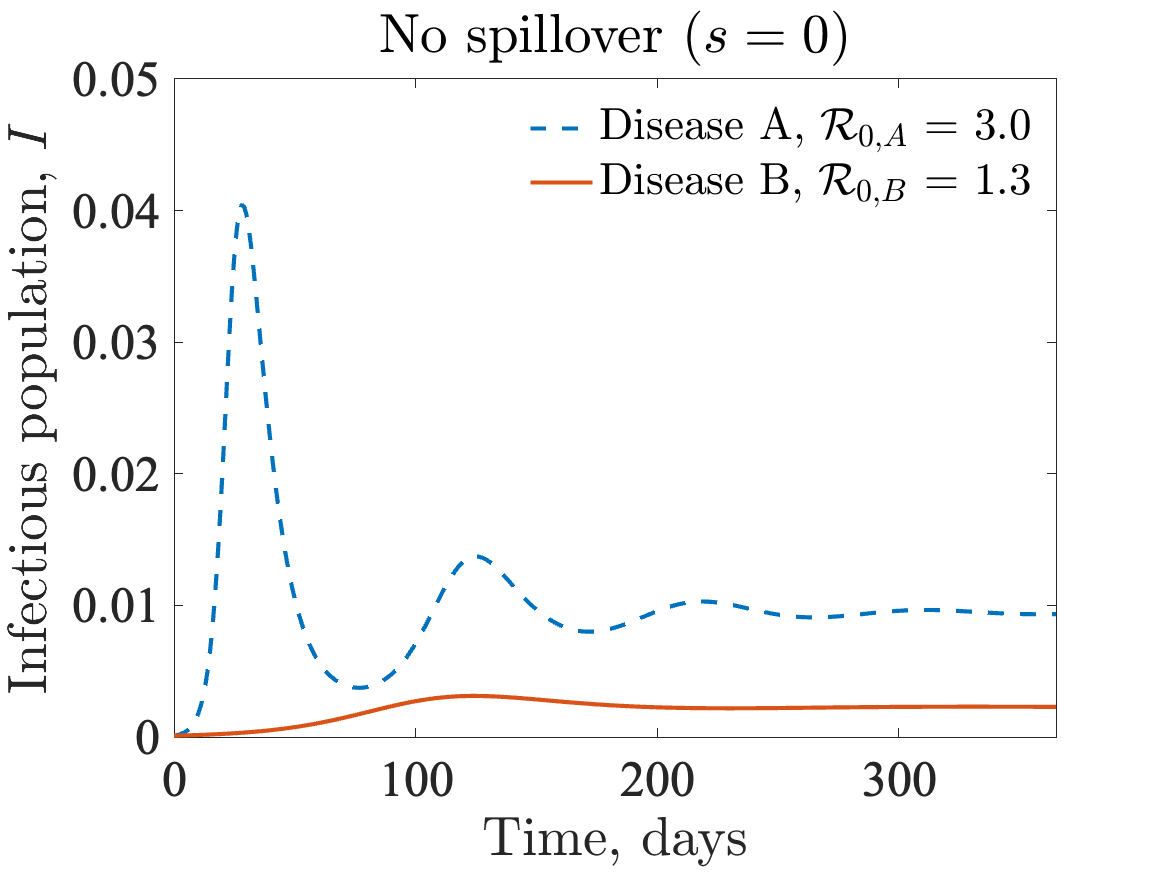}
\caption{}
\label{fig:s0b13}
\end{subfigure}
\begin{subfigure}[t]{.32\textwidth}
\centering \includegraphics[width=\linewidth]{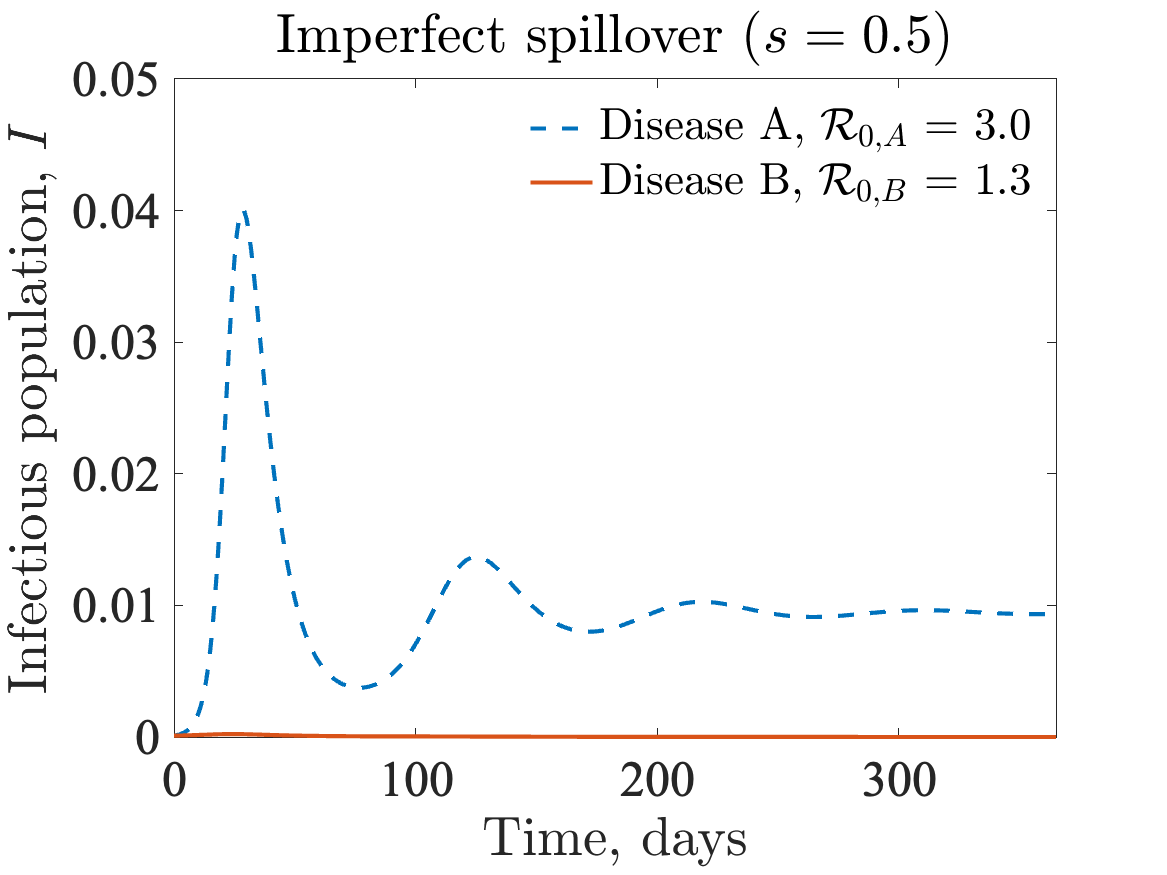}
\caption{}
\label{fig:s50b13}
\end{subfigure}
\begin{subfigure}[t]{.32\textwidth}
\centering \includegraphics[width=\linewidth]{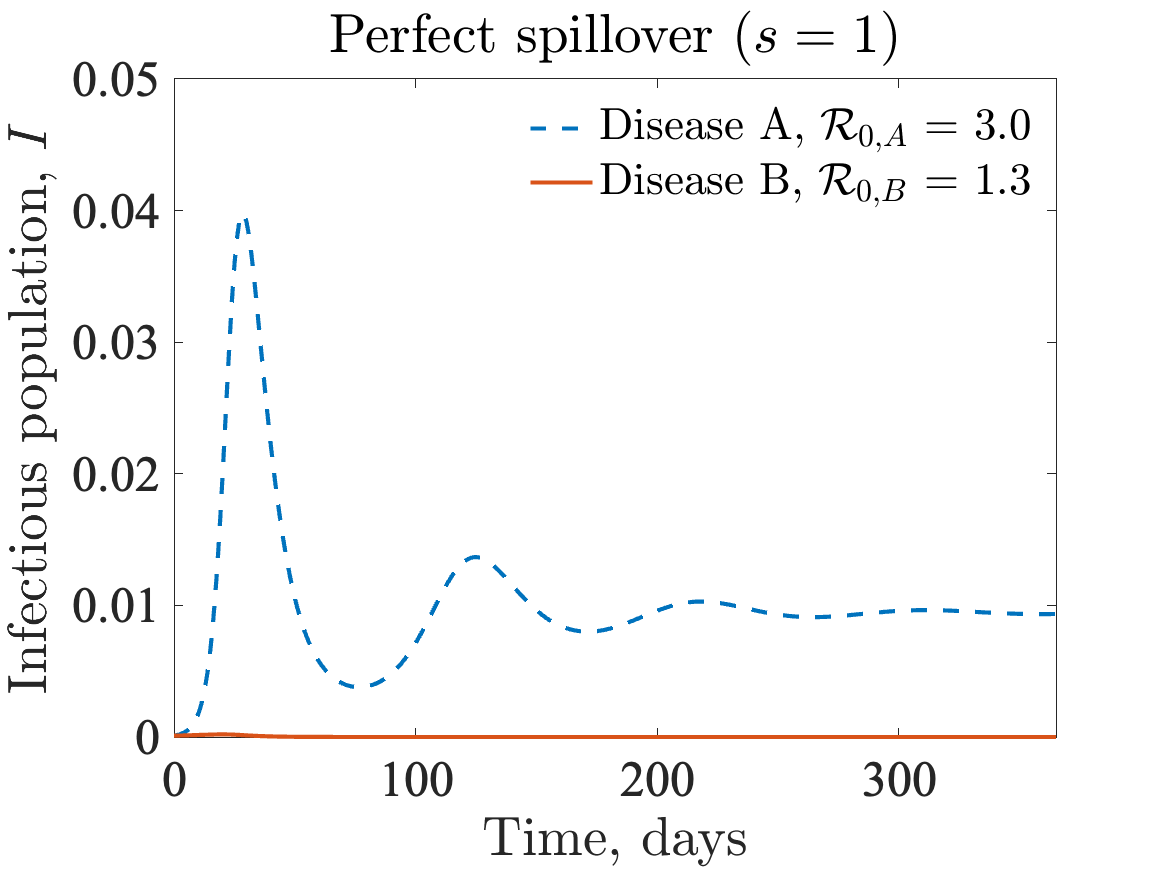}
\caption{}
\label{fig:s100b13}
\end{subfigure}\\

\caption{Dynamics with three levels of spillover and three values of the basic reproduction number of disease $B$. For (a), (d) and (g), there is no spillover, $s=0$, which corresponds to the two diseases spreading independently (i.e., the curve is identical to that of the single disease spreading in the population); for (b), (e), (h), there is  imperfect spillover, $s=0.5$; and for (c), (f), (i), there is perfect spillover, $s=1$. For (a)-(c), $\mathcal{R}_{0,B}=2.9$; for (d)-(f), $\mathcal{R}_{0,B}=2$; and for (g)-(i) $\mathcal{R}_{0,B}=1.3$). Other parameters as in Table \ref{table:params} with $\mathcal{R}_{0,A}=3$.}
    \label{fig:dyanmics_spillover}
\end{figure}

To examine the dynamics of our behaviorally-coupled system, we numerically simulate our system in MATLAB R2021a \cite{MATLAB21a} using {\tt ode45} and observe the dynamics over short (one year) and long (multiple years) time scales. Figure \ref{fig:dyanmics_spillover} shows results from nine simulation experiments which vary the level of spillover and the basic reproduction number of disease $B$ (with the basic reproduction number of disease $A$ fixed at 3).  A helpful quantity for explaining the behavior in the system is the effective reproduction number, $\mathcal{R}_e$, which tells if the disease is expanding or contracting at the current time.
For our two disease system, we define the effective reproduction number of disease $i$, $\mathcal{R}_{e,i}$, as
\begin{align*}
    \mathcal{R}_{e,i} &= \beta_i S_i \tau_I,\\
    &= m_i (1-s(1-m_j)) \beta_{0,i} S_i \tau_I,\\
    &= m_i (1-s(1-m_j)) S_i R_{0,i},\\
\end{align*}
where $m_j$ denotes dependence on disease $j$ ($j\neq i$). Without loss of generality, we consider the scenario where $\mathcal{R}_{0,A}>\mathcal{R}_{0,B}>1$. Here, the spillover level, $s$, falls between zero and one, where zero indicates no spillover (independent diseases) and one indicates perfect spillover.

\begin{table}[ht!]
\centering
\caption{Symbol, description and values for parameter and initial conditions. The values are chosen to represent spillover between two similar diseases (e.g., COVID-19 and Influenza).} 
\label{table:params}
\begin{tabular}{|c|l|c|c|c|}
\hline
\textbf{Symbol} & \textbf{Description} & \textbf{Value} & Source  \\
\hline
$\mathcal{R}_{0, A}$ & basic reproduction number of disease $A$ & 3 & \cite{d2020assessment}\\
$\mathcal{R}_{0, B}$ & basic reproduction number of disease $B$  & [1,3] & \cite{nikbakht2019comparison}\\
$\beta_{0,A}$ & infectivity  rate  of disease $A$ & $\frac{\mathcal{R}_{0, A}}{\tau_I}$ & \cite{d2020assessment, CDC_FluView_YYYY}\\
$\beta_{0,B}$ & infectivity rate of disease $B$ & $\frac{\mathcal{R}_{0, B}}{\tau_I}$ & \cite{nikbakht2019comparison, CDC_FluView_YYYY}\\
$\tau_I$ & infection period  & 7 days  & \cite{CDC_FluView_YYYY} \\
$\tau_P$ & delay to adjust risk perception  & 30 days & \cite{rahmandad2022enhancing} \\
$\tau_R$ & immunity period  & 100 days & \cite{stephens2020covid}\\
$k$ & sensitivity to risk  & 100 & \cite{ghaffarzadegan2021simulation} \\
$s$ & spillover constant & 0 & varied in [0,1] \\
\hline
$S_i(0)$ & initial proportion of disease-$i$ susceptible individuals & 0.9999 & assumed \\
$I_i(0)$ & initial proportion of disease-$i$ infectious individuals & 0.0001 & assumed \\
$R_i(0)$ & initial proportion of disease-$i$ recovered individuals &  0 & assumed \\
$\widetilde{I}_i(0)$ & initial proportion of disease-$i$ perceived infectious individuals &  0 & assumed \\
\hline
\end{tabular}
\end{table}

Without spillover (under our assumption $\mathcal{R}_{0,A}>\mathcal{R}_{0,B}>1$), there are two regimes of behavior: (i) diseases $A$ and $B$ co-exist but disease $A$ prevalence always remains above disease $B$ prevalence and (ii) diseases $A$ and $B$ co-exist but the diseases alternate whose prevalence is higher. As depicted in the left column of Figure \ref{fig:dyanmics_spillover} (Panels \ref{fig:s0b29}, \ref{fig:s0b20}, \ref{fig:s0b13}), both diseases persist as the diseases are spreading independently and (by assumption)  the basic reproduction number of each disease is above one. As the diseases are not interacting in any way, the endemic level of the disease is dictated solely by its basic reproduction number. Furthermore, the dynamics are qualitatively similar for the two diseases and are quantitatively identical when the diseases have the same traits.  As cases grow for disease $i$, individuals respond, leading to reductions in the transmission rate and a drop in $\mathcal{R}_{e,i}$. Ultimately, $\mathcal{R}_{e,i}$ settles to one such that cases reach a non-zero equilibrium (under our assumption  $\mathcal{R}_{0,i}>1$). As the basic reproduction number of disease $B$ is lowered, the initial peak occurs later and reaches a lower level (Figure \ref{fig:dyanmics_spillover}adg, red lines as $\mathcal{R}_{0,B}$ is lowered from 2.9 to 2 to 1.3). This is because the natural disease infectivity is lower, reducing the transmission rate even without behavioral response. Thus, for a range of differences in $\mathcal{R}_{0,i}$ for the two diseases -- even in the absence of interdependent behavioral feedback (i.e., diseases spread entirely independently) -- the diseases may alternate in prevalence levels, as seen in  Figure \ref{fig:dyanmics_spillover}ad, an example of regime (ii). When  $\mathcal{R}_{0,A}>>\mathcal{R}_{0,B}$ (Figure \ref{fig:dyanmics_spillover}g), while the disease co-exist, disease $A$ always has higher prevalence, an example of regime (i).

In the case of behavioral spillover, the dynamics become more interesting. The middle and right columns of Figure \ref{fig:dyanmics_spillover} show results from the scenarios of imperfect and perfect spillover, respectively. With our assumption of a higher basic reproduction number for disease $A$ ($\mathcal{R}_{0,A}>\mathcal{R}_{0,B}$), cases initially grow for disease $A$ faster than for disease $B$ and thus, people react to disease $A$ sooner than to disease $B$. As a result, the $\mathcal{R}_{e,A}$ falls, but depending on $\mathcal{R}_{0,B}$ and the level of spillover, this may provide the opportunity for disease $B$ to flourish. Here, there are now three regimes of behavior:  (i) diseases $A$ and $B$ co-exist and disease $A$ prevalence always remains above disease $B$ prevalence, (ii) diseases $A$ and $B$  co-exist but alternate whose prevalence is higher, and (iii) only disease $A$ persists, as disease $B$ is unable to persist despite a basic reproduction number above one.

When there is behavioral spillover, disease $A$ and disease $B$ impact each other's spread. This alters the region of co-existence of both diseases. Under our assumption of $\mathcal{R}_{0,A}>\mathcal{R}_{0,B}>1$, disease $A$ always persists but the region where disease $B$ also persists is determined by the combination of spillover ($0<s\leq 1$) and the transmission capacity of disease $B$, as determined by $\mathcal{R}_{0,B}$. The impact of disease $A$ on disease $B$ potentially leads to its extinction even when $\mathcal{R}_{0,B}>1$. For example, regime (iii) is seen in  Figure \ref{fig:dyanmics_spillover}fhi. 

\begin{figure}
    \centering
    \begin{subfigure}[t]{.49\textwidth}
\centering
    \includegraphics[width=\linewidth]{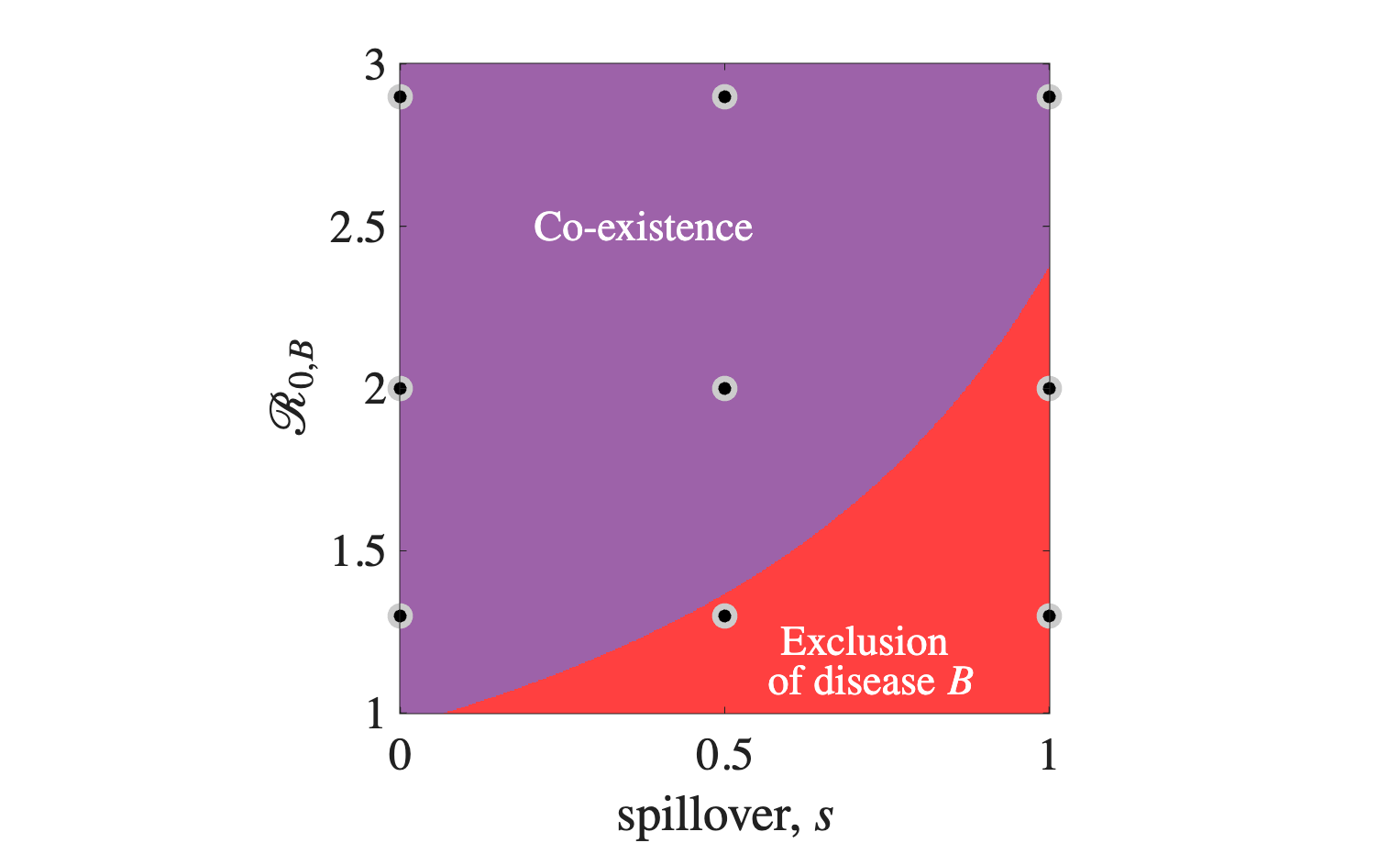}
    \caption{Co-existence vs. exclusion (simulation)}\label{fig:persistence_exceed_A}
\end{subfigure}%
    \begin{subfigure}[t]{.49\textwidth}
\centering
    \includegraphics[width=\linewidth]{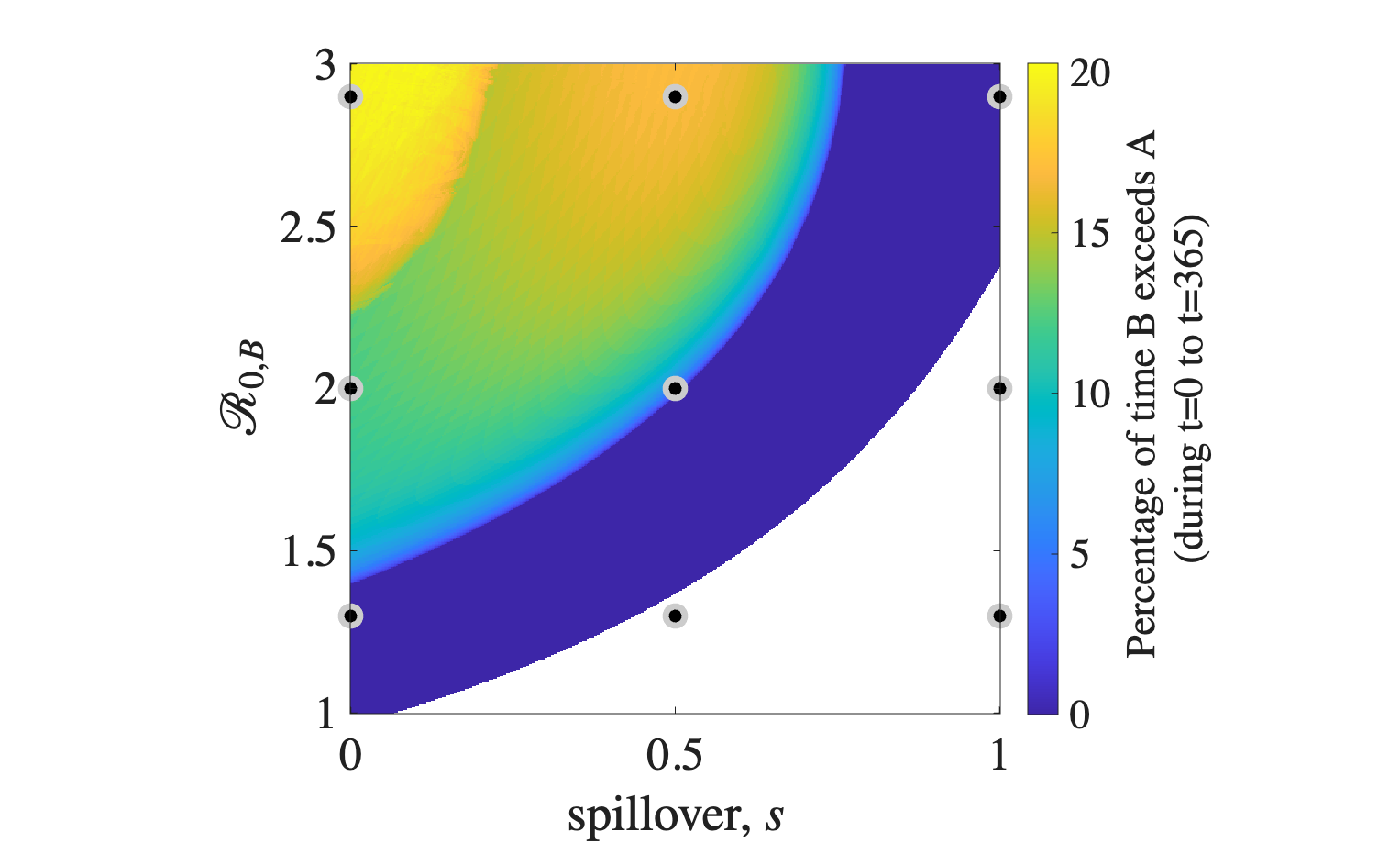}
    \caption{Time disease $B$ exceeds disease  $A$}\label{fig:persistence_exceed_B}
\end{subfigure}
    \begin{subfigure}[t]{.49\textwidth}
\centering
    \includegraphics[width=\linewidth]{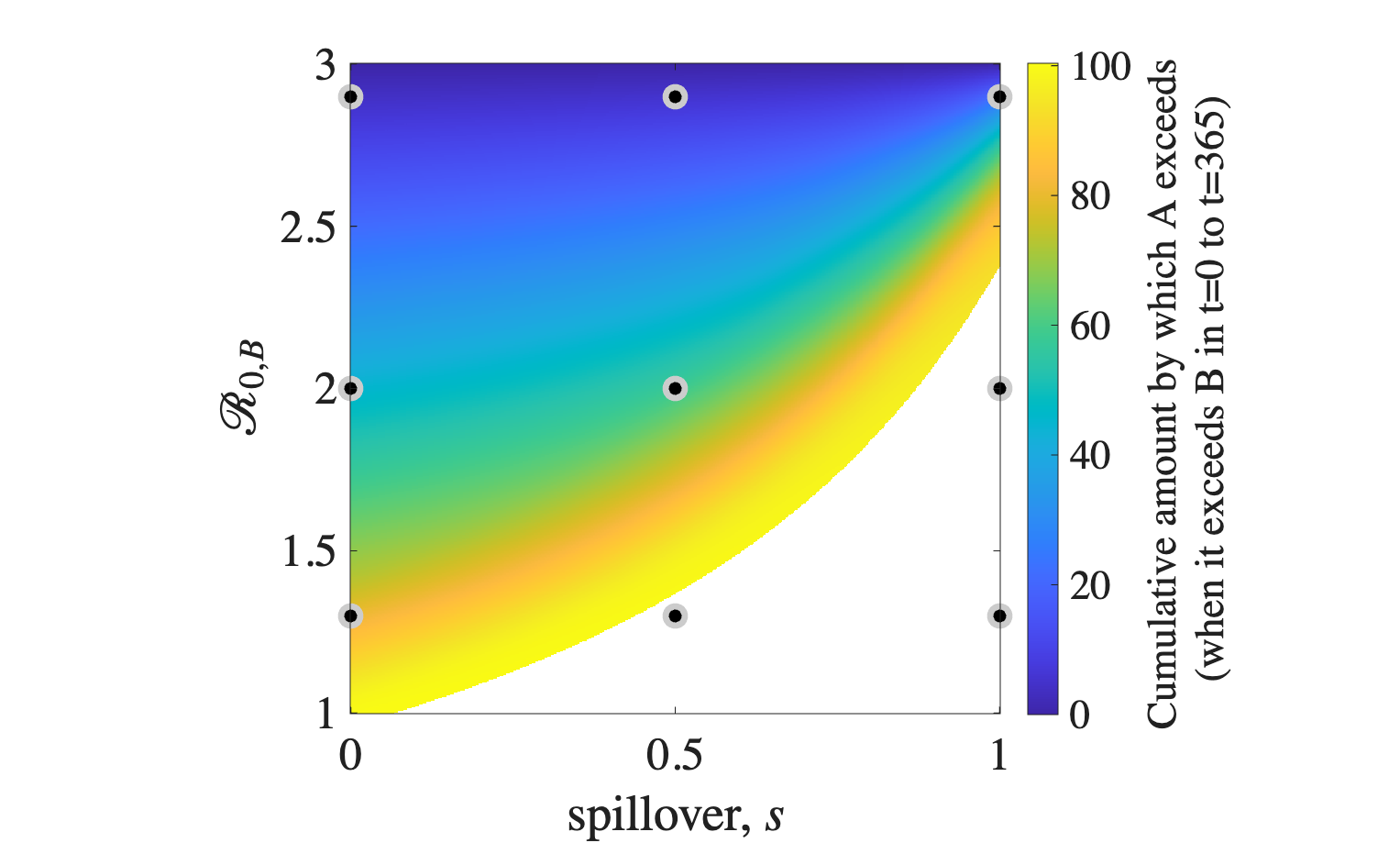}
    \caption{Amount disease $A$ exceeds disease  $B$}\label{fig:persistence_A_exceed_B}
\end{subfigure}% 
    \begin{subfigure}[t]{.49\textwidth}
\centering
    \includegraphics[width=\linewidth]{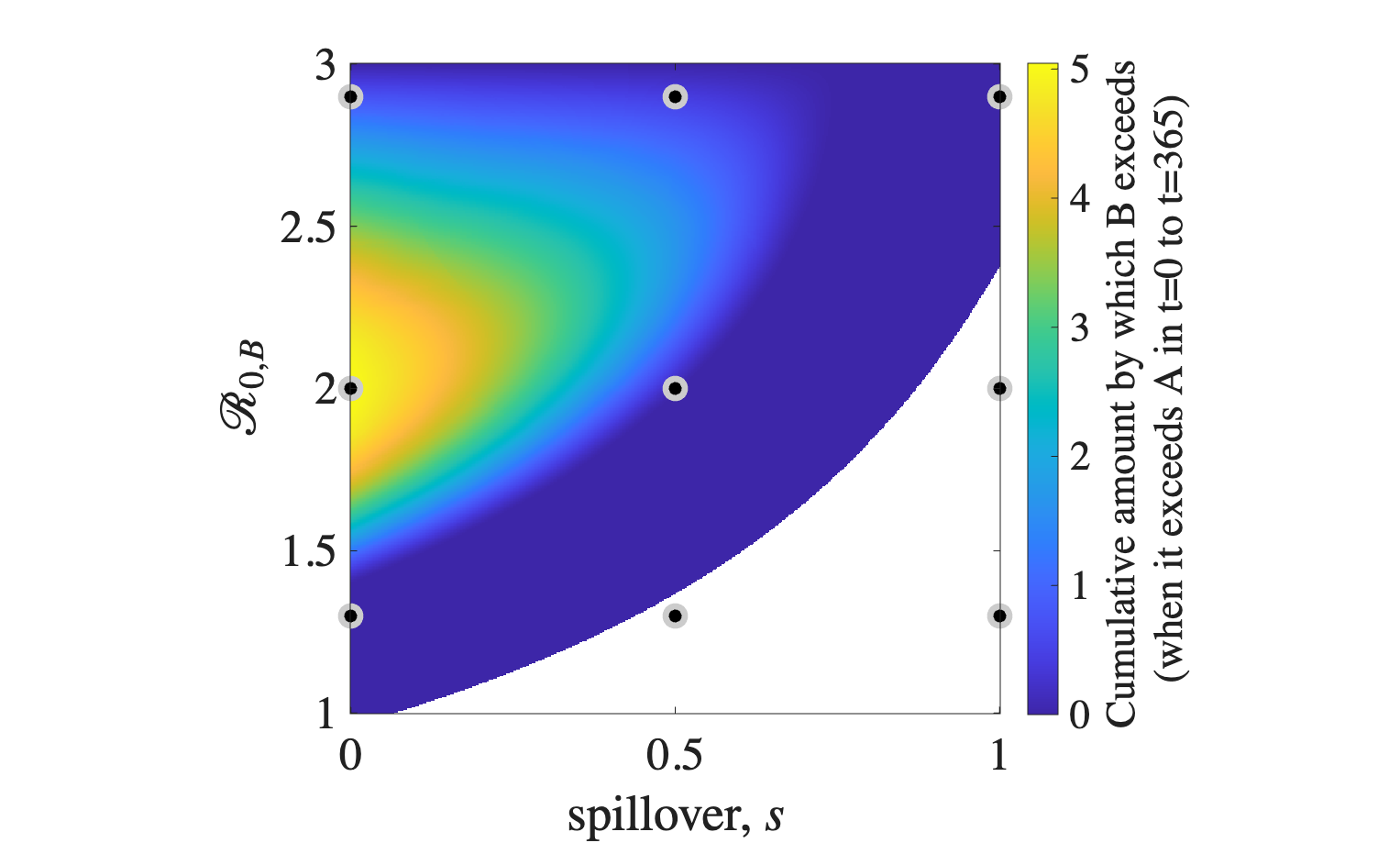}
    \caption{Amount disease $B$ exceeds disease  $A$}\label{fig:persistence_B_exceed_A}
\end{subfigure}% 
    \caption{Persistence and superiority of diseases across one year as computed through numerical simulation. (a) Persistence of both diseases (purple) or only disease $A$ (red), which has the higher basic reproduction number. (b) Percentage of the first year that disease $B$ prevalence is above disease $A$ prevalence. For $\mathcal{R}_{0,A}=\mathcal{R}_{0,B}=3$ the diseases produce independent, identical outbreaks so neither disease exceeds the other one (white line at the tope). White space at the bottom right corresponds to when there is persistence of only disease $A$.  (c) Cumulative amount during the first year that disease $A$ exceeds disease $B$. (d) Cumulative amount during the first year that disease $B$ exceeds disease $A$. Note the different color scaling in panels (b)-(d). Dots correspond to combination of spillover ($s$) and basic reproduction number of disease $B$ ($\mathcal{R}_{0,B}$) used for plots in Figure \ref{fig:dyanmics_spillover}. %Note that for $\mathcal{R}_{0,A}=\mathcal{R}_{0,B}=3$, the strains are independent and identical, regardless of the level of spillover, so neither are considered dominant.
    }
    \label{fig:persistence_exceed}
\end{figure}

We conduct a more systematic analysis, varying values of $\mathcal{R}_{0,B}$ and $s$ in small increments, and examine the final simulation results for outcome regimes.  Specifically, $\mathcal{R}_{0,B}$ is iterated between 1 and 3 (with $\mathcal{R}_{0,A}=3$) and $s$ is iterated between 0 (no spillover) and 1 (perfect spillover). Figure \ref{fig:persistence_exceed} shows simulated outcomes over the first year of an epidemic for various combinations of the two parameter values. As seen in Figure \ref{fig:persistence_exceed}a, as $\mathcal{R}_{0,B}$ increases, a larger amount of spillover (larger $s$) is required for exclusion of disease $B$. In Figure \ref{fig:persistence_exceed}b, we quantify the percentage of time in the first year where the prevalence of disease $B$ exceeds that of disease $A$.  However, when disease $B$ prevalence exceeds disease $A$, it may not be by much. Thus, in Figures \ref{fig:persistence_exceed}c and d we quantify the cumulative amount across one year by which one disease exceeds the other. The specified nine points marked in Figure \ref{fig:persistence_exceed} correspond to the scenarios presented in Figure \ref{fig:dyanmics_spillover}.
Considering the entire range of possible spillover and basic reproduction numbers of disease $B$  above one but below that of disease $A$ ($\mathcal{R}_{0,A}>\mathcal{R}_{0,B}>1$), co-existence is always possible for low spillover, but requires a large basic reproduction number for strain $B$ with high spillover (Figure \ref{fig:persistence_exceed}a). Furthermore, only for high enough $\mathcal{R}_{0,B}$ and low enough spillover can strain $B$ rise in prevalence above strain $A$ at any point (Figure \ref{fig:persistence_exceed}b) and the level by which it exceeds disease $A$ may be quite low (Figure \ref{fig:persistence_exceed}d).
Our results on the impact of behavioral spillover are robust to another choice of the behavioral feedback, as demonstrated by the use of a fractional rather than exponential formulation (see \cite{lejeune2025formulating} for discussion of various functional formulations for behavioral feedback). In particular, we reproduce Figs. \ref{fig:dyanmics_spillover} and \ref{fig:persistence_exceed} in Section \ref{sec:app:frac} with a discussion of the fractional formulation.

\subsubsection{Approximation of spillover}

\begin{figure}
    \centering
        \begin{subfigure}[t]{.49\textwidth}
\centering
    \includegraphics[width=\linewidth]{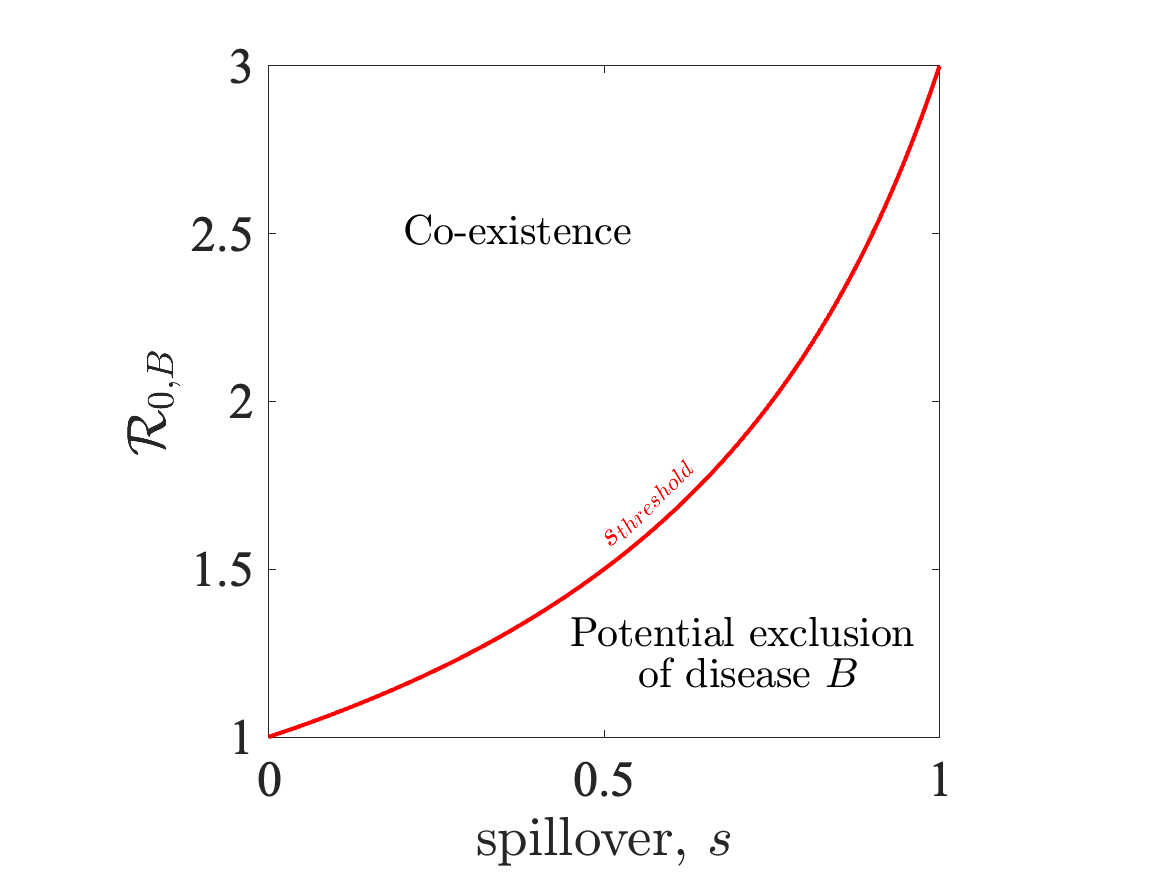}
    \caption{Co-existence vs. exclusion (approximation)}\label{fig:coexistence_A}
\end{subfigure}%
    \begin{subfigure}[t]{.49\textwidth}
\centering
    \includegraphics[width=\linewidth]{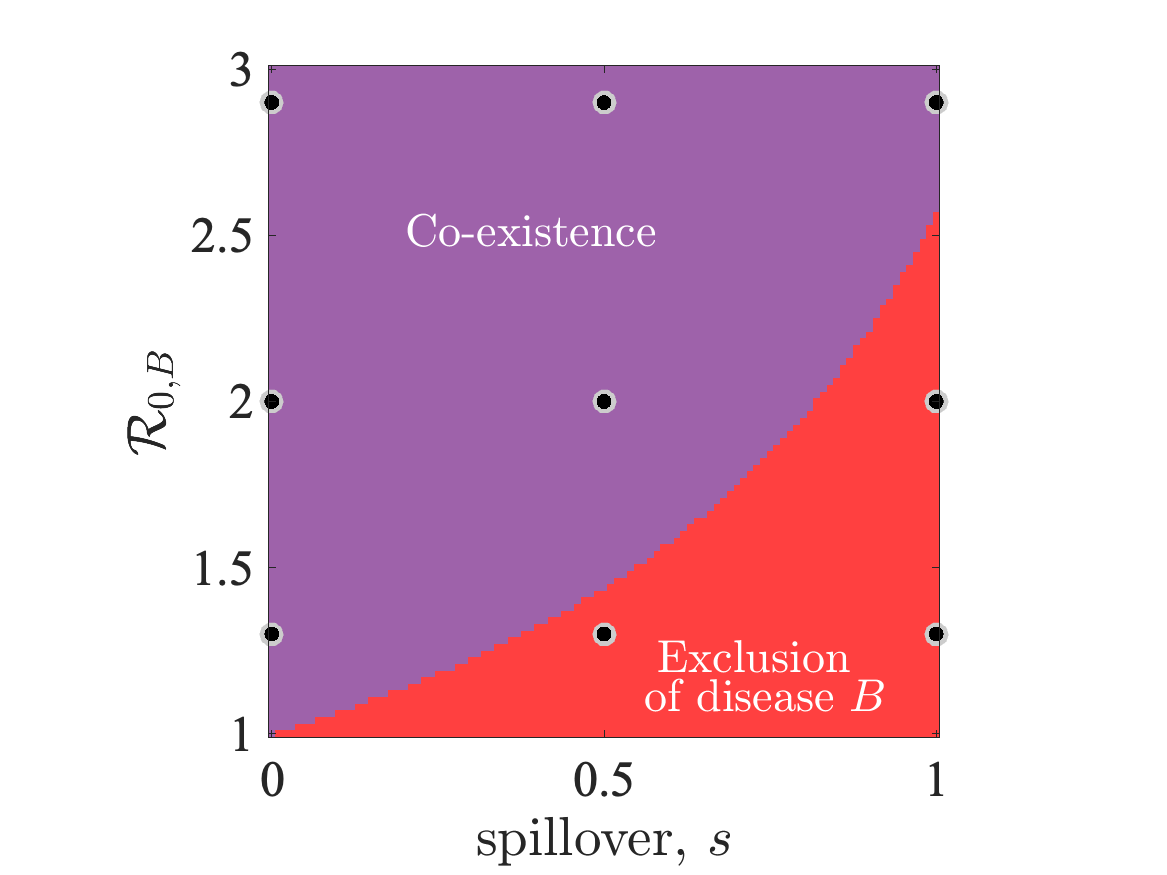}
    \caption{Co-existence vs. exclusion (analytical)}\label{fig:coexistence_B}
\end{subfigure}%  
    \caption{Analytical approximation and analytical solution of regions of disease persistence.  (a) Approximated threshold for co-existence of both diseases as determined in Theorem \ref{thm:approximation}. (b) Analytical solution (numerically computed) of persistence of both diseases (purple) or only disease $A$ (red). Recall that disease $A$ has, by assumption, the higher basic reproduction number ($\mathcal{R}_{0,A}=3$). Dots correspond to combination of spillover ($s$) and basic reproduction number of disease $B$ ($\mathcal{R}_{0,B}$) used for plots in Figure \ref{fig:dyanmics_spillover}.}
    \label{fig:coexistence}
\end{figure}

To better understand the role of spillover in co-existence of both diseases, we assess our analytical conditions of stability equilibria (see Section \ref{sec:results:equilibria} for their calculation) as well as provide an approximation of the minimum spillover required in order for disease $A$ to exclude disease $B$. See Theorem \ref{thm:approximation} below. 
The combination of the level of spillover and the basic reproduction number for disease $B$ (under our assumption that $\mathcal{R}_{0,A}>\mathcal{R}_{0,B}>1$) in order for co-existence of both diseases is seen in Figure \ref{fig:coexistence}. Compared to our numerical simulations in Figure \ref{fig:persistence_exceed}a and our numerically-evaluated analytical conditions in Figure \ref{fig:coexistence}b (which are visually nearly identical), our approximation (shown in Figure \ref{fig:coexistence}a) under-estimates the parameter regime for which there is co-existence, i.e., co-existence is still possible for larger spillover and smaller basic reproduction number than the approximation suggests. Furthermore, there are slight differences between our simulation based results (Figure \ref{fig:persistence_exceed}a) and our numerically-calculated analytical results (Figure \ref{fig:coexistence}b), such as for low $\mathcal{R}_{0,B}$ and low $s$, as we assess persistence of disease $B$ in the simulation results as falling below the threshold of $10^{-4}$ at one year. 
 
\begin{thm} \textbf{Approximated threshold of spillover for co-existence of diseases.} Without loss of generality, assuming $\mathcal{R}_{0,A}>\mathcal{R}_{0,B}>1$, the spillover fraction, $s$, needs to be above $s_{\text{threshold}}$, i.e.,
\begin{equation*}
    s_{\text{threshold}} = \frac{1-\frac{1}{\mathcal{R}_{0,B}}}{1-\frac{1}{\mathcal{R}_{0,A}}}<s<1,
\end{equation*}
for the possibility of exclusion of strain $B$.
\label{thm:approximation}    
\end{thm}

To provide an intuition about the threshold, consider ${\mathcal{R}_{0,A}}=3$ and ${\mathcal{R}_{0,B}}=1.3$ (as in the bottom row of Figure \ref{fig:dyanmics_spillover}), which may resemble the case of SARS-CoV-2 and an influenza-like virus. In this case, $s_{threshold} \approx 0.35$, such that for higher values of 
spillover ($s>0.35$), there is the possibility of complete suppression of disease $B$, the influenza-like virus. We observe this suppression when $s=0.5$ in Figure \ref{fig:dyanmics_spillover}, bottom row, middle column, where ${\mathcal{R}_{0,A}}=3$ and ${\mathcal{R}_{0,B}}=1.3$.

\subsection{Equilibria and their stability}
\label{sec:results:equilibria}

Even though epidemics exhibit complex dynamics and such systems may rarely remain at equilibrium due to external perturbations or seasonal forcing, equilibrium analysis provides insight into the underlying mechanisms of disease spread and helps guide intuition about the system’s dynamics.

As suggested by Figure \ref{fig:dyanmics_spillover}, the system exhibits multiple equilibria. In this section, we establish conditions for the existence, uniqueness, and stability of all possible equilibria for the system under the cases when $s = 0$ (no spillover), $s=1$ (perfect spillover), and $s\in(0,1)$ (imperfect spillover). Theorem \ref{thm:equilibria} establishes the existence of the equilibria for the three scenarios. Theorem \ref{thm:stabilityDFE} describes the stability of the disease-free equilibrium (DFE) and  its association to the basic reproduction numbers of each disease. Theorems \ref{thm:BE} and \ref{thm:EEstability} establish the stability of the boundary equilibria and endemic equilibrium, respectively. The proofs for each theorem are found in the appendix. We summarize the existence and uniqueness conditions as well as the stability conditions for no spillover in Table \ref{tab:equilcond_none}, for perfect spillover in Table \ref{tab:equilcond_perf}, and for imperfect spillover in Table \ref{tab:equilcond_imperf}.

\begin{thm} 
\label{thm:equilibria}
Existence of equilibria for the three scenarios

When considering two diseases ($A$ and $B$), given the biologically relevant (positive) parameters and realistic (non-negative and bounded) initial conditions, each of the three scenarios (no spillover, perfect spillover, and imperfect spillover) has four equilibria, given below, written with the compartments in the following order
\[
\left( S_A,I_A,R_A,\widetilde{I}_A,S_B,I_B,R_B,\widetilde{I}_B\right).
\]
Solutions for each of the three scenarios -- no spillover, perfect spillover, and imperfect spillover -- are denoted by $\mathcal{E}^*$, $\mathcal{\hat{E}}$, and $\mathcal{\bar{E}}$, respectively. We denote the equilibrium values for the infectious populations as $I_i = I_i^*$ for $\mathcal{E}^*$; $I_i = \hat{I}_i$ for $\mathcal{\hat{E}}$; and $I_i = \bar{I}_i$ for $\mathcal{\bar{E}}$. Here, $i\in\{A,~B,~A_E,~B_E\}$, with $I_A$, $I_B$ representing the equilibrium values where only disease $A$ or $B$ are endemic, respectively, and $I_{A_E}$, $I_{B_E}$ representing the equilibrium values where both diseases are endemic.
\begin{enumerate}
    \item The disease-free equilibrium for each scenario ($\mathcal{E}_0^*$, $\hat{\mathcal{E}}_0$, $\bar{\mathcal{E}}_0$, respectively)
    \[(1,0,0,0,1,0,0,0)\] always exists.
    \item The disease-\textit{B}-free and disease-\textit{A}-endemic equilibrium (disease-$A$ ``boundary equilibrium") for each scenario ($\mathcal{E}_A^*$, $\hat{\mathcal{E}}_A$, $\bar{\mathcal{E}}_A$, respectively)
    \[
    \left(1-I_A-\dfrac{\tau_R}{\tau_I}I_A, I_A,\dfrac{\tau_R}{\tau_I}I_A,I_A,1,0,0,0 \right)
    \]
    exists if and only if $\beta_{0,A}\tau_I>1$. Further, it follows that $\beta_{0,A}\tau_I>\exp(kI^*_A)$, $\beta_{0,A}\tau_I>\exp(k\hat{I}_A)$, and $\beta_{0,A}\tau_I>\exp(k\bar{I}_A)$ when these equilibria exist. 
    \item The disease-\textit{A}-free and disease-\textit{B}-endemic equilibrium (disease-$B$ ``boundary equilibrium") for each scenario ($\mathcal{E}_B^*$, $\hat{\mathcal{E}}_B$, $\bar{\mathcal{E}}_B$, respectively)
    \[
    \left(1,0,0,0, 1-I_B-\dfrac{\tau_R}{\tau_I}I_B, I_B,\dfrac{\tau_R}{\tau_I}I_B,I_B \right)
    \]
    exists if and only if $\beta_{0,B}\tau_I>1$. Further, it follows that $\beta_{0,B}\tau_I>\exp(kI^*_B)$, $\beta_{0,B}\tau_I>\exp(k\hat{I}_B)$, and $\beta_{0,B}\tau_I>exp(k\bar{I}_B)$ when these equilibria exist.
    \item The endemic equilibrium for each scenario ($\mathcal{E}_{A,B}^*$, $\hat{\mathcal{E}}_{A,B}$, $\bar{\mathcal{E}}_{A,B}$, respectively)
    \[
      \left(1-{I}_{A_E}-\dfrac{\tau_R}{\tau_I}{I}_{A_E}, {I}_{A_E},\dfrac{\tau_R}{\tau_I}{I}_{A_E},{I}_{A_E},1-{I}_{B_E}-\dfrac{\tau_R}{\tau_I}{I}_{B_E}, {I}_{B_E},\dfrac{\tau_R}{\tau_I}{I}_{B_E},{I}_{B_E} \right),
    \]
    exists if and only if 
    \begin{enumerate}
        \item $\beta_{0,A}\tau_I>1$ and $\beta_{0,B}\tau_I>1$ for $\mathcal{E}_{A,B}^*$;
        \item $\beta_{0,A}\tau_I>e^{k\hat{I}_{B_E}}$ and $\beta_{0,B}\tau_I>e^{k\hat{I}_{A_E}}$ for $\hat{\mathcal{E}}_{A,B}$; and
        \item $\beta_{0,A}\tau_I>\dfrac{\exp(k\bar{I}_{B_E})}{\exp(k\bar{I}_{B_E})-s(\exp(k\bar{I}_{B_E})-1)}$ and $\beta_{0,B}\tau_I>\dfrac{\exp(k\bar{I}_{A_E})}{\exp(k\bar{I}_{A_E})-s(\exp(k\bar{I}_{A_E})-1)}$ for $\bar{\mathcal{E}}_{A,B}$.
    \end{enumerate}  
\end{enumerate}
\end{thm} 

From Theorem \ref{thm:equilibria}, we observe that all disease-specific endemic boundary equilibrium values are equal (along with the endemic equilibrium values for the case of no spillover ($s=0$)), regardless of spillover; that is, $I^*_A = \hat{I}_A = \bar{I}_A = I^*_{A_E}$, and $I^*_B = \hat{I}_B = \bar{I}_B= I^*_{B_E}$.

\begin{table}[h!]
\centering
\caption{Conditions for existences and uniqueness as well as stability of equilibria with no spillover ($s=0$).}
\label{tab:equilcond_none}
\begin{tabular}{|c|c|c|c|c|c|c|}
\hline
\multirow{3}{*}{} & \multicolumn{4}{|c|}{\textbf{No spillover ($s=0$)}}\\ 
\cline{2-5} 
 & $\mathcal{E}_0^*$	& $\mathcal{E}_A^*$ & $\mathcal{E}_B^*$ & $\mathcal{E}_{A,B}^*$ \\
 \cline{2-5}
& $I_A=0$, & $I_A=I^*_A>0$, & $I_A=0$,  & $I_A=I^*_{A_E}>0$, \\
&  $I_B=0$ & $I_B=0$ &  $I_B=I^*_B>0$ &  $I_B=I^*_{B_E}>0$ \\
\hline 
\multirow{3}{*}{\shortstack[c]{Existence \& \\ uniqueness}} & \multirow{3}{*}{Always exists} & \multirow{3}{*}{$\mathcal{R}_{0,A}>1$} & \multirow{3}{*}{$\mathcal{R}_{0,B}>1$} &  \multirow{3}{*}{\shortstack[c]{~\\$\mathcal{R}_{0,A}>1$, \\
$\mathcal{R}_{0,B}>1$, 
 }}  \\ 
 & & & &\\
 & & & &\\
\hline 
\multirow{3}{*}{Stability}   & \multirow{3}{*}{\shortstack[c]{$\mathcal{R}_{0,A}<1$, \\$\mathcal{R}_{0,B}<1$}} & \multirow{3}{*}{\shortstack[c]{$\mathcal{R}_{0,A}>1$,\\$\mathcal{R}_{0,B}<1$}} &  \multirow{3}{*}{\shortstack[c]{$\mathcal{R}_{0,A}<1$, \\$\mathcal{R}_{0,B}>1$}} & \multirow{3}{*}{
\shortstack[c]{$\mathcal{R}_{0,A}>1$ \\
   $\mathcal{R}_{0,B}>1$}}\\
 &  &  &  & \\
 &  &  &  & \\
\hline
\end{tabular}
\end{table}

For Theorems \ref{thm:stabilityDFE} - \ref{thm:EEstability}, regarding stability of equilibria, we simplify our system by noting that $R_i=1-S_i-I_i$, as $R_i$ can be determined by $S_i$ and $I_i$, and  omit $R_i$ when carrying out the analysis. To this end, we make the substitution for $R_i$ in each of the three scenarios and consider only $\dfrac{dS_i}{dt}$, $\dfrac{dI_i}{dt}$, and $\dfrac{d\widetilde{I}_i}{dt}$ in each, reducing to the six-dimensional subsystem. 

\begin{thm} 
\label{thm:stabilityDFE} 
Stability of disease-free equilibrium (DFE) $\mathcal{E}_0^*$, $\hat{\mathcal{E}}_0$, $\bar{\mathcal{E}}_0$

Given biologically relevant parameters and realistic initial conditions, the disease-free equilibrium is locally asymptotically stable for all three scenarios if $\mathcal{R}_{0,A}:=\beta_{0,A}\tau_I<1$ and $\mathcal{R}_{0,B}:=\beta_{0,B}\tau_I<1$.
\end{thm}

\begin{table}[h!]
\centering
\caption{Conditions for existences and uniqueness as well as stability for equilibria with perfect spillover ($s=1$). Note that the condition for the stability of $\bar{\mathcal{E}}_{A,B}$ is a conjecture.}
\label{tab:equilcond_perf}
\begin{tabular}{|c|c|c|c|c|c|c|}
\hline
\multicolumn{1}{|c|}{} & \multicolumn{4}{|c|}{\textbf{Perfect Spillover ($s=1$})}  \\ 
\cline{2-5}
 & $\hat{\mathcal{E}}_0$	& $\hat{\mathcal{E}}_A$ & $\hat{\mathcal{E}}_B$ & $\hat{\mathcal{E}}_{A,B}$ \\
 \cline{2-5}
& $I_A=0$, & $I_A=\hat{I}_A>0$, & $I_A=0$,  & $I_A=\hat{I}_{A_E}>0$, \\
&  $I_B=0$ & $I_B=0$ &  $I_B=\hat{I}_B>0$ &  $I_B=\hat{I}_{B_E}>0$ \\
\hline 
\multirow{3}{*}{\shortstack[c]{Existence \& \\ uniqueness}} & \multirow{3}{*}{Always exists} & \multirow{3}{*}{$\mathcal{R}_{0,A}>1$} & \multirow{3}{*}{$\mathcal{R}_{0,B}>1$} &  \multirow{3}{*}{\shortstack[c]{~\\$\mathcal{R}_{0,A}>e^{k\hat{I}_{B}}$, \\
$\mathcal{R}_{0,B}>e^{k\hat{I}_{A}}$ 
}}  \\ 
 & & & &\\
 & & & &\\
\hline 
\multirow{3}{*}{Stability}   & \multirow{3}{*}{\shortstack[c]{$\mathcal{R}_{0,A}<1$, \\$\mathcal{R}_{0,B}<1$}} & \multirow{3}{*}{\shortstack[c]{$\mathcal{R}_{0,A}>1$,\\$\mathcal{R}_{0,B}<e^{k\hat{I}_{A}}$,\\ ($\implies \mathcal{R}_{0,A}>\mathcal{R}_{0,B}$)}} &  \multirow{3}{*}{\shortstack[c]{$\mathcal{R}_{0,A}<e^{k\hat{I}_{B}}$, \\$\mathcal{R}_{0,B}>1$,\\ ($\implies \mathcal{R}_{0,B}>\mathcal{R}_{0,A}$)}} & \multirow{3}{*}{\shortstack[c]{$\mathcal{R}_{0,A}>e^{k\hat{I}_{B}}$, \\$\mathcal{R}_{0,B}>e^{k\hat{I}_{A}}$}}\\
 &  &  &  & \\
 &  &  &  & \\
\hline
\end{tabular}
\end{table}

\begin{thm} 
\label{thm:BE} 
Stability of boundary equilibria $\mathcal{E}_A^*,\hat{\mathcal{E}}_A,\bar{\mathcal{E}}_A$ and $\mathcal{E}_B^*,\hat{\mathcal{E}}_B,\bar{\mathcal{E}}_B$ 

Given biologically relevant parameters and realistic initial conditions, for each scenario, we have the following conditions for stability of the boundary equilibria:
\begin{enumerate}
    \item $\mathcal{E}_A^*$ will be locally asymptotically stable if $\mathcal{R}_{0,A}>1$ and $\mathcal{R}_{0,B}<1$;
    \item $\mathcal{E}_{B}^*$ will be locally asymptotically stable if $\mathcal{R}_{0,B}>1$ and $\mathcal{R}_{0,A}<1$;
    \item $\hat{\mathcal{E}}_A$ will be locally asymptotically stable if $\mathcal{R}_{0,A}>1$ and $\mathcal{R}_{0,B}<\exp(k\hat{I}_A)$, which, along with the existence condition, implies $\mathcal{R}_{0,A}>\mathcal{R}_{0,B}$;
    \item $\hat{\mathcal{E}}_{B}$ will be locally asymptotically stable if $\mathcal{R}_{0,B}>1$ and $\mathcal{R}_{0,A}<\exp(k\hat{I}_B)$, which, along with the existence condition,  implies $\mathcal{R}_{0,B}>\mathcal{R}_{0,A}$;
    \item $\bar{\mathcal{E}}_{A}$ will be locally asymptotically stable if $\mathcal{R}_{0,A}>1$ and $\mathcal{R}_{0,B}<\dfrac{\exp(k\bar{I}_A)}{\exp(k\bar{I}_A)-s(\exp(k\bar{I}_A)-1)}$, which, along with the existence condition, implies $\mathcal{R}_{0,A}>\mathcal{R}_{0,B}$; and
    \item $\bar{\mathcal{E}}_{B}$ will be locally asymptotically stable if $\mathcal{R}_{0,B}>1$ and $\mathcal{R}_{0,A}<\dfrac{\exp(k\bar{I}_B)}{\exp(k\bar{I}_B)-s(\exp(k\bar{I}_B)-1)}$, which, along with the existence condition,  implies $\mathcal{R}_{0,B}>\mathcal{R}_{0,A}$.
\end{enumerate}
\end{thm}

Note that, in the case of spillover (either perfect or imperfect), we can have both $\mathcal{R}_{0,A},\mathcal{R}_{0,B}>1$, a condition which typically results in long-term establishment of both diseases, with the boundary equilibria still stable. In other words, the disease with the lower basic reproduction number will not be successful in persisting and will eventually become extinct, as long as its $\mathcal{R}_0$ is smaller than the respective endemic equilibrium existence condition for $0<s\leq 1$, i.e.,
\[\frac{e^{k\bar{I}_{B}}}{e^{k\bar{I}_{B}}-s(e^{k\bar{I}_{B}}-1)}>\mathcal{R}_{0,A}>1 \qquad \text{or} \qquad \frac{e^{k\bar{I}_{A}}}{e^{k\bar{I}_{A}}-s(e^{k\bar{I}_{A}})-1)}>\mathcal{R}_{0,B}>1.\]
These conditions being strictly greater than one follow from the denominator always being smaller than the numerator for any amount of spillover, $0<s\leq 1$.
For example, the boundary equilibrium with only disease $A$ for imperfect spillover, $\bar{\mathcal{E}}_A$, will be locally asymptotically stable, such that disease $B$ will return to its disease-free state after an initial outbreak, if  $\mathcal{R}_{0,B}<\dfrac{e^{k\bar{I}_{A}}}{e^{k\bar{I}_{A}}-s(e^{k\bar{I}_{A}})-1)}$,  even though both $\mathcal{R}_{0,A}>1$ and $\mathcal{R}_{0,B}>1$. Note that when this condition is met, it also holds that $\mathcal{R}_{0,A}>\mathcal{R}_{0,B}$. We see this behavior exhibited in Figure \ref{fig:dyanmics_spillover}h for $s=0.5$ with $\mathcal{R}_{0,B}=1.3$. Similar behavior occurs with perfect spillover as in Figure \ref{fig:dyanmics_spillover}f for $s=1$ with $\mathcal{R}_{0,B}=2$; and in Figure \ref{fig:dyanmics_spillover}i for $s=1$ with $\mathcal{R}_{0,B}=1.3$.

\begin{thm}
\label{thm:EEstability} 
Stability of $\mathcal{E}^*_{A,B}$

Given biologically relevant parameters and realistic initial conditions, $\mathcal{E}^*_{A,B}$ will be locally asymptotically stable if $\mathcal{R}_{0,A}>1$ and $\mathcal{R}_{0,B}>1$.
\end{thm}

Theorem \ref{thm:EEstability} only addresses the stability of the endemic equilibrium, $\mathcal{E}^*_{A,B}$, for no spillover. In addition, we present conjectures for the stability of the endemic equilibria in the case of perfect spillover (Table \ref{tab:equilcond_perf}) and imperfect spillover (Table \ref{tab:equilcond_imperf}).

\begin{table}[h!]
\centering
\caption{Conditions for existences and uniqueness as well as stability for equilibria with imperfect spillover ($0<s<1$). Note that the condition for the stability of $\bar{\mathcal{E}}_{A,B}$ is a conjecture.}
\label{tab:equilcond_imperf}
\begin{tabular}{|c|c|c|c|c|c|c|}
\hline
& \multicolumn{4}{|c|}{\textbf{Imperfect Spillover ($0<s<1$)}}  \\ 
\cline{2-5}
 & $\bar{\mathcal{E}}_0$	& $\bar{\mathcal{E}}_A$ & $\bar{\mathcal{E}}_B$ & $\bar{\mathcal{E}}_{A,B}$ \\
 \cline{2-5}
& $I_A=0$, & $I_A=\bar{I}_A>0$, & $I_A=0$,  & $I_A=\bar{I}_{A_E}>0$, \\
&  $I_B=0$ & $I_B=0$ &  $I_B=\bar{I}_B>0$ &  $I_B=\bar{I}_{B_E}>0$ \\
\hline 
  \multirow{3}{*}{\shortstack[c]{Existence \& \\ uniqueness}} & \multirow{3}{*}{Always exists} & \multirow{3}{*}{$\mathcal{R}_{0,A}>1$} & \multirow{3}{*}{$\mathcal{R}_{0,B}>1$} &  \multirow{3}{*}{\shortstack[c]{~\\$\mathcal{R}_{0,A}>\frac{e^{k\bar{I}_{B}}}{e^{k\bar{I}_{B}}-s(e^{k\bar{I}_{B}}-1)}$, \\ $\mathcal{R}_{0,B}>\frac{e^{k\bar{I}_{A}}}{e^{k\bar{I}_{A}}-s(e^{k\bar{I}_{A}})-1)}$
}}  \\ 
 & & & &\\
 & & & &\\
\hline 
 \multirow{4}{*}{Stability}   & \multirow{4}{*}{\shortstack[c]{$\mathcal{R}_{0,A}<1$, \\$\mathcal{R}_{0,B}<1$}} & \multirow{4}{*}{\shortstack[c]{$\mathcal{R}_{0,A}>1$,\\$\mathcal{R}_{0,B}<\frac{e^{k\bar{I}_{A}}}{e^{k\bar{I}_{A}}-s(e^{k\bar{I}_{A}}-1)}$,\\ ($\implies \mathcal{R}_{0,A}>\mathcal{R}_{0,B}$)}} &  \multirow{4}{*}{\shortstack[c]{$\mathcal{R}_{0,A}<\frac{e^{k\bar{I}_{B}}}{e^{k\bar{I}_{B}}-s(e^{k\bar{I}_{B}}-1)}$, \\$\mathcal{R}_{0,B}>1$,\\ ($\implies \mathcal{R}_{0,B}>\mathcal{R}_{0,A}$)}} & \multirow{4}{*}{\shortstack[c]{$\mathcal{R}_{0,A}>\frac{e^{k\bar{I}_{B}}}{e^{k\bar{I}_{B}}-s(e^{k\bar{I}_{B}}-1)}$, \\$\mathcal{R}_{0,B}>\frac{e^{k\bar{I}_{A}}}{e^{k\bar{I}_{A}}-s(e^{k\bar{I}_{A}}-1)}$}}\\
 &  &  &  & \\
 &  &  &  & \\ 
 &  &  &  & \\ 
\hline
\end{tabular}
\end{table}

\subsection{Identifiability}
\label{sec:results:identifiability}

In order to determine whether the parameter values can be uniquely estimated from the model formulation, we conduct an identifiability analysis. First, we perform structural identifiability analysis using the software STRIKE-GOLDD 4.2 to determine whether the parameters are globally identifiable, locally identifiable, or unidentifiable \cite{diaz2023controllability} (Section \ref{sec:results:structural}). For scenarios that are structurally identical, we consider a practical identifiability analysis using the Monte Carlo approach to establish which parameters are strongly identifiable, weakly identifiable, or unidentifiable (Section \ref{sec:results:practical}).

\subsubsection{Structural Identifiability}
\label{sec:results:structural}

Mathematically, we say that a model is \emph{globally structurally identifiable} for a set of parameters $\mathbf{p_{1}}$ for output $\mathbf{y}(t)$, if for every parameter set $\mathbf{p_{2}}$, the relationship $y(t, \mathbf{p_{1}})=y(t, \mathbf{p_{2}}) \hspace{0.2cm} \text{implies} \hspace{0.2cm} \mathbf{p_{1}}=\mathbf{p_{2}}.$ If there is a neighborhood in which the equation above holds true, we say the model is \emph{locally identifiable}. We consider three outputs: prevalence (level of infectious classes), incidence (entrance into the infectious classes), and `recognized prevalence' (a scaled version of perceived prevalence, which relates to the infectious level following a delay). Here we assess parameters: $\beta_{A}, \beta_{B}, k, \tau_{R}, \tau_{I}, \tau_{P}$ for prevalence ($y_{1}(t)=I_A(t)$ and $y_{2}(t)=I_B(t)$), incidence ($y_{1}(t)=\beta_{A}S_{A}(t)I_A(t)$ and $y_{2}(t)=\beta_{B}S_B(t)I_B(t) $), and recognized prevalence ($y_{1}(t)=K_A\widetilde{I}_A(t) $ and $y_{2}(t)=K_B\widetilde{I}_B(t)$, where $K_A$ and $K_B$ are a fraction of perceived prevalence).

We employ STRIKE-GOLDD 4.2 which uses Lie derivatives to obtain information about the state variables of nonlinear systems from the output measures \cite{villaverde2016structural}. Since our models contain an exponential term, we cannot use the traditional software such as DAISY or Julia as these tools are dependent on the differential algebra methodology \cite{dong2023differential,bellu2007daisy}. STRIKE-GOLDD 4.2 is implemented in MATLAB R2021a.

\begin{table}[ht!]
\centering
\caption{The identifiability results using the three formulations described in Table \ref{table:beta-v2} using the toolbox STRIKE-GOLDD 4.2. Note that {\color{red}\text{*}} means MATLAB ran out of memory to complete the simulation.} 
\label{table:IdentifabilityResults} 
\begin{tabular}{|c|l|l|}
\hline \textbf{Case}
& \textbf{Outputs}  & \textbf{Structural Identifiability Results} \\
\hline
\multirow{9}{*}{\shortstack[c]{\textbf{Independent diseases}\\\textbf{No spillover ($s=0$)}}} 
&
\multirow{3}{*}{\shortstack[c]{$\begin{aligned}
y_{1}(t)&=I_A(t)\\y_{2}(t)&=I_B(t)
\end{aligned}$}} & \multirow{3}{*}{\shortstack[l]{$\beta_{A}, \beta_{B}, k, \tau_{R}, \tau_{I}, \tau_{P}$ are \\all locally identifiable.}}\\
& & \\
& & \\
\cline{2-3}
& \multirow{3}{*}{\shortstack[c]{$\begin{aligned}
y_{1}(t)&=\beta_{A}S_{A}(t)I_A(t)
\\y_{2}(t)&=\beta_{B}S_B(t)I_B(t) \end{aligned}$}}       & \multirow{3}{*}{\shortstack[l]{$\beta_{A}, \beta_{B}, k, \tau_{R}, \tau_{I}, \tau_{P}$ are \\ all locally identifiable.}    } \\
& & \\
& & \\          
\cline{2-3}
& \multirow{3}{*}{\shortstack[c]{$\begin{aligned}
y_{1}(t)&=K_A\widetilde{I}_A(t) \\ y_{2}(t)&=K_B\widetilde{I}_B(t)   \end{aligned}$}}        & \multirow{3}{*}{\shortstack[l]{$\beta_{A}, \beta_{B}, k, K_{A}, K_{B}, \tau_{R}, \tau_{I}, \tau_{P}$ are \\all locally identifiable.}} \\  
& & \\
& & \\
\hline
\multirow{9}{*}{\shortstack[c]{\textbf{Behaviorally-coupled diseases} \\
\textbf{Perfect spillover ($s=1$)}}} & 
\multirow{3}{*}{\shortstack[c]{$\begin{aligned}
y_{1}(t)&=I_A(t)\\y_{2}(t)&=I_B(t) \end{aligned}$}}                 & \multirow{3}{*}{\shortstack[l]{$\beta_{A}, \beta_{B}, k, \tau_{R}, \tau_{I}, \tau_{P} $ are\\ all locally identifiable.\\ $\widetilde{I}_A(0)$ and $\widetilde{I}_{B}(0)$ are unidentifiable.}}   \\
& & \\
& & \\
\cline{2-3}
&
\multirow{3}{*}{\shortstack[c]{$\begin{aligned}
y_{1}(t)&=\beta_{A}S_{A}(t)I_A(t) \\ y_{2}(t&)=\beta_{B}S_B(t)I_B(t)  \end{aligned}$}}       &  \multirow{3}{*}{\color{red}{\text{*}}  }           \\
& & \\
& & \\
\cline{2-3}
& 
\multirow{3}{*}{\shortstack[c]{$\begin{aligned}
y_{1}(t)&=K_A\widetilde{I}_A(t) \\ y_{2}(t)&=K_B\widetilde{I}_B(t) \end{aligned}$}}            & \multirow{3}{*}{\shortstack[l]{$\beta_{A}, \beta_{B}, k,  K_{A}, K_{B}, \tau_{R}, \tau_{I}, \tau_{P} $ are \\all locally identifiable.}}   \\
& & \\
& & \\
%$\beta_B=m_Am_B\beta_{0,B}$   \\
\hline
\multirow{9}{*}{\shortstack[c]{\textbf{Behaviorally-coupled diseases} \\
\textbf{Imperfect spillover ($0<s<1$)}}}
&
\multirow{3}{*}{\shortstack[c]{$\begin{aligned}
y_{1}(t)&=I_A(t) \\ y_{2}(t)&=I_B(t)
\end{aligned}$}}
& \multirow{3}{*}{\shortstack[l]{$\beta_{A}, \beta_{B}, k, \tau_{R}, \tau_{I}, \tau_{P}$ are \\all locally identifiable.}} \\
& & \\
& & \\
\cline{2-3}
&
\multirow{3}{*}{\shortstack[c]{$\begin{aligned}
y_{1}(t)&=\beta_{A}S_{A}(t)I_A(t) \\ y_{2}(t)&=\beta_{B}S_B(t)I_B(t) \end{aligned}$}}          &  \multirow{3}{*}{\color{red}{\text{*}}} \\
& & \\
& & \\
\cline{2-3}
&
\multirow{3}{*}{\shortstack[c]{$\begin{aligned}
y_{1}(t)&=K_A\widetilde{I}_A(t) \\ y_{2}(t)&=K_B\widetilde{I}_B(t) \end{aligned}$}}          & \multirow{3}{*}{\shortstack[l]{$\beta_{A}, \beta_{B}, k,  K_{A}, K_{B}, \tau_{R}, \tau_{I}, \tau_{P}$ are \\all locally identifiable. }} \\
& & \\
& & \\
\hline
\end{tabular}
\end{table}

The output of STRIKE-GOLDD 4.2 indicates that all the parameters are locally identifiable for prevalence and recognized prevalence. All parameters are locally identifiable for incidence when there is no spillover. We note that when incidence was used as an output measure in the cases for perfect spillover and imperfect spillover, MATLAB ran out of memory and did not produce any results. Thus, we do use incidence as an output measure when performing practical identifiability. All of the initial conditions of the state variables are identifiable except in the case of perfect spillover when using prevalence as the output measure.

\subsubsection{Practical Identifiability}
\label{sec:results:practical}

When fitting to data, it is important to verify if the model is formulated in a way that allows for the model parameters to be uniquely identified in the context of imperfect data to determine practical identifiability.
Using prevalence data and recognized prevalence data, we considered ten scenarios to determine which parameters are strongly identifiable, weakly identifiable, or nonidentifiable. We say that a parameter is strongly identifiable if the average relative error (ARE) for all noise levels satisfies 0 $\leq$ ARE $\leq  \sigma$, where $\sigma$ is the measurement error added to the data set. A parameter is weakly identifiable if $\sigma < $ ARE $ \leq 10\times\sigma$ for all noise levels, and a parameter is practically nonidentifiable if ARE $>10\times\sigma$ for at least one noise level. The ARE values were generated through the use of the Monte Carlo (MC) approach, a sampling method that uses probability distributions to test the practical identifiability of a model (outlined in \ref{Appendix:PI}) \cite{saucedo2024comparative,miao2011identifiability}. We use the parameter values in Table \ref{table:params} and set the parameters $K_A$ and $K_B$ both equal to 0.5.

{The MC approach indicates that for prevalence as an output measures, regardless of the level of spillover, infectivity rates ($\beta_A$, $\beta_B$) and the delay to adjust risk perception ($\tau_P$) are identifiable while sensitivity to risk ($k$) and infection period ($\tau_I$) are weakly identifiable. For both $k$ and $\tau_I$, the MC approach indicates the parameters are close to being identifiable, as ARE exceeds the noise level ($\sigma$) for only a few instances (see Tables \ref{tab:Practical_Prevalence_None}, \ref{tab:Practical_Prevalence_Imperfect_0.1}, \ref{tab:Practical_Prevalence_Imperfect_0.5}, \ref{tab:Practical_Prevalence_Imperfect_0.9}, \ref{tab:Practical_Prevalence_Perfect}). The immunity period ($\tau_R$) ranges from weakly identifiable to non-identifiable depending on the level of spillover. For all the scenarios, $\tau_R$ is the only parameter in which we observe a drastic jump in the ARE value. For example, there is a 184.73\% increase of the ARE value for $\tau_R$ when moving from 5\% to 10\% noise level in Table \ref{tab:Practical_Prevalence}. As we only considered data across a single year and the immunity period is assumed to be 100 days, it is possible that identifiability results would improve for longer periods of data. 

For recognized prevalence as an output measure, regardless of the level of spillover, the infectivity of disease $A$ ($\beta_A$) and the delay to adjust risk perception ($\tau_P$) are identifiable. Additionally, the infectivity of disease $B$ ($\beta_B$), the sensitivity to risk ($k$), infection period ($\tau_I$), and the fraction of reports cases for both diseases ($K_A$, $K_B$) are weakly identifiable across all spillover values tested. For the weakly identifiable parameters, the MC approach indicates the parameters are close to being identifiable, as ARE exceeds the noise level ($\sigma$) for only a few instances (see Tables \ref{tab:Practical_Perceived_Prevalence_None}, \ref{tab:Practical_Perceived_Prevalence_Imperfect_0.1}, \ref{tab:Practical_Perceived_Prevalence_Imperfect_0.5}, \ref{tab:Practical_Perceived_Prevalence_Imperfect_0.9}, \ref{tab:Practical_Perceived_Prevalence_Perfect}).  For recognized prevalence, the immunity period ($\tau_R$) is only weakly identifiable for no spillover ($s=0$) and non-identifiable for any amount of spillover. Similar to the prevalence scenario, $\tau_R$ is the only parameter that exhibits a substantial increase in the ARE value as the noise level escalates.  In Table \ref{tab:Practical_Perceived_Prevalence_None}, the ARE value for $\tau_{R}$ increases by 361.92\% from noise level 5\% to 10\%.

\begin{table}[]
    \centering
    \caption{Practical identifiability for System \ref{eqn:SIRS} with prevalence data under various spillover values. Results obtained using MC approach. See Tables \ref{tab:Practical_Prevalence_None}, \ref{tab:Practical_Prevalence_Imperfect_0.1}, \ref{tab:Practical_Prevalence_Imperfect_0.5}, \ref{tab:Practical_Prevalence_Imperfect_0.9}, \ref{tab:Practical_Prevalence_Perfect} for average relative error (ARE) values used to determine practical identifiability.}
\label{tab:Practical_Prevalence}
    \begin{tabular}{|c|c|c|c|c|c|c|}
    \hline Spillover & $\beta_A$	& $\beta_B$ & $k$ & $\tau_R$	& $\tau_I$ & $\tau_P$	 \\
    \hline
       $s=0$  &  Yes & Yes & Weakly & Weakly & Weakly & Yes  \\
        $s=0.1$ & Yes & Yes & Weakly & Weakly & Weakly & Yes \\
        $s=0.5$ & Yes & Yes & Weakly & No & Weakly & Yes \\
        $s=0.9$ &  Yes & Yes & Weakly & No & Weakly & Yes \\
        $s=1$ &  Yes & Yes & Weakly & Weakly & Weakly & Yes\\
        \hline
    \end{tabular}    
\end{table}

\begin{table}[]
    \centering
    \caption{Practical identifiability for System \ref{eqn:SIRS} with recognized prevalence data under various spillover levels. Results obtained using MC approach. See Tables \ref{tab:Practical_Perceived_Prevalence_None}, \ref{tab:Practical_Perceived_Prevalence_Imperfect_0.1}, \ref{tab:Practical_Perceived_Prevalence_Imperfect_0.5}, \ref{tab:Practical_Perceived_Prevalence_Imperfect_0.9}, \ref{tab:Practical_Perceived_Prevalence_Perfect} for average relative error (ARE) values used to determine practical identifiability.}
    \label{tab:Practical_Perceived_Prevalence}
    \begin{tabular}{|c|c|c|c|c|c|c|c|c|}
    \hline Spillover & $\beta_A$	& $\beta_B$ & $k$ & $\tau_R$	& $\tau_I$ & $\tau_P$	 & $K_A$ & 	$K_B$ \\
    \hline
       $s=0$  &  Yes & Weakly & Weakly & Weakly & Weakly & Yes & Weakly & Weakly \\
        $s=0.1$ & Yes & Weakly & Weakly & No & Weakly & Yes & Weakly & Weakly \\
        $s=0.5$ & Yes & Weakly & Weakly & No & Weakly & Yes & Weakly & Weakly \\
        $s=0.9$ &  Yes & Weakly & Weakly & No & Weakly & Yes & Weakly & Weakly \\
        $s=1$ &  Yes & Weakly & Weakly & No & Weakly & Yes & Weakly & Weakly \\
        \hline
    \end{tabular}
\end{table}

\section{Discussion}

This paper examines the dynamics of interdependencies across multiple pathogens in the absence of cross-immunity. While previous studies have often attributed inter-pathogen interactions to cross-immunity -- where infection by one pathogen confers immunity against another \cite{bhattacharyya2015cross, gog2002dynamics} -- this study relaxes that assumption and instead highlights a behavioral mechanism that can couple the spread of two pathogens. This perspective is motivated by observations during the COVID-19 pandemic, where a surge in COVID-19 cases coincided with a marked decline in Influenza A and B cases \cite{rubin2020happens,jones2020covid}. This decline was largely due to public adherence to NPIs implemented to control COVID-19, which also inadvertently suppressed the spread of influenza \cite{solomon2020influenza}.  This observation raises an important question: under what conditions can behavioral responses to one deadly disease contain the spread of others?
We do not aim  to reproduce
the full complexity of epidemic dynamics here; rather, this work is  theoretical investigation into behaviorally mediated interactions between diseases in a framework where analysis is possible. 

To this end, we develop a mathematical model to examine behavioral interdependencies across two pathogens. Building on a previously established behavioral epidemic model (SEIRb) \cite{lejeune2024mathematical, rahmandad2022enhancing}, we extend the framework to capture the spread of two pathogens, where behavioral responses to one can influence the transmission dynamics of the other. We refer to this form of interdependence as behavioral spillover. Our model represents two diseases of different infectivity, and varying degrees of behavioral interdependence. Specifically, the model investigates the extent to which NPIs targeting one pathogen affect the spread of another.

Our analysis explores a range of simulation outcomes for different values of the basic reproduction number and spillover strength. The results indicate that NPIs targeted at a pathogen with a relatively high basic reproduction number can significantly suppress the spread of other pathogens in the short term. This leads to scenarios in which one pathogen exceeds the rest while others maintain effective reproduction numbers below one, i.e., are fully suppressed. For example, in the presence of complete behavioral spillover, we observed that the population may only experience the disease with the highest $\mathcal{R}_0$, as if the ecosystem, including the diseases and humans, selects the pathogen with the highest $\mathcal{R}_0$. Additionally, we identify parameter regimes in which both pathogens coexist but spread at reduced rates compared to their behavior in isolation. With imperfect spillover, where NPIs designed for one disease only partially affect the spread of others, we estimate the range of values $\mathcal{R}_0$ in which only one of the pathogens spreads.

Moreover, our study explored equilibria properties of the new model given in System \eqref{eqn:SIRS}. %as well as structural and practical identifiability.  
Each scenario has four possible equilibria, depending on the values of the $\mathcal{R}_{0,i}$: one disease-free, two boundary (endemic in one disease and disease-free in the other), and one endemic equilibrium. While results regarding the disease-free equilibria were typical for disease models, i.e., the basic reproductive number functioned as a threshold value at one, results for the boundary and endemic equilibria when behavioral spillover was present diverged from this norm as one was not the threshold for stability for the equilibria. Within these scenarios, it is possible to have $\mathcal{R}_{0,A},\mathcal{R}_{0,B}>1$ with a stable boundary equilibrium and unstable endemic equilibrium; that is, the disease with the greater reproduction number will persist, while the other disease will die out despite both having a basic reproduction number above one, the typical threshold for persistence \cite{blackwood2018introduction, lejeune2023effect, hernandez1999threshold}. In cases with spillover, transmission depends on both diseases; analytically, the threshold condition for the newly introduced disease thus depends on the value of the infectious equilibrium of the endemic disease. If this threshold is large enough (i.e. if the presence of the endemic disease is large enough), then the newly introduced disease will be unsuccessful in establishing itself endemically. We identify the thresholds \OS{that} determine the persistence disease(s) for all three scenarios.

As a critical component of working with a novel model is the parameterization with respect to real world data. This procedure usually involves a thorough literature review or model fitting. To assess the potential utility of our novel model to be fit to data, we conducted  a structural identifiability analysis \cite{saucedo2024comparative, tuncer2018structural}.  Our identifiability analysis demonstrated that all the parameters are structurally locally identifiable with respect to prevalence and recognized prevalence data. The only scenario in which we obtained unidentifiable results was for the initial conditions of perceived prevalence, $\widetilde{I}_A(0)$ and $\widetilde{I}_B(0)$, when using prevalence data in the case of perfect spillover. Although $\widetilde{I}_A(0)$ and $\widetilde{I}_B(0)$ are unidentifiable, the initial condition values are known, which provided confidence to proceed with the practical identifiability analysis. From the practical identifiability analysis, we observed that the majority of the parameters were identifiable in some form (i.e. strongly or weakly) with respect to prevalence and recognized prevalence. The immunity period, $\tau_{R}$, was the only parameter that was not identifiable. Since we used $\tau_{R}=100$ representing an immunity period of 100 days, and we simulated the diseases dynamics for 100 days, it is probable that $\tau_{R}$ is identifiable if the simulations were lengthened considerably or if the  immunity period was shortened.

This study contributes to the fields of mathematical biology and epidemiological modeling, particularly to the emerging area of behavioral epidemic modeling \cite{funk2015nine}. Firstly, it adds to the body of literature on multi-pathogen dynamics \cite{bhattacharyya2015cross, gog2002dynamics} by offering a contrasting perspective: that interdependent epidemic patterns can arise even in the absence of cross-immunity, purely due to behavioral interdependencies. This paper introduced the term behavioral spillover to describe this phenomenon, in which behavioral responses to one pathogen influence the transmission of others. An important implication is that ignoring such behavioral interdependencies may lead to misleading conclusions—such as over- or under-estimating pathogen transmission rate or wrongly inferring the presence of cross-immunity between pathogens. 
Second, this study contributes to \OS{the} endogenous formulation of human behavior \cite{lejeune2025formulating}. It builds on prior behavioral epidemic models that incorporate a risk-response feedback loop, wherein rising (or falling) prevalence increases (or decreases) adherence to non-pharmaceutical interventions (NPIs), which in turn influences transmission \cite{rahmandad2022missing,espinoza2024adaptive,mashayekhi2025dynamics}. The paper extends this framework to demonstrate how multiple pathogens can become dynamically coupled through shared behavioral responses.

However, this study has several limitations that point to directions for future research. First, the current model assumes a homogeneous population and does not account for heterogeneity in risk or behavioral response. For example, while elderly populations are at higher risk for both COVID-19 and Influenza A/B, other diseases like Lyme disease primarily affect different demographics, such as children. Incorporating population heterogeneity could yield more nuanced insights. Second, the model does not include vaccination dynamics, which can significantly alter behavioral interdependencies. For instance, after COVID-19 vaccines became available, increased perceived safety led to greater social activity, potentially raising the risk of influenza spread. Accounting for such vaccine-induced behavioral shifts would be a valuable extension. Finally, the model assumes infected individuals do not alter their behavior and does not include biological cross-immunity between pathogens. While these assumptions simplify the system, they allow us to isolate behavioral mechanisms, showing that apparent protection against one disease can arise purely from changes in perceived risk of the other. Relaxing these assumptions—by incorporating quarantine upon infection or partial cross-immunity—could further strengthen the interdependencies and accelerate the practical extinction of diseases with lower reproduction numbers.

In conclusion, while the current model focuses on two pathogens for clarity, it can be generalized to consider multiple specific pathogens, different variants, and combinations of cross-immunity and behavioral spillover. Such extensions could offer broader insights into epidemic interactions in complex real-world settings. Furthermore, our focus was on offering a mathematical formulation of relevant models and the effort should be extended to include empirical tests using data across different regions and over a long period of time.

%\section{Acknowledgments}
%This work was supported by the US National Science Foundation, Division of Mathematical Sciences and Division of Social and Economic Sciences, Award \# 2229819.

\section{Disclosure Statement}
The authors report there are no competing interests to declare.

\section{Funding}
This work was supported by the US National Science Foundation, Division of Mathematical Sciences and Division of Social and Economic Sciences under Grant \# 2229819.

\section{Data Availability Statement}
All code files will be posted to GitHub upon acceptance. A zipped file containing code is available for review.

\clearpage
\newpage
\bibliographystyle{vancouver}
\bibliography{ref}

@incollection{allen2008introduction,
  title={An introduction to stochastic epidemic models},
  author={Allen, Linda JS},
  booktitle={Mathematical epidemiology},
  pages={81--130},
  year={2008},
  publisher={Springer}
}

@article{novak2021lyme,
  title={Lyme disease in the era of COVID-19: a delayed diagnosis and risk for complications},
  author={Novak, Cheryl B and Scheeler, Verna M and Aucott, John N},
  journal={Case reports in infectious diseases},
  volume={2021},
  number={1},
  pages={6699536},
  year={2021},
  publisher={Wiley Online Library}
}

@article{mccormick2021effects,
  title={Effects of COVID-19 pandemic on reported Lyme disease, United States, 2020},
  author={McCormick, David W and Kugeler, Kiersten J and Marx, Grace E and Jayanthi, Praveena and Dietz, Stephanie and Mead, Paul and Hinckley, Alison F},
  journal={Emerging infectious diseases},
  volume={27},
  number={10},
  pages={2715},
  year={2021}
}

@article{borșan2021recreational,
  title={Recreational behaviour, risk perceptions, and protective practices against ticks: a cross-sectional comparative study before and during the lockdown enforced by the COVID-19 pandemic in Romania},
  author={Borșan, Silvia-Diana and Trif, Sabina Ramona and Mihalca, Andrei Daniel},
  journal={Parasites \& vectors},
  volume={14},
  number={1},
  pages={423},
  year={2021},
  publisher={Springer}
}

@article{jones2024lyme,
  title={Lyme Disease Under-Ascertainment During the COVID-19 Pandemic in the United States: Retrospective Study},
  author={Jones, Brie S and DeWitt, Michael E and Wenner, Jennifer J and Sanders, John W},
  journal={JMIR Public Health and Surveillance},
  volume={10},
  number={1},
  pages={e56571},
  year={2024},
  publisher={JMIR Publications Inc., Toronto, Canada}
}

@article{jore2023outdoor,
  title={Outdoor recreation, tick borne encephalitis incidence and seasonality in Finland, Norway and Sweden during the COVID-19 pandemic (2020/2021)},
  author={Jore, Solveig and Viljugrein, Hildegunn and Hjertqvist, Marika and Dub, Timoth{\'e}e and M{\"a}kel{\"a}, Henna},
  journal={Infection Ecology \& Epidemiology},
  volume={13},
  number={1},
  pages={2281055},
  year={2023},
  publisher={Taylor \& Francis}
}

@article{ghaffarzadegan2021simulation,
  title={Simulation-based what-if analysis for controlling the spread of Covid-19 in universities},
  author={Ghaffarzadegan, Navid},
  journal={PloS one},
  volume={16},
  number={2},
  pages={e0246323},
  year={2021},
  publisher={Public Library of Science San Francisco, CA USA}
}

@article{nikbakht2019comparison,
  title={Comparison of methods to estimate basic reproduction number (R0) of influenza, using Canada 2009 and 2017-18 A (H1N1) data},
  author={Nikbakht, Roya and Baneshi, Mohammad Reza and Bahrampour, Abbas and Hosseinnataj, Abolfazl},
  journal={Journal of Research in Medical Sciences},
  volume={24},
  number={1},
  pages={67},
  year={2019},
  publisher={Medknow}
}

@article{d2020assessment,
  title={Assessment of the SARS-CoV-2 basic reproduction number, R0, based on the early phase of COVID-19 outbreak in Italy},
  author={D'Arienzo, Marco and Coniglio, Angela},
  journal={Biosafety and health},
  volume={2},
  number={2},
  pages={57--59},
  year={2020},
  publisher={Elsevier}
}

@article{stephens2020covid,
  title={COVID-19 and the Path to Immunity},
  author={Stephens, David S and McElrath, M Juliana},
  journal={Jama},
  volume={324},
  number={13},
  pages={1279--1281},
  year={2020},
  publisher={American Medical Association}
}

@article{hernandez1999threshold,
  title={Threshold parameters and metapopulation persistence},
  author={Hern{\'a}ndez-Su{\'a}rez, Carlos M and Marquet, Pablo A and Velasco-Hern{\'a}ndez, Jorge X},
  journal={Bulletin of Mathematical Biology},
  volume={61},
  number={2},
  pages={341--353},
  year={1999},
  publisher={Elsevier}
}

@misc{MATLAB21a,
  author  = {The MathWorks Inc.},
  title   = {MATLAB version: 9.10.0 (R2021a)}, 
  publisher = {The MathWorks Inc.},
  address = {Natick, Massachusetts, United States},
  year    = {2021}, 
  url     = {https://www.mathworks.com}
}

@book{allen2007introduction,
  title={An Introduction to Mathematical Biology},
  author={Allen, Linda J.S.},
  year={2007},
  publisher={Pearson/Prentice Hall},
  address={Upper Saddle River, NJ}
}

@article{bellu2007daisy,
  title={{DAISY: a new software tool to test global identifiability of biological and physiological systems}},
  author={Bellu, Giuseppina and Saccomani, Maria Pia and Audoly, Stefania and D’Angi{\`o}, Leontina},
  journal={Computer Methods and Programs in Biomedicine},
  volume={88},
  number={1},
  pages={52--61},
  year={2007},
  publisher={Elsevier},
  DOI = {10.1016/j.cmpb.2007.07.002}
}

@article{tuncer2018structural,
  title={Structural and practical identifiability analysis of outbreak models},
  author={Tuncer, Necibe and Le, Trang T},
  journal={Mathematical biosciences},
  volume={299},
  pages={1--18},
  year={2018},
  publisher={Elsevier}
}

@article{dong2023differential,
  title={Differential elimination for dynamical models via projections with applications to structural identifiability},
  author={Dong, Ruiwen and Goodbrake, Christian and Harrington, Heather A and Pogudin, Gleb},
  journal={SIAM Journal on Applied Algebra and Geometry},
  volume={7},
  number={1},
  pages={194--235},
  year={2023},
DOI={https://doi.org/10.1137/22M1469067},
  publisher={SIAM}
}

@article{dong2020interactive,
  title={An interactive web-based dashboard to track COVID-19 in real time},
  author={Dong, Ensheng and Du, Hongru and Gardner, Lauren},
  journal={The Lancet Infectious Diseases},
  volume={20},
  number={5},
  pages={533--534},
  year={2020},
  publisher={Elsevier}
}

@misc{CDC_FluView_YYYY,
  author = {{Centers for Disease Control and Prevention (CDC)}},
  title = {{FluView: A Weekly Influenza Surveillance Report}},
publisher = {CDC},
  howpublished = "Available at \url{https://gis.cdc.gov/grasp/fluview/fluportaldashboard.html}",
  year = {2025},
  note = "Accessed: July 30, 2025"
}

@article{osi2025simultaneous,
  title={A simultaneous simulation of human behavior dynamics and epidemic spread: A multi-country study amidst the {COVID-19 }pandemic},
  author={Osi, Ann and Ghaffarzadegan, Navid},
  journal={Mathematical Biosciences},
  volume={380},
  pages={109368},
  year={2025},
  publisher={Elsevier}
}

@article{nikin2018unraveling,
  title={Unraveling within-host signatures of dengue infection at the population level},
  author={Nikin-Beers, Ryan and Blackwood, Julie C and Childs, Lauren M and Ciupe, Stanca M},
  journal={Journal of Theoretical Biology},
  volume={446},
  pages={79--86},
  year={2018},
  publisher={Elsevier}
}

@article{miao2011identifiability,
  title={On identifiability of nonlinear ODE models and applications in viral dynamics},
  author={Miao, Hongyu and Xia, Xiaohua and Perelson, Alan S and Wu, Hulin},
  journal={SIAM Review},
  volume={53},
  number={1},
  pages={3--39},
  year={2011},
  publisher={SIAM}
}

@article{mashayekhi2025dynamics,
  title={Dynamics of COVID-19: Exploring Behavioral Responsiveness},
  author={Mashayekhi, Ali and Gordon, Daniel and Tomoaia-Cotisel, Andrada and Bahaddin, Babak and Kim, Hyunjung and Luna-Reyes, Luis Felipe and Andersen, David F},
  journal={System Dynamics Review},
  volume={41},
  number={3},
  pages={e70006},
  year={2025},
  publisher={Wiley Online Library}
}

@article{saucedo2024comparative,
  title={Comparative analysis of practical identifiability methods for an SEIR model},
  author={Saucedo, Omar and Laubmeier, Amanda and Tang, Tingting and Levy, Benjamin and Asik, Lale and Pollington, Tim and Feldman, Olivia Prosper},
  journal={AIMS Mathematics},
  volume={9},
  number={9},
  pages={24722--24761},
  year={2024}
}

@article{gog2002dynamics,
  title={Dynamics and selection of many-strain pathogens},
  author={Gog, Julia R and Grenfell, Bryan T},
  journal={Proceedings of the National Academy of Sciences},
  volume={99},
  number={26},
  pages={17209--17214},
  year={2002},
  publisher={National Academy of Sciences}
}

@article{villaverde2016structural,
  title={Structural identifiability of dynamic systems biology models},
  author={Villaverde, Alejandro F and Barreiro, Antonio and Papachristodoulou, Antonis},
  journal={PLoS computational biology},
  volume={12},
  number={10},
  pages={e1005153},
  year={2016},
  publisher={Public Library of Science San Francisco, CA USA}
}

@article{diaz2023controllability,
  title={Controllability and accessibility analysis of nonlinear biosystems},
  author={D{\'\i}az-Seoane, Sandra and Blas, Antonio Barreiro and Villaverde, Alejandro F},
  journal={Computer methods and programs in biomedicine},
  volume={242},
  pages={107837},
  year={2023},
  publisher={Elsevier}
}

@article{lejeune2023effect,
  title={Effect of cross-immunity in a two-strain cholera model with aquatic component},
  author={LeJeune, Leah and Browne, Cameron},
  journal={Mathematical Biosciences},
  volume={365},
  pages={109086},
  year={2023},
  publisher={Elsevier}
}

@article{qiu2022understanding,
  title={Understanding the coevolution of mask wearing and epidemics: A network perspective},
  author={Qiu, Zirou and Espinoza, Baltazar and Vasconcelos, Vitor V and Chen, Chen and Constantino, Sara M and Crabtree, Stefani A and Yang, Luojun and Vullikanti, Anil and Chen, Jiangzhuo and Weibull, J{\"o}rgen and others},
  journal={Proceedings of the National Academy of Sciences},
  volume={119},
  number={26},
  pages={e2123355119},
  year={2022},
  publisher={National Academy of Sciences}
}

@article{epstein2021triple,
  title={Triple contagion: a two-fears epidemic model},
  author={Epstein, Joshua M and Hatna, Erez and Crodelle, Jennifer},
  journal={Journal of the Royal Society Interface},
  volume={18},
  number={181},
  pages={20210186},
  year={2021},
  publisher={The Royal Society}
}

@article{espinoza2024adaptive,
  title={Adaptive human behaviour modulates the impact of immune life history and vaccination on long-term epidemic dynamics},
  author={Espinoza, Baltazar and Saad-Roy, Chadi M and Grenfell, Bryan T and Levin, Simon A and Marathe, Madhav},
  journal={Proceedings B},
  volume={291},
  number={2033},
  pages={20241772},
  year={2024},
  publisher={The Royal Society}
}

@article{saad2023dynamics,
  title={Dynamics in a behavioral--epidemiological model for individual adherence to a nonpharmaceutical intervention},
  author={Saad-Roy, Chadi M and Traulsen, Arne},
  journal={Proceedings of the National Academy of Sciences},
  volume={120},
  number={44},
  pages={e2311584120},
  year={2023},
  publisher={National Academy of Sciences}
}

@article{pant2024mathematical,
  title={{Mathematical assessment of the role of human behavior changes on SARS-CoV-2 transmission dynamics in the united states}},
  author={Pant, Binod and Safdar, Salman and Santillana, Mauricio and Gumel, Abba B},
  journal={Bulletin of Mathematical Biology},
  volume={86},
  number={8},
  pages={92},
  year={2024},
  publisher={Springer}
}

@article{abbas2022evolution,
  title={{Evolution and consequences of individual responses during the COVID-19 outbreak}},
  author={Abbas, Wasim and MA, Masud and Park, Anna and Parveen, Sajida and Kim, Sangil},
  journal={PLoS One},
  volume={17},
  number={9},
  pages={e0273964},
  year={2022},
  publisher={Public Library of Science San Francisco, CA USA}
}

@article{rahmandad2022missing,
  title={{A missing behavioural feedback in COVID-19 models is the key to several puzzles}},
  author={Rahmandad, Hazhir and Xu, Ran and Ghaffarzadegan, Navid},
  journal={BMJ Global Health},
  volume={7},
  number={10},
  pages={e010463},
  year={2022},
  publisher={BMJ Specialist Journals}
}

@article{lim2023similar,
  title={Why Similar Policies Resulted In Different COVID-19 Outcomes: How Responsiveness And Culture Influenced Mortality Rates},
  author={Lim, Tse Yang and Xu, Ran and Ruktanonchai, Nick and Saucedo, Omar and Childs, Lauren M and Jalali, Mohammad S and Rahmandad, Hazhir and Ghaffarzadegan, Navid},
  journal={Health Affairs},
  volume={42},
  number={12},
  pages={1637--1646},
  year={2023}
}

@article{verelst2016behavioural,
  title={Behavioural change models for infectious disease transmission: a systematic review (2010--2015)},
  author={Verelst, Frederik and Willem, Lander and Beutels, Philippe},
  journal={Journal of The Royal Society Interface},
  volume={13},
  number={125},
  pages={20160820},
  year={2016},
  publisher={The Royal Society}
}

@article{funk2010modelling,
  title={Modelling the influence of human behaviour on the spread of infectious diseases: a review},
  author={Funk, Sebastian and Salath{\'e}, Marcel and Jansen, Vincent AA},
  journal={Journal of the Royal Society Interface},
  volume={7},
  number={50},
  pages={1247--1256},
  year={2010},
  publisher={The Royal Society}
}

@article{funk2015nine,
  title={Nine challenges in incorporating the dynamics of behaviour in infectious diseases models},
  author={Funk, Sebastian and Bansal, Shweta and Bauch, Chris T and Eames, Ken TD and Edmunds, W John and Galvani, Alison P and Klepac, Petra},
  journal={Epidemics},
  volume={10},
  pages={21--25},
  year={2015},
  publisher={Elsevier}
}

@article{kryazhimskiy2007state,
  title={{On state-space reduction in multi-strain pathogen models, with an application to antigenic drift in influenza A}},
  author={Kryazhimskiy, Sergey and Dieckmann, Ulf and Levin, Simon A and Dushoff, Jonathan},
  journal={PLoS Computational Biology},
  volume={3},
  number={8},
  pages={e159},
  year={2007},
  publisher={Public Library of Science San Francisco, USA}
}

@article{eletreby2020effects,
  title={The effects of evolutionary adaptations on spreading processes in complex networks},
  author={Eletreby, Rashad and Zhuang, Yong and Carley, Kathleen M and Ya{\u{g}}an, Osman and Poor, H Vincent},
  journal={Proceedings of the National Academy of Sciences},
  volume={117},
  number={11},
  pages={5664--5670},
  year={2020},
  publisher={National Academy of Sciences}
}

@article{yaqinuddin2020cross,
  title={{Cross-immunity between respiratory coronaviruses may limit COVID-19 fatalities}},
  author={Yaqinuddin, Ahmed},
  journal={Medical Hypotheses},
  volume={144},
  pages={110049},
  year={2020},
  publisher={Elsevier}
}

@article{welsh2010heterologous,
  title={Heterologous immunity between viruses},
  author={Welsh, Raymond M and Che, Jenny W and Brehm, Michael A and Selin, Liisa K},
  journal={Immunological Reviews},
  volume={235},
  number={1},
  pages={244--266},
  year={2010},
  publisher={Wiley Online Library}
}

@article{bhattacharyya2015cross,
  title={Cross-immunity between strains explains the dynamical pattern of paramyxoviruses},
  author={Bhattacharyya, Samit and Gesteland, Per H and Korgenski, Kent and Bj{\o}rnstad, Ottar N and Adler, Frederick R},
  journal={Proceedings of the National Academy of Sciences},
  volume={112},
  number={43},
  pages={13396--13400},
  year={2015},
  publisher={National Academy of Sciences}
}

@article{welsh2002no,
  title={{No one is naive: the significance of heterologous T-cell immunity}},
  author={Welsh, Raymond M and Selin, Liisa K},
  journal={Nature Reviews Immunology},
  volume={2},
  number={6},
  pages={417--426},
  year={2002},
  publisher={Nature Publishing Group UK London}
}

@article{kiseleva2021covid,
  title={{COVID-19 shuts doors to flu but keeps them open to rhinoviruses}},
  author={Kiseleva, Irina and Ksenafontov, Andrey},
  journal={Biology},
  volume={10},
  number={8},
  pages={733},
  year={2021},
  publisher={MDPI}
}

@article{perez2020dramatic,
  title={{Dramatic decrease of laboratory-confirmed influenza A after school closure in response to COVID-19}},
  author={Perez-Lopez, Andres and Hasan, Mohammad and Iqbal, Muhammad and Janahi, Mohammed and Roscoe, Diane and Tang, Patrick},
  journal={Pediatric Pulmonology},
  volume={55},
  number={9},
  pages={2233},
  year={2020}
}

@article{solomon2020influenza,
  title={{Influenza in the COVID-19 Era}},
  author={Solomon, Daniel A and Sherman, Amy C and Kanjilal, Sanjat},
  journal={JAMA},
  volume={324},
  number={13},
  pages={1342--1343},
  year={2020},
  publisher={American Medical Association}
}

@article{jones2020covid,
  title={{How COVID-19 is changing the cold and flu season}},
  author={Jones, Nicola and others},
  journal={Nature},
  volume={588},
  number={7838},
  pages={388--390},
  year={2020},
  publisher={Springer Science and Business Media LLC}
}

@article{rubin2020happens,
  title={{What happens when COVID-19 collides with flu season?}},
  author={Rubin, Rita},
  journal={JAMA},
  volume={324},
  number={10},
  pages={923--925},
  year={2020},
  publisher={American Medical Association}
}

@article{moore2014meteorological,
  title={{Meteorological influences on the seasonality of Lyme disease in the United States}},
  author={Moore, Sean M and Eisen, Rebecca J and Monaghan, Andrew and Mead, Paul},
  journal={The American Journal of Tropical Medicine and Hygiene},
  volume={90},
  number={3},
  pages={486},
  year={2014}
}

@article{lofgren2007influenza,
  title={Influenza seasonality: underlying causes and modeling theories},
  author={Lofgren, Eric and Fefferman, Nina H and Naumov, Yuri N and Gorski, Jack and Naumova, Elena N},
  journal={Journal of Virology},
  volume={81},
  number={11},
  pages={5429--5436},
  year={2007},
  publisher={American Society for Microbiology}
}

@article{neumann2022seasonality,
  title={Seasonality of influenza and other respiratory viruses},
  author={Neumann, Gabriele and Kawaoka, Yoshihiro},
  journal={EMBO Molecular Medicine},
  volume={14},
  number={4},
  pages={e15352},
  year={2022}
}

@article{he2020seir,
  title={{SEIR modeling of the COVID-19 and its dynamics}},
  author={He, Shaobo and Peng, Yuexi and Sun, Kehui},
  journal={Nonlinear Dynamics},
  volume={101},
  pages={1667--1680},
  year={2020},
  publisher={Springer}
}

@article{li2020substantial,
  title={{Substantial undocumented infection facilitates the rapid dissemination of novel coronavirus (SARS-CoV-2)}},
  author={Li, Ruiyun and Pei, Sen and Chen, Bin and Song, Yimeng and Zhang, Tao and Yang, Wan and Shaman, Jeffrey},
  journal={Science},
  volume={368},
  number={6490},
  pages={489--493},
  year={2020},
  publisher={American Association for the Advancement of Science}
}

@article{coburn2009modeling,
  title={Modeling influenza epidemics and pandemics: insights into the future of swine flu {(H1N1)}},
  author={Coburn, Brian J and Wagner, Bradley G and Blower, Sally},
  journal={BMC Medicine},
  volume={7},
  pages={1--8},
  year={2009},
  publisher={Springer}
}

@article{cramer2022evaluation,
  title={{Evaluation of individual and ensemble probabilistic forecasts of COVID-19 mortality in the United States}},
  author={Cramer, Estee Y and Ray, Evan L and Lopez, Velma K and Bracher, Johannes and Brennen, Andrea and Castro Rivadeneira, Alvaro J and Gerding, Aaron and Gneiting, Tilmann and House, Katie H and Huang, Yuxin and others},
  journal={Proceedings of the National Academy of Sciences},
  volume={119},
  number={15},
  pages={e2113561119},
  year={2022},
  publisher={National Academy of Sciences}
}

@article{traulsen2023individual,
  title={{Individual costs and societal benefits of interventions during the COVID-19 pandemic}},
  author={Traulsen, Arne and Levin, Simon A and Saad-Roy, Chadi M},
  journal={Proceedings of the National Academy of Sciences},
  volume={120},
  number={24},
  pages={e2303546120},
  year={2023},
  publisher={National Academy of Sciences}
}

@article{mcbryde2020role,
  title={{Role of modelling in COVID-19 policy development}},
  author={McBryde, Emma S and Meehan, Michael T and Adegboye, Oyelola A and Adekunle, Adeshina I and Caldwell, Jamie M and Pak, Anton and Rojas, Diana P and Williams, Bridget M and Trauer, James M},
  journal={Paediatric Respiratory Reviews},
  volume={35},
  pages={57--60},
  year={2020},
  publisher={Elsevier}
}

@article{blackwood2018introduction,
  title={An introduction to compartmental modeling for the budding infectious disease modeler},
  author={Blackwood, Julie C and Childs, Lauren M},
  journal={Letters in Biomathematics},
  volume={5},
  number={1},
  pages={195--221},
  year={2018}
}

@article{anderson1979population,
  title={Population biology of infectious diseases: Part I},
  author={Anderson, Roy M and May, Robert M},
  journal={Nature},
  volume={280},
  number={5721},
  pages={361--367},
  year={1979},
  publisher={Nature Publishing Group}
}

@article{kermack1927contribution,
  title={A contribution to the mathematical theory of epidemics},
  author={Kermack, William Ogilvy and McKendrick, Anderson G},
  journal={Proceedings of the Royal Society of London. Series A, Containing papers of a mathematical and physical character},
  volume={115},
  number={772},
  pages={700--721},
  year={1927},
  publisher={The Royal Society London}
}

@article{ross1917applicationIII,
  title={An application of the theory of probabilities to the study of a priori pathometry.—Part III},
  author={Ross, Ronald and Hudson, Hilda P},
  journal={Proceedings of the Royal Society of London. Series A, Containing papers of a mathematical and physical character},
  volume={93},
  number={650},
  pages={225--240},
  year={1917},
  publisher={The Royal Society London}
}

@article{rahmandad2022enhancing,
  title={Enhancing long-term forecasting: Learning from {COVID-19} models},
  author={Rahmandad, Hazhir and Xu, Ran and Ghaffarzadegan, Navid},
  journal={PLOS Computational Biology},
  volume={18},
  number={5},
  pages={e1010100},
  year={2022},
  publisher={Public Library of Science San Francisco, CA USA}
}

@article{lejeune2025formulating,
  title={Formulating human risk response in epidemic models: exogenous vs endogenous approaches},
  author={LeJeune, Leah and Ghaffarzadegan, Navid and Childs, Lauren M and Saucedo, Omar},
  journal={European Journal of Operational Research},
  year={2025},
  publisher={Elsevier}
}

@article{lejeune2024mathematical,
  title={Mathematical analysis of simple behavioral epidemic models},
  author={LeJeune, Leah and Ghaffarzadegan, Navid and Childs, Lauren M and Saucedo, Omar},
  journal={Mathematical Biosciences},
  volume={375},
  pages={109250},
  year={2024},
  publisher={Elsevier}
}

@article{nerlove1958,
    author = {Nerlove, Marc},
    title = {Adaptive Expectations and Cobweb Phenomena},
    journal = {The Quarterly Journal of Economics},
    volume = {72},
    number = {2},
    pages = {227-240},
    year = {1958},
    doi = {10.2307/1880597}
}

@book{sterman2000,
  author    = {John D. Sterman},
  title     = {Business Dynamics: Systems Thinking and Modeling for a Complex World},
  publisher = {Irwin/McGraw-Hill},
  year      = {2000},
  address   = {Boston}
}

@article{thompson2007polio,
  author  = {Thompson, Kimberly M. and Tebbens, Radboud J. Duintjer},
  title   = {Eradication versus control for poliomyelitis: An economic analysis},
  journal = {The Lancet},
  volume  = {369},
  number  = {9570},
  pages   = {1363--1371},
  year    = {2007},
  doi     = {10.1016/S0140-6736(07)60532-4}
}

@article{kraemer2025artificial,
  title={Artificial intelligence for modelling infectious disease epidemics},
  author={Kraemer, Moritz UG and Tsui, Joseph L-H and Chang, Serina Y and Lytras, Spyros and Khurana, Mark P and Vanderslott, Samantha and Bajaj, Sumali and Scheidwasser, Neil and Curran-Sebastian, Jacob Liam and Semenova, Elizaveta and others},
  journal={Nature},
  volume={638},
  number={8051},
  pages={623--635},
  year={2025},
  publisher={Nature Publishing Group UK London}
}

@article{huang2022game,
  title={Game-theoretic frameworks for epidemic spreading and human decision-making: A review},
  author={Huang, Yunhan and Zhu, Quanyan},
  journal={Dynamic Games and Applications},
  volume={12},
  number={1},
  pages={7--48},
  year={2022},
  publisher={Springer}
}

@article{chang2020game,
  title={Game theoretic modelling of infectious disease dynamics and intervention methods: a review},
  author={Chang, Sheryl L and Piraveenan, Mahendra and Pattison, Philippa and Prokopenko, Mikhail},
  journal={Journal of biological dynamics},
  volume={14},
  number={1},
  pages={57--89},
  year={2020},
  publisher={Taylor \& Francis}
}

@article{dandekar2020machine,
  title={A machine learning-aided global diagnostic and comparative tool to assess effect of quarantine control in COVID-19 spread},
  author={Dandekar, Raj and Rackauckas, Chris and Barbastathis, George},
  journal={Patterns},
  volume={1},
  number={9},
  year={2020},
  publisher={Elsevier}
}

@article{chen2025immunity,
  title={Immunity Debt for Seasonal Influenza After the COVID-19 Pandemic and as a Result of Nonpharmaceutical Interventions: An Ecological Analysis and Cohort Study},
  author={Chen, Li and Guo, Yuchen and L{\'o}pez-G{\"u}ell, Kim and Ma, Jun and Dong, Yanhui and Xie, Junqing and Alhambra, Daniel Prieto},
  journal={Advanced Science},
  volume={12},
  number={20},
  pages={2410513},
  year={2025},
  publisher={Wiley Online Library}
}

\newpage
\renewcommand{\thesection}{A\arabic{section}}
\setcounter{section}{0}
\renewcommand{\thefigure}{A\arabic{figure}}
\setcounter{figure}{0}
\renewcommand{\thetable}{A\arabic{table}}
\setcounter{table}{0}
\renewcommand{\theequation}{A\arabic{equation}}
\setcounter{equation}{0}

\section{Appendix}

\subsection{SIRS system for each scenario:}
\label{sec:app:eq}

\textbf{No spillover ($s=0$):}

\begin{equation}
\begin{aligned}\label{eqn:independent}
& \frac{d S_A}{d t}=-\beta_{0,A}\exp(-k\widetilde{I}_A) S_A I_A+\frac{R_A}{\tau_R}, \\
& \frac{d I_A}{d t}=\beta_{0,A}\exp(-k\widetilde{I}_A) S_A I_A-\frac{I_A}{\tau_I}, \\
& \frac{d R_A}{d t}=\frac{I_A}{\tau_I} -\frac{R_A}{\tau_R},\\
&\frac{d\widetilde{I}_A}{dt}=\frac{I_A-\widetilde{I}_A}{\tau_P}\\
& \frac{d S_B}{d t}=-\beta_{0,B}\exp(-k\widetilde{I}_B) S_B I_B+\frac{R_B}{\tau_R}, \\
& \frac{d I_B}{d t}=\beta_{0,B}\exp(-k\widetilde{I}_B) S_B I_B-\frac{I_B}{\tau_I}, \\
& \frac{d R_B}{d t}=\frac{I_B}{\tau_I} -\frac{R_B}{\tau_R},\\
&\frac{d\widetilde{I}_B}{dt}=\frac{I_B-\widetilde{I}_B}{\tau_P},\\
& S_i+I_i+R_i=1.
\end{aligned} 
\end{equation}

\textbf{Perfect spillover ($s=1$):}

\begin{equation}
\begin{aligned}\label{eqn:perfect}
& \frac{d S_A}{d t}=-\beta_{0,A}\exp(-k(\widetilde{I}_A+\widetilde{I}_B)) S_A I_A+\frac{R_A}{\tau_R}, \\
& \frac{d I_A}{d t}=\beta_{0,A}\exp(-k(\widetilde{I}_A+\widetilde{I}_B)) S_A I_A-\frac{I_A}{\tau_I}, \\
& \frac{d R_A}{d t}=\frac{I_A}{\tau_I} -\frac{R_A}{\tau_R},\\
&\frac{d\widetilde{I}_A}{dt}=\frac{I_A-\widetilde{I}_A}{\tau_P}\\
& \frac{d S_B}{d t}=-\beta_{0,B}\exp(-k(\widetilde{I}_A+\widetilde{I}_B)) S_B I_B+\frac{R_B}{\tau_R}, \\
& \frac{d I_B}{d t}=\beta_{0,B}\exp(-k(\widetilde{I}_A+\widetilde{I}_B)) S_B I_B-\frac{I_B}{\tau_I}, \\
& \frac{d R_B}{d t}=\frac{I_B}{\tau_I} -\frac{R_B}{\tau_R},\\
&\frac{d\widetilde{I}_B}{dt}=\frac{I_B-\widetilde{I}_B}{\tau_P},\\
& S_i+I_i+R_i=1.
\end{aligned} 
\end{equation}

\textbf{Imperfect spillover ($0<s<1$):}

\begin{equation}
\begin{aligned}\label{eqn:imperfect}
& \frac{d S_A}{d t}=-\beta_{0,A}\exp(-k\widetilde{I}_A)(1-s(1-\exp(-k\widetilde{I}_B))) S_A I_A+\frac{R_A}{\tau_R}, \\
& \frac{d I_A}{d t}=\beta_{0,A}\exp(-k\widetilde{I}_A)(1-s(1-\exp(-k\widetilde{I}_B))) S_A I_A-\frac{I_A}{\tau_I}, \\
& \frac{d R_A}{d t}=\frac{I_A}{\tau_I} -\frac{R_A}{\tau_R},\\
&\frac{d\widetilde{I}_A}{dt}=\frac{I_A-\widetilde{I}_A}{\tau_P}\\
& \frac{d S_B}{d t}=-\beta_{0,B}\exp(-k\widetilde{I}_B)(1-s(1-\exp(-k\widetilde{I}_A))) S_B I_B+\frac{R_B}{\tau_R}, \\
& \frac{d I_B}{d t}=\beta_{0,B}\exp(-k\widetilde{I}_B)(1-s(1-\exp(-k\widetilde{I}_A))) S_B I_B-\frac{I_B}{\tau_I}, \\
& \frac{d R_B}{d t}=\frac{I_B}{\tau_I} -\frac{R_B}{\tau_R},\\
&\frac{d\widetilde{I}_B}{dt}=\frac{I_B-\widetilde{I}_B}{\tau_P},\\
& S_i+I_i+R_i=1.
\end{aligned} 
\end{equation}

\subsection{Proofs of Theorems}

\noindent \textbf{Proof of Theorem \ref{thm:approximation}.}\\\\
With behavioral spillover, the transmission rates of our two diseases are
\begin{align*}
    \beta_A &= m_A(1-s(1-m_B))\beta_{0,A},\\
    \beta_B &= (1-s(1-m_A))m_B\beta_{0,B}.
\end{align*}
Multiplying both equations by $\tau_I$ gives
\begin{align*}
    %\mathcal{R}_{e,A} 
    \beta_A \tau_I &= m_A(1-s(1-m_B))\mathcal{R}_{0,A},\\
    %\mathcal{R}_{e,B} 
    \beta_B \tau_I &= (1-s(1-m_A))m_B\mathcal{R}_{0,B}.
\end{align*}
Multiplying both equations by the respective susceptible population and noticing that the left hand side gives effective reproduction values gives
\begin{align*}
    \mathcal{R}_{e,A} &= m_A(1-s(1-m_B))S_A\mathcal{R}_{0,A},\\
    \mathcal{R}_{e,B} &= (1-s(1-m_A))m_B S_B\mathcal{R}_{0,B}.
\end{align*}
Assuming persistence of strain $A$ such that $\mathcal{R}_{e,A}\approx 1$, then 
\begin{align*}
    m_A(1-s(1-m_B)) = \frac{1}{S_A\mathcal{R}_{0,A}},
\end{align*}
and loss of strain $B$ such that $\mathcal{R}_{e,B}<1$ which gives $I_B\approx 0$, $m_B\approx 1$ and $S_B\approx 1$ then
\begin{align*}
    m_A &= \frac{1}{S_A\mathcal{R}_{0,A}},\\
    \mathcal{R}_{e,B} &= (1-s(1-m_A)) \mathcal{R}_{0,B} < 1.
\end{align*}
Plugging the former into the latter  gives
\begin{align*}
    1-s+s \frac{1}{S_A\mathcal{R}_{0,A}} < \frac{1}{\mathcal{R}_{0,B} },
\end{align*}
which after rearranging becomes 
\begin{equation*}
    \frac{1-\frac{1}{\mathcal{R}_{0,B}}}{1-\frac{1}{S_A \mathcal{R}_{0,A}}}<s<1.
\end{equation*}
As smaller values of $S_A$ imply larger values of $\frac{1-\frac{1}{\mathcal{R}_{0,B}}}{1-\frac{1}{S_A \mathcal{R}_{0,A}}}$, we bound this by assuming $S_A \approx 1$ such that \begin{equation*}
    s_{\text{threshold}} = \frac{1-\frac{1}{\mathcal{R}_{0,B}}}{1-\frac{1}{\mathcal{R}_{0,A}}}<s<1.
\end{equation*}
Thus, $s_{\text{threshold}}$ is a necessary but not sufficient condition for exclusion of disease $B$.

\noindent \textbf{Proof of Theorem \ref{thm:equilibria}}
\\\\
To solve for equilibria, we set all equations with behavior equal to zero. In doing so, the equations for $\frac{d R_i}{d t}$ and $\frac{d\widetilde{I}_i}{dt}$ will always yield $R_i=\frac{\tau_R}{\tau_I}I_i$ and $\widetilde{I}_i=I_i$. Substituting these into the equations for $\frac{d S_i}{d t}$ (or $\frac{d I_i}{d t}$) yields the equations
\begin{align*}
    0&=(\beta_A(\cdot) S_A -\frac{1}{\tau_I})I_A \text{ and }\\
    0&=(\beta_B(\cdot) S_B -\frac{1}{\tau_I})I_B.
\end{align*}

\noindent \textbf{1. Disease-free equilibrium:} For each of the three scenarios, when $I_A=I_B=0$, we have the disease-free equilibrium given by
\[
(1,0,0,0,1,0,0,0).
\]

\noindent \textbf{2. Disease-$B$-free, disease-$A$-endemic boundary equilibrium:} For each of the three scenarios, when $I_B=0$ and $I_A>0$, then $\beta_A(\cdot) = \beta_{0,A}\exp(-k\widetilde{I}_A)$. The equilibrium is given by $S_B=1$, $I_B=R_B=\widetilde{I}_B=0$ and $S_A=1-I_A-R_A$, $R_A=\frac{\tau_R}{\tau_I}I_A$, $\widetilde{I}_A=I_A$ where $I_A$ satisfies 
\[0=\beta_{0,A}\exp(-kI_A) S_A -\frac{1}{\tau_I}.\]
Since $S_A$ depends on $I_A$, we cannot obtain a closed-form expression for $I_A$.
 We follow the same argument as in Theorem 2 of \cite{lejeune2024mathematical} to determine conditions of existence for this equilibrium.
We rewrite the equation (substituting in the expressions for $S_A$ and $R_A$) as
\begin{equation}\label{eqn:exp_R0}
    \exp(kI_A) =\beta_{0,A}\tau_I -\beta_{0,A}(\tau_I+\tau_R)I_A.
\end{equation}
Under the assumption of positive (i.e. biologically relevant) parameters, the left hand side of the expression is increasing monotonically from one while the right hand side is decreasing linearly from $\beta_{0,A}\tau_I$ (see Figure \ref{fig:equil_illustrations}). When $\beta_{0,A}\tau_I<1$ (Figure \ref{fig:No_EE}), the two curves will not intersect in the positive quadrant, resulting in no equilibrium of this type in a biologically relevant region. When $\beta_{0,A}\tau_I=1$ (Figure \ref{fig:DFE_BE}), the curves will intersect at one, which is exactly when $I_A=0$. This is equivalent to the disease-free equilibrium, so $\mathcal{E}_A^*$, $\hat{\mathcal{E}}_A$, $\bar{\mathcal{E}}_A$ do not exist when $\beta_{0,A}\tau_I\leq 1$. When $\beta_{0,A}\tau_I>1$ (Figure \ref{fig:EE}), the two curves intersect at exactly one point in the positive quadrant ($I_A = I^*_A, \hat{I}_A,\text{ or } \bar{I}_A$ at the point of intersection, for each respective scenario), resulting in a unique equilibrium ($\mathcal{E}_A$, $\hat{\mathcal{E}}_A$, $\bar{\mathcal{E}}_A$ for each respective scenario). Since the endemic equilibrium value for $I_A$ gives equality of the two sides of Equation \eqref{eqn:exp_R0}, it follows that $\beta_{0,A}\tau_I>\exp(kI^*_A)$, $\beta_{0,A}\tau_I>\exp(k\hat{I}_A)$, and $\beta_{0,A}\tau_I>\exp(k\bar{I}_A)$, for each respective scenario.\\

\noindent \textbf{3. Disease-$A$-free, disease-$B$-endemic boundary equilibrium:} Similar reasoning to the disease-$B$-free, disease-$A$-endemic setting shows that $\mathcal{E}_B^*$, $\hat{\mathcal{E}}_B$, and $\bar{\mathcal{E}}_B$ exist if and only if $\beta_{0,B}\tau_I >1$. When these equilibria exist, it follows that $\beta_{0,B}\tau_I>\exp(kI^*_B)$, $\beta_{0,B}\tau_I>\exp(k\hat{I}_B)$, and $\beta_{0,B}\tau_I>\exp(k\bar{I}_B)$.\\

\begin{figure}[h!]
\begin{subfigure}[t]{.32\textwidth}
\centering
    \includegraphics[width=\linewidth]{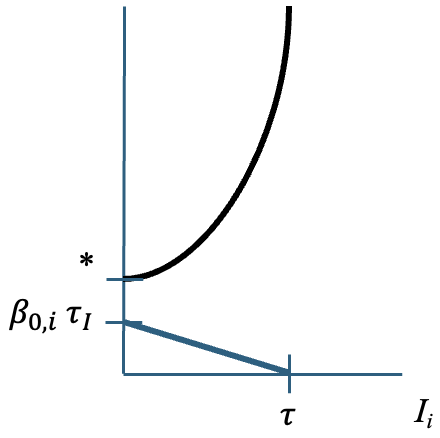}
    \caption{No boundary equilibrium}\label{fig:No_EE}
\end{subfigure}
\begin{subfigure}[t]{.32\textwidth}
\centering
    \includegraphics[width=\linewidth]{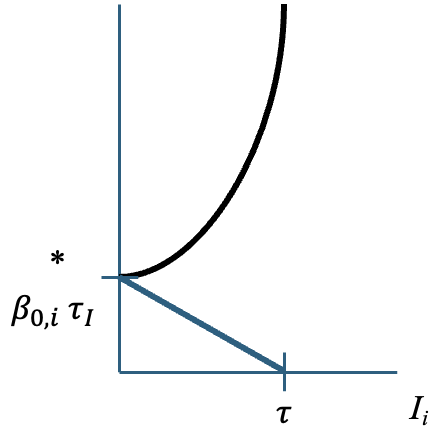}
    \caption{No boundary equilibrium}\label{fig:DFE_BE}
\end{subfigure}
\begin{subfigure}[t]{.32\textwidth}
\centering
    \includegraphics[width=\linewidth]
    {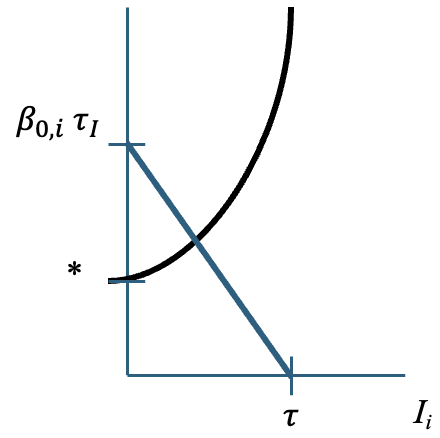}
    \caption{Boundary equilibrium exists }\label{fig:EE}
\end{subfigure}%
\caption{Illustration of the intersection of the general linear curve with the general exponential curve for strain $i$. Existence of the endemic equilibrium is determined by the relationship of $\beta_{0,i}\tau_I$ with *, the respective threshold value. At the intersection in the positive quadrant, $I_i = I^*_i,~ \hat{I}_i,~\bar{I}_i$, depending on the scenario. For (b), while the curves intersect, it occurs at $I_i=0$, so is the disease-free equilibrium.}
\label{fig:equil_illustrations}
\end{figure}

\noindent \textbf{4. Endemic equilibrium:} When $I_A>0$ and $I_B>0$, we consider the general case of imperfect spillover $(0<s<1)$. Results from this case can be extended to the cases of no spillover and perfect spillover by considering $s=0$ or $s=1$, respectively.

\noindent \textbf{Imperfect spillover:} If both $I_A\neq0$ and $I_B\neq0$, then the following two equations are simultaneously satisfied:
\begin{equation}\label{eqn:imperf0A}
 0 = \beta_{0,A}\exp(-kI_A)(1-s(1-\exp(-kI_B)))(1-I_A-\dfrac{\tau_R}{\tau_I}I_A)-\dfrac{1}{\tau_I}   
\end{equation}
and
\begin{equation}\label{eqn:imperf0B}
    0 = \beta_{0,B}\exp(-kI_B)(1-s(1-\exp(-kI_A)))(1-I_B-\dfrac{\tau_R}{\tau_I}I_B)-\dfrac{1}{\tau_I}.
\end{equation}
These equations can be rewritten, respectively, as 
\begin{equation}\label{eq:l2.1}
\dfrac{\exp(k(I_A+I_B))}{\exp(kI_B)-s(\exp(kI_B)-1)} = \beta_{0,A}\tau_I-\beta_{0,A}(\tau_I+\tau_R)I_A
\end{equation}
and
\begin{equation}\label{eq:l2.2}
\dfrac{\exp(k(I_A+I_B))}{\exp(kI_A)-s(\exp(kI_A)-1)} = \beta_{0,B}\tau_I-\beta_{0,B}(\tau_I+\tau_R)I_B.
\end{equation}
We work with Eqn. \eqref{eq:l2.1} (results for Eqn. \eqref{eq:l2.2} will follow in a similar manner). 
The left hand side can be rewritten to give
\begin{equation}\label{eqn:imperf}
\dfrac{\exp(kI_B)}{\exp(kI_B)-s(\exp(kI_B)-1)}\exp(kI_A) = \beta_{0,A}\tau_I-\beta_{0,A}(\tau_I+\tau_R)I_A.
\end{equation}
For Eqn. \eqref{eqn:imperf}, the right hand side is a plane that only depends on $I_A$ while the left hand side is a surface that grows exponentially with $I_A$ and $I_B$. By observation, these two surfaces are distinct. Thus, if these two surfaces intersect, they do so along a curve. The right hand side obtains its maximum value of $\mathcal{R}_{0,A}=\beta_{0,A}\tau_I$ when $I_A=0$ (for all values of $I_B$) and decreases as $I_A$ increases, ultimately becoming negative. The left hand side obtains its minimum value (with respect to $I_A$) of $\dfrac{\exp(kI_B)}{\exp(kI_B)-s(\exp(kI_B)-1)}$ when $I_A=0$ and is increasing with respect to $I_A$ and $I_B$. Thus, as long as 
\begin{equation}\label{eqn:thresholdcondition}
\mathcal{R}_{0,A}=\beta_{0,A}\tau_I>\dfrac{\exp(kI_B)}{\exp(kI_B)-s(\exp(kI_B)-1)},
\end{equation}
these surfaces will intersect for ($I_A,I_B$) pairs that satisfy Equation \eqref{eqn:imperf0A}. We can project the curve that results from the intersection of the surface and plane in Equation \eqref{eqn:imperf} onto the $I_A-I_B$ plane.

To obtain a qualitative picture of the curve, we first rewrite Equation \eqref{eqn:imperf} as
\begin{equation}\label{eqn:imperf1}
\dfrac{\exp(kI_B)}{\exp(kI_B)-s(\exp(kI_B)-1)} = \exp(-kI_A)(\beta_{0,A}\tau_I-\beta_{0,A}(\tau_I+\tau_R)I_A),
\end{equation}
such that the left hand side is only dependent on $I_B$ and the right hand side is only dependent on $I_A$.

The linear portion of right hand side (RHS) of Equation \eqref{eqn:imperf1} is a decreasing function of $I_A$. Since $0 \leq I_A \leq 1$, this linear factor is bounded above by $\beta_{0,A}\tau_I$ and bounded below by $-\beta_{0,A}\tau_R$. The exponential portion of the RHS is also decreasing in $I_A$, bounded above by 1 and below by $e^{-k}$. Hence,
\[
-e^{-k}\beta_{0,A}\tau_R \leq RHS \leq \beta_{0,A}\tau_I.
\]
As $I_A$ increases, the RHS decreases, so the LHS must also decrease to maintain equality. The derivative of LHS with respect to $I_B$ is positive, so LHS decreases as $I_B$ decreases.

When $I_B=0$, Equation \eqref{eqn:imperf1} reduces to
\begin{equation}\label{eqn:boundary}
1 = \exp(-kI_A)(\beta_{0,A}\tau_I-\beta_{0,A}(\tau_I+\tau_R)I_A),
\end{equation}
which has the unique solution $\bar{I}_{A}$ (the disease $A$ boundary equilibrium which is identical for any level of spillover), which exists if $\beta_{0,A}\tau_I>1$. Thus, the projected curve intersects the $I_A$-axis at $(\bar{I}_A,0)$. As $I_B$ increases from 0, $I_A$ must decrease, moving towards the $I_B$-axis.

When $I_A=0$, Equation \eqref{eqn:imperf1} reduces to
\begin{equation}\label{eqn:A11}
\dfrac{\exp(kI_B)}{\exp(kI_B)-s(\exp(kI_B)-1)} = \beta_{0,A}\tau_I,
\end{equation}
which has the closed-form solution,
\begin{equation}\label{eqn:I1B}
I^1_B=\dfrac{1}{k}\ln\left(\dfrac{s\beta_{0,A}\tau_I}{1-\beta_{0,A}\tau_I(1-s)}\right)
\end{equation}
if $s\neq0$ and $\beta_{0,A}\tau_I(1-s)<1$, which can be rewritten as $1-\frac{1}{\mathcal{R}_{0,A}}<s$.

For $s=0$, Equation \eqref{eqn:imperf1} reduces to Equation \eqref{eqn:boundary}, independent of the values of $I_A$ and $I_B$, so the surfaces in \eqref{eqn:imperf1} intersect along the $\bar{I}_A$ line for all $I_B$. For $0<s\leq 1-\frac{1}{\mathcal{R}_{0,A}}$, \eqref{eqn:A11} does not have a real solution, so the surfaces intersect along a curve that asymptotes along some value $I^t_A$, where $0<I^t_A<\bar{I}_A$. It is known that $I^t_A<\bar{I}_A$ since the projected curves intersects the $I_A$ axis at $(\bar{I}_A,0)$ and increases in $I_B$ as $I_A$ decreases, i.e., moving towards the $I_B$-axis (and thus approaches $I^t_A$).
For $s> 1-\frac{1}{\mathcal{R}_{0,A}}$, the projected curve intersects the $I_B$-axis at $(0,I^1_B)$.

Following the same line of reasoning for Equation \eqref{eq:l2.2}, as long as 
\begin{equation}\label{eqn:thresholdcondition1}
\mathcal{R}_{0,B}=\beta_{0,B}\tau_I>\dfrac{\exp(kI_A)}{\exp(kI_A)-s(\exp(kI_A)-1)},
\end{equation}
the projection of the curve resulting from the intersection of the respective surface and plane hits the $I_B$-axis at the disease $B$ boundary equilibrium, $(0,\bar{I}_B)$. For $s=0$, the surfaces intersect along the $\bar{I}_B$ line for all $I_A$. For $0<s\leq 1-\frac{1}{\mathcal{R}_{0,B}}$, the surfaces intersect along a curve that asymptotes along some $I^t_B>0$, which is smaller than $\bar{I}_B$ (decreases in $I_B$ lead to increases in $I_A$ as discussed above).

To summarize, for no spillover ($s=0$), if $\mathcal{R}_{0,A},\mathcal{R}_{0,B}>1$, there is always an intersection of the projected curves, and thus a unique endemic equilibrium (Fig. \ref{fig:Orthog}).

For non-zero spillover, there are multiple cases to consider. For $0<s\leq 1-\frac{1}{\mathcal{R}_{0,B}}$ and $0<s\leq 1-\frac{1}{\mathcal{R}_{0,A}}$, there is always an intersection of the projected curves, and thus a unique endemic equilibrium (Fig. \ref{fig:2asym}). For $s>1-\frac{1}{\mathcal{R}_{0,B}}$ but $s\leq 1-\frac{1}{\mathcal{R}_{0,A}}$, there is an intersection if $\bar{I}_A<I_A^1$ (Fig. \ref{fig:1asym}). For $s>1-\frac{1}{\mathcal{R}_{0,A}}$ but $s\leq 1-\frac{1}{\mathcal{R}_{0,B}}$, there is an intersection if $\bar{I}_B<I_B^1$. For both $s>1-\frac{1}{\mathcal{R}_{0,B}}$ and $s>1-\frac{1}{\mathcal{R}_{0,A}}$, there is an intersection if $\bar{I}_A<I_A^1$ and $\bar{I}_B<I_B^1$ (Fig. \ref{fig:Noasym}). When both $s>1-\frac{1}{\mathcal{R}_{0,B}}$ and $s>1-\frac{1}{\mathcal{R}_{0,A}}$, but $\bar{I}_A>I_A^1$ or $\bar{I}_B>I_B^1$, the curves do not intersect and the endemic equilibrium does not exist (Fig. \ref{fig:NoCross}).

\begin{figure}[h!]
\begin{subfigure}[t]{.32\textwidth}
\centering
\includegraphics[width=\linewidth]{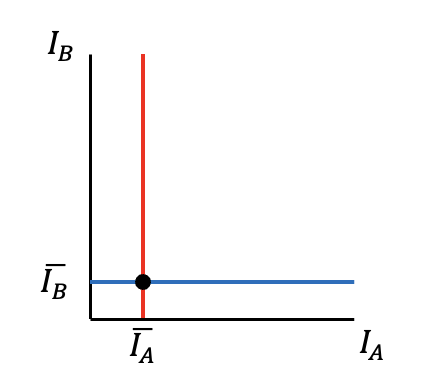}
    \caption{$s=0$ (no spillover)}\label{fig:Orthog}
\end{subfigure}
\begin{subfigure}[t]{.32\textwidth}
\centering
\includegraphics[width=\linewidth]{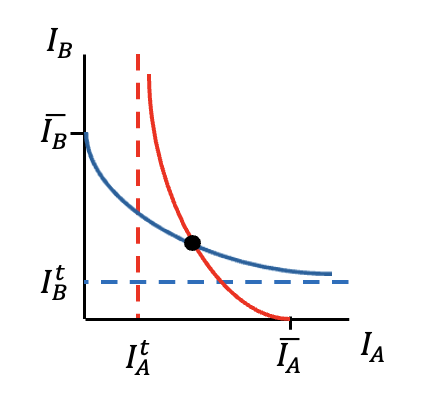}
    \caption{$0<s\leq 1-\frac{1}{\mathcal{R}_{0,A}},1-\frac{1}{\mathcal{R}_{0,B}}$}\label{fig:2asym}
\end{subfigure}
\begin{subfigure}[t]{.32\textwidth}
\centering
\includegraphics[width=\linewidth]{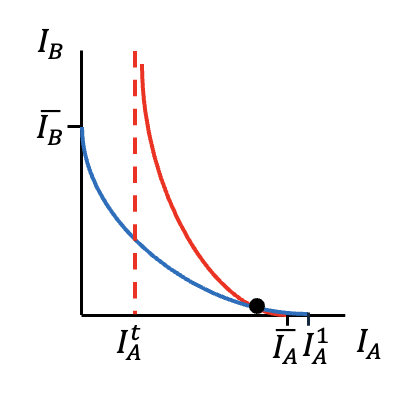}
    \caption{$1-\frac{1}{\mathcal{R}_{0,B}}<s\leq 1-\frac{1}{\mathcal{R}_{0,A}}$}\label{fig:1asym}
\end{subfigure}
\hspace{10cm}
\begin{subfigure}[t]{.32\textwidth}
\centering
\includegraphics[width=\linewidth]{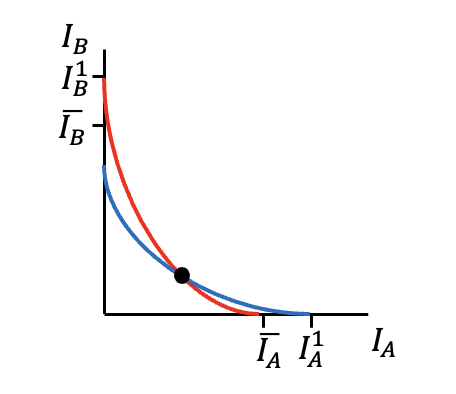}
    \caption{$1-\frac{1}{\mathcal{R}_{0,A}},1-\frac{1}{\mathcal{R}_{0,B}}<s$ and $\bar{I}_A<I^1_A$}\label{fig:Noasym}
\end{subfigure}
\begin{subfigure}[t]{.32\textwidth}
\centering
\includegraphics[width=\linewidth]{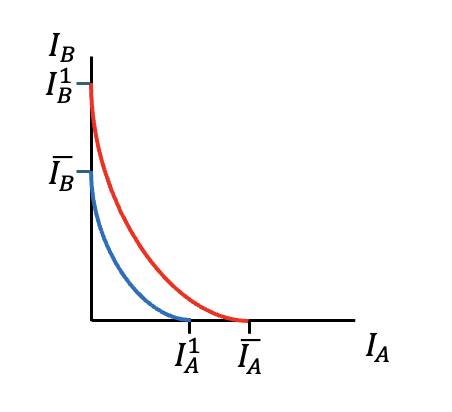}
    \caption{$1-\frac{1}{\mathcal{R}_{0,A}},1-\frac{1}{\mathcal{R}_{0,B}}<s$ and $\bar{I}_A>I^1_A$}\label{fig:NoCross}
\end{subfigure}
\caption{Illustration of the intersection of the projection into $I_A$-$I_B$ space of the curves created by the intersecting surfaces of Equations \eqref{eq:l2.1} (red lines) and \eqref{eq:l2.2} (blue lines) under the assumption that $\mathcal{R}_{0,A}>\mathcal{R}_{0,B}$ and for varying values of $s$. The intersection of the curves represents the endemic equilibrium values $(\bar{I}_{A_E},\bar{I}_{B_E})$.}
\label{fig:equil_illustrations2}
\end{figure}

From Equation \eqref{eqn:A11}, $\bar{I}_B<I_B^1$ implies
\begin{equation*}
\dfrac{\exp(k\bar{I}_B)}{\exp(k\bar{I}_B)-s(\exp(k\bar{I}_B)-1)} < \beta_{0,A}\tau_I,
\end{equation*}
since $\frac{\exp(kI_B)}{\exp(kI_B)-s(\exp(kI_B)-1)}$ is decreasing in $I_B$.

Thus, for the intersection of the projected curves, the conditions only require that Equation \eqref{eqn:thresholdcondition} is satisfied for $\bar{I}_B$, i.e.,
\[
\mathcal{R}_{0,A}=\beta_{0,A}\tau_I>\dfrac{\exp(k\bar{I}_B)}{\exp(k\bar{I}_B)-s(\exp(k\bar{I}_B)-1)},
\]
and that Equation \eqref{eqn:thresholdcondition1} is satisfied for $\bar{I}_A$, i.e., 
\[
\mathcal{R}_{0,B}=\beta_{0,B}\tau_I>\dfrac{\exp(k\bar{I}_A)}{\exp(k\bar{I}_A)-s(\exp(k\bar{I}_A)-1)}.
\]

For perfect spillover ($s=1$), these conditions simplify to $\mathcal{R}_{0,A}>\exp(k\bar{I}_B)$ and $\mathcal{R}_{0,B}>\exp(k\bar{I}_A)$.

\noindent\textbf{Proof of Theorem \ref{thm:stabilityDFE}}\\\\
We derive the basic reproduction numbers $\mathcal{R}_{0,i}$. The Jacobian matrices for each subsystem, evaluated at the disease free equilibrium, are equivalent and can be written as
\[
J(DFE) = \begin{pmatrix}
        - \frac{1}{\tau_R} & \beta_{0,A} - \frac{1}{\tau_R} & 0 & 0 & 0 & 0\\
        0 & \beta_{0,A}-\frac{1}{\tau_I} & 0 & 0 & 0  & 0\\
        0 & \frac{1}{\tau_P} & -\frac{1}{\tau_P} & 0 & 0 & 0\\
        0 & 0 & 0 & -\frac{1}{\tau_R} & \beta_{0,B}-\frac{1}{\tau_R} & 0\\
        0 & 0 & 0 & 0 & \beta_{0,B}-\frac{1}{\tau_I} & 0\\
        0 & 0 & 0 & 0 & \frac{1}{\tau_P} & -\frac{1}{\tau_P}\\
    \end{pmatrix}.
\]
The only nonnegative eigenvalues are and $\beta_{0,i}-\dfrac{1}{\tau_I}$. It is sufficient that $\beta_{0,i}-\dfrac{1}{\tau_I}<0$ for all $i$ in order for the DFE to be stable.

\noindent \textbf{Proof of Theorem \ref{thm:BE}}\\\\
Without loss of generality, we assume $\mathcal{R}_{0,A}>1$ and consider the disease $B$-free, disease-$A$-endemic equilibrium for each subsystem ($\mathcal{E}_{A}$, $\hat{\mathcal{E}}_{A}$, $\bar{\mathcal{E}}_{A}$),
    \[
    \left(1-I_A-\dfrac{\tau_R}{\tau_I}I_A, I_A,I_A,1,0,0 \right).
    \] 
    The Jacobian matrices of the subsystems ($J$, $\hat{J}$ and $\bar{J}$, respectively), evaluated at the respective equilibrium ($\mathcal{E}_{A}$, $\hat{\mathcal{E}}_{A}$, $\bar{\mathcal{E}}_{A}$), can be written as 
    \[
    J(\mathcal{E}_A) =
    \begin{pmatrix}
        \mathbf{J_1} &  \mathbf{J_2} \\
          \mathbf{0} &   \mathbf{J_3}
    \end{pmatrix},
    ~J(\hat{\mathcal{E}}_A) =
    \begin{pmatrix}
         \mathbf{J_1} &  \mathbf{\hat{J}_2} \\
        \mathbf{0} &   \mathbf{\hat{J}_3}
    \end{pmatrix},
    \text{ and }~
    J(\bar{\mathcal{E}}_A) =
    \begin{pmatrix}
         \mathbf{J_1} &   \mathbf{\bar{J}_2}\\
         \mathbf{0} &  \mathbf{\bar{J}_3}
    \end{pmatrix}.
    \]
    Since each block Jacobian matrix is upper triangular, to determine the stability of the respective equilibrium, we need only consider the eigenvalues of the block matrices on the diagonals. The upper left block matrix is the same for all three Jacobians and is given by
\[
 \mathbf{J_1} = 
\begin{pmatrix}
    -\beta_{0,A}\exp(-kI_A)I_A-\frac{1}{\tau_R} & -\beta_{0,A}\exp(-kI_A)S_A -\frac{1}{\tau_R} & k\beta_{0,A}\exp(-kI_A)S_AI_A \\
    \beta_{0,A}\exp(-kI_A)I_A & \beta_{0,A}\exp(-kI_A)S_A-\frac{1}{\tau_I} & -k\beta_{0,A}\exp(-kI_A)S_AI_A \\
    0 & \frac{1}{\tau_P} & -\frac{1}{\tau_P} \\
\end{pmatrix},
\]
where $S_A=1-I_A-\dfrac{\tau_R}{\tau_I}I_A$, with $I_A=I^*_A$ for $\mathcal{E}_A$, $I_A=\hat{I}_A$ for $\hat{\mathcal{E}}_A$, and $I_A=\bar{I}_A$ for $\bar{\mathcal{E}}_A$. Since we are assuming disease $A$ is endemic, it follows that $\beta_{0,A}\exp(-kI_A)S_A-\dfrac{1}{\tau_I}=0$. Hence, we can simplify $\mathbf{J_1}$:
\[
 \mathbf{J_1} = 
\begin{pmatrix}
    -\beta_{0,A}\exp(-kI_A)I_A-\frac{1}{\tau_R} & -\frac{1}{\tau_I}-\frac{1}{\tau_R} & \frac{kI_A}{\tau_I} \\
    \beta_{0,A}\exp(-kI_A)I_A & 0 & -\frac{kI_A}{\tau_I} \\
    0 & \frac{1}{\tau_P} & -\frac{1}{\tau_P} \\
\end{pmatrix}.
\]
The characteristic polynomial of $\mathbf{J_1}$ is
\begin{align*}
0= \lambda^3+a_1\lambda^2+a_2\lambda+a_3
\end{align*}
where
\begin{align*}
&a_1 = \frac{1}{\tau_P}+\frac{1}{\tau_R}+\beta_{0,A}\exp\left(-kI_A\right)I_A,\\
&a_2 = \frac{kI_A}{\tau_I\tau_P}+\frac{1}{\tau_R\tau_P}+\frac{1}{\tau_P}\beta_{0,A}\exp\left(-kI_A\right)I_A+\frac{\beta_{0,A}\exp\left(-kI_A\right)I_A}{\tau_I}+\frac{\beta_{0,A}\exp\left(-kI_A\right)I_A}{\tau_R},\\
&a_3=\frac{\beta_{0,A}\exp\left(-kI_A\right)I_A}{\tau_I\tau_P}+\frac{\beta_{0,A}\exp\left(-kI_A\right)I_A}{\tau_R\tau_P}+\frac{kI_A}{\tau_I\tau_P\tau_R}.
\end{align*}
Since all coefficients are positive and $a_1a_2-a_3>0$, by the Routh-Hurwitz criteria \cite{allen2007introduction}, the remaining eigenvalues will have negative real part. We next consider the lower right block matrices:
\[
 \mathbf{J_3}   =\begin{pmatrix}
        -\frac{1}{\tau_R} & -\beta_{0,B}-\frac{1}{\tau_R} & 0\\
        0 & \beta_{0,B}-\frac{1}{\tau_I} & 0\\
        0 & \frac{1}{\tau_P} & -\frac{1}{\tau_P}\\
    \end{pmatrix},
\]
\[
\mathbf{\hat{J}_3}=    \begin{pmatrix}
        -\frac{1}{\tau_R} & -\beta_{0,B}\exp(-k\hat{I}_A)-\frac{1}{\tau_R}  & 0\\
        0 & \beta_{0,B}\exp(-k\hat{I}_A)-\frac{1}{\tau_I} & 0\\
        0 & \frac{1}{\tau_P} & -\frac{1}{\tau_P}\\
    \end{pmatrix},
\]
and
\[
\mathbf{\bar{J}_3}=\begin{pmatrix}
        -\frac{1}{\tau_R} & -\beta_{0,B}(1-s(1-\exp(-k\bar{I}_A)))-\frac{1}{\tau_R} & 0\\
        0 & \beta_{0,B}(1-s(1-\exp(-k\bar{I}_A)))-\frac{1}{\tau_I} & 0\\
        0 & \frac{1}{\tau_P} & -\frac{1}{\tau_P}\\
    \end{pmatrix}.
\]
From $ \mathbf{J_1}$ and $ \mathbf{J_3}$, $\mathcal{E}_{A}$ will be locally asymptotically stable if $\mathcal{R}_{0,B}<1$. From $ \mathbf{J_1}$ and $ \mathbf{\hat{J}_3}$, $\hat{\mathcal{E}}_{A}$ will be locally asymptotically stable if $\mathcal{R}_{0,B}<\exp(k\hat{I}_A)$, and thus $\mathcal{R}_{0,A}>\mathcal{R}_{0,B}$ (by Theorem \ref{thm:equilibria}, $\mathcal{R}_{0,A}>\exp(k\hat{I}_A)$ when $\hat{I}_A$ exists). From $ \mathbf{J_1}$ and $ \mathbf{\bar{J}_3}$, $\bar{\mathcal{E}}_{A}$ will be locally asymptotically stable if $\mathcal{R}_{0,B}<\dfrac{\exp(k\bar{I}_A)}{\exp(k\bar{I}_A)-s(\exp(k\bar{I}_A)-1)}$, and thus $\mathcal{R}_{0,A}>\mathcal{R}_{0,B}$ (by Theorem \ref{thm:equilibria}, $\mathcal{R}_{0,A}>\exp(k\hat{I}_A)$ when $\bar{I}_A$ exists).

A similar argument shows that, assuming $\mathcal{R}_{0,B}>1$, $\mathcal{E}_{B}$ will be locally asymptotically stable if $\mathcal{R}_{0,A}<1$; $\hat{\mathcal{E}}_{B}$ will be locally asymptotically stable if $\mathcal{R}_{0,A}<\exp(k\hat{I}_B)$, with $\mathcal{R}_{0,B}>\mathcal{R}_{0,A}$; and $\bar{\mathcal{E}}_{B}$ will be locally asymptotically stable if $\mathcal{R}_{0,A}<\dfrac{\exp(k\bar{I}_B)}{\exp(k\bar{I}_B)-s(\exp(k\bar{I}_B)-1)}$, with $\mathcal{R}_{0,B}>\mathcal{R}_{0,A}$. \\
~\\

\noindent \textbf{Proof of Theorem \ref{thm:EEstability}}\\\\

Assume $\mathcal{R}_{0,A}>1$ and $\mathcal{R}_{0,B}>1$, so $\mathcal{E}^*_{A,B}$ exists. Note that this will result in \[\beta_{0,A}e^{-kI^*_{A_E}}S^*_{A_E}-\frac{1}{\tau_I}=0\]
and 
\[\beta_{0,B}e^{-kI^*_{B_E}}S^*_{B_E}-\frac{1}{\tau_I}=0\]
(where $S^*_{A_E}=1-I^*_{A_E}-\frac{\tau_R}{\tau_I}I^*_{A_E}$ and $S^*_{B_E}=1-I^*_{B_E}-\frac{\tau_R}{\tau_I}I^*_{B_E}$). The Jacobian matrix evaluated at $\mathcal{E}^*_{A,B}$ is given below.
\begin{align*}
 &\mathbf{J(\mathcal{E}^*_{A,B}):=J^*_{A,B}}
&=\begin{pmatrix}
  \begin{matrix}
-\beta_{0,A}e^{-kI_{A_E}^*}I_{A_E}^*-\frac{1}{\tau_R} & -\frac{1}{\tau_I}-\frac{1}{\tau_R} & \frac{kI_{A_E}^*}{\tau_I} \\
    \beta_{0,A}e^{-kI_{A_E}^*}I_{A_E}^* & 0 & -\frac{kI_{A_E}^*}{\tau_I} \\
    0 & \frac{1}{\tau_P} & -\frac{1}{\tau_P}
  \end{matrix}
  & \rvline & \bigzero \\
\hline
  \bigzero & \rvline &
    \begin{matrix}
  -\beta_{0,B}e^{-kI_{B_E}^*}I_{B_E}^*-\frac{1}{\tau_R} & -\frac{1}{\tau_I}-\frac{1}{\tau_R} & \frac{kI_{B_E}^*}{\tau_I} \\
    \beta_{0,B}e^{-kI_{B_E}^*}I_{B_E}^* & 0 & -\frac{kI_{B_E}^*}{\tau_I} \\
    0 & \frac{1}{\tau_P} & -\frac{1}{\tau_P}
  \end{matrix}
\end{pmatrix} 
\end{align*}
It is straightforward to see that applying the Routh-Hurwitz criteria to the characteristic polynomial of each of the block diagonal matrices in $J^*_{A,B}$, in the same manner as in Theorem \ref{thm:BE} results in its six eigenvalues having negative real part. Thus, if $\mathcal{R}_{0,A}>1$ and $\mathcal{R}_{0,B}>1$, then the endemic equilibrium for the scenario with no spillover will be locally asymptotically stable.

\newpage

\subsection{Practical Identifiability Method} \label{Appendix:PI}

The steps below describe the Monte Carlo approach we used to determine the practical identifiability of the model parameters. 
\begin{itemize}
    \item[1)] Solve the system of ordinary differential equations using the true parameter vector $\hat{p}$ to obtain an output vector $g(x(t),\hat{p})$ at discrete time points $\{ t_{i}\}_{i=1}^{n}$.
    \item[2)] Construct $N=1,000$ datasets using an assigned measurement error. In other words, the data are described by $$y_{i}=g(x(t_{i}),\hat{p})+g(x(t_{i}),\hat{p})\epsilon, $$ where $\epsilon \sim  \mathcal{N}(0, \sigma)$ beginning with $\sigma=0$.
    \item[3)] Estimate the parameter set $p_{j}$ with respect to the generated datasets of prevalence and recognized prevalence. The minimizing function is the sum of squared errors between the model output and the datasets constructed in Step 2. We used the function \emph{fminsearchbnd} in MATLAB to estimate parameters. 
    \item[4)] After the optimization procedure, we computed the ARE values via 
        \begin{equation}
       \mathrm{ARE} \big(p^{(k)}\big)=100 \%  \times \frac{1}{N}  \sum_{j=1}^{N} \frac{\left| p^{(k)}- p_j^{(k)}\right|}{\left| p^{(k)}\right| },
       \label{eqn:ARE_eqn}
    \end{equation}
\item[5)] After the procedure is completed for $\sigma=0$, the process is repeated by increasing the measurement error by $\sigma=\{1\%, 5\%, 10\%, 20\%, 30\%\}.$
\end{itemize}

\subsection{Practical Identifiability Results}

\begin{table}[h!]
\centering
\caption{MC approach results for prevalence data with no spillover ($s=0$). The values in the table represent the average relative error (ARE) for each parameter with respect to the noise level added to the data. }\label{tab:Practical_Prevalence_None}
\begin{tabular}{|c||c|c|c|c|c|c||}
\hline
\multicolumn{1}{|c||}{} & \multicolumn{6}{|c||}{\textbf{No Spillover ($s=0$})}  \\ 
\hline $\sigma$ & $\beta_A$	& $\beta_B$ & $k$ & $\tau_R$ & $\tau_I$ & $\tau_P$  \\
\hline 0\% & 0 & 0 & 0 & 0 & 0 & 0 \\
\hline 1\% & 0.2763 & 0.4924 & 1.0918 & 2.5759 & 1.0455 & 0.1852 \\
\hline 5\% & 1.2771 & 2.2904 & 5.039 & 16.5385 & 4.7842 & 1.0242 \\
\hline 10\% & 2.8153 & 4.8763 & 10.3588 & 47.0884 & 9.6123 & 2.2142  \\
\hline 20\% & 7.3598 & 12.3883 & 18.7359 & 130.5298 & 18.3179 & 6.0153  \\
\hline 30\% & 14.852 & 23.7868 & 27.2016 & 241.342 & 26.3271 & 11.4872  \\
\hline Identifiable? & Yes & Yes & Weakly & Weakly & Weakly & Yes  \\
\hline
\end{tabular}
\end{table}

\begin{table}[h!]
\centering
\caption{MC approach results for prevalence data with imperfect spillover where $s=0.1$. The values in the table represent the average relative error (ARE) for each parameter with respect to the noise level added to the data. }\label{tab:Practical_Prevalence_Imperfect_0.1}
\begin{tabular}{|c||c|c|c|c|c|c||}
\hline
\multicolumn{1}{|c||}{} & \multicolumn{6}{|c||}{\textbf{Imperfect Spillover $s=0.1$}} \\ 
\hline $\sigma$  & $\beta_A$	& $\beta_B$ & $k$ & $\tau_R$	& $\tau_I$ & $\tau_P$	 \\
\hline 0\% & 0 & 0 & 0 & 0 & 0 & 0 \\
\hline 1\% & 0.2864 & 0.5145 & 1.1715 & 2.8399 & 1.0773 & 0.1823 \\
\hline 5\% & 1.3845 & 2.447 & 5.6331 & 19.0718 & 5.1352 & 0.9918 \\
\hline 10\% & 2.9838 & 5.2547 & 11.1374 & 59.2036 & 10.2647 & 2.3034 \\
\hline 20\% & 7.0616 & 11.8189 & 19.0378 & 104.0596 & 18.2241 & 5.0284 \\
\hline 30\% & 14.3226 & 23.1563 & 27.0026 & 254.4839 & 26.4345 & 11.1761 \\
\hline Identifiable? & Yes & Yes & Weakly & Weakly & Weakly & Yes \\
\hline
\end{tabular}
\end{table}

\begin{table}[h!]
\centering
\caption{MC approach results for prevalence data with imperfect spillover where $s=0.5$. The values in the table represent the average relative error (ARE) for each parameter with respect to the noise level added to the data. }\label{tab:Practical_Prevalence_Imperfect_0.5}
\begin{tabular}{|c||c|c|c|c|c|c||}
\hline
\multicolumn{1}{|c||}{} & \multicolumn{6}{|c||}{\textbf{Imperfect Spillover $s=0.5$}} \\ 
\hline $\sigma$  & $\beta_A$	& $\beta_B$ & $k$ & $\tau_R$	& $\tau_I$ & $\tau_P$	 \\
\hline 0\% & 0 & 0 & 0 & 0 & 0 & 0 \\
\hline 1\% & 0.2845 & 0.4787 & 1.1067 & 2.9084 & 1.0886 & 0.2145 \\
\hline 5\% & 1.4744 & 2.4773 & 6.0112 & 27.0649 & 5.6221 & 1.2729 \\
\hline 10\% & 3.149 & 5.054 & 10.8949 & 67.8649 & 10.7817 & 2.7958 \\
\hline 20\% & 9.7035 & 14.1726 & 19.9611 & 253.7215 & 20.8321 & 8.2631 \\
\hline 30\% & 16.9918 & 23.5491 & 27.0075 & 378.8752 & 27.8657 & 14.0836 \\
\hline Identifiable? & Yes & Yes & Weakly & No & Weakly & Yes \\
\hline
\end{tabular}
\end{table}

\begin{table}[h!]
\centering
\caption{MC approach results for prevalence data with imperfect spillover where $s=0.9$. The values in the table represent the average relative error (ARE) for each parameter with respect to the noise level added to the data. }\label{tab:Practical_Prevalence_Imperfect_0.9}
\begin{tabular}{|c||c|c|c|c|c|c||}
\hline
\multicolumn{1}{|c||}{} & \multicolumn{6}{|c||}{\textbf{Imperfect Spillover $s=0.9$}} \\ 
\hline $\sigma$  & $\beta_A$	& $\beta_B$ & $k$ & $\tau_R$	& $\tau_I$ & $\tau_P$	 \\
\hline 0\% & 0 & 0 & 0 & 0 & 0 & 0 \\
\hline 1\% & 0.2953 & 0.5157 & 1.1369 & 3.0146 & 1.1028 & 0.2297 \\
\hline 5\% & 1.5602 & 2.675 & 5.9749 & 21.761 & 5.62 & 1.2048 \\
\hline 10\% & 3.1565 & 5.1839 & 11.0611 & 62.4287 & 10.8384 & 2.6086 \\
\hline 20\% & 9.3393 & 13.6041 & 19.15 & 245.6167 & 20.3468 & 8.1421 \\
\hline 30\% & 17.4082 & 22.9349 & 27.2479 & 393.7193 & 28.2785 & 14.3078 \\
\hline Identifiable? & Yes & Yes & Weakly & No & Weakly & Yes \\
\hline
\end{tabular}
\end{table}

\begin{table}[h!]
\centering
\caption{MC approach results for prevalence data with perfect spillover ($s=1$). The values in the table represent the average relative error (ARE) for each parameter with respect to the noise level added to the data.  }\label{tab:Practical_Prevalence_Perfect}
\begin{tabular}{|c||c|c|c|c|c|c||}
\hline
\multicolumn{1}{|c||}{}  & \multicolumn{6}{|c||}{\textbf{Perfect Spillover ($s=1$})} \\ 
\hline $\sigma$  & $\beta_A$	& $\beta_B$ & $k$ & $\tau_R$ & $\tau_I$ & $\tau_P$ \\
\hline 0\% & 0 & 0 & 0 & 0 & 0 & 0 \\
\hline 1\% & 0.3144 & 0.5796 & 1.1973 & 3.1499 & 1.1672 & 0.1914 \\
\hline 5\% & 1.5014 & 2.7624 & 5.6142 & 17.2976 & 5.4727 & 0.9932 \\
\hline 10\% & 3.1466 & 5.5832 & 10.8647 & 39.5551 & 10.418 & 2.025 \\
\hline 20\% & 7.0348 & 11.8155 & 18.4124 & 94.4371 & 18.1673 & 5.0349 \\
\hline 30\% & 13.9267 & 21.9621 & 26.8615 & 202.176 & 26.1597 & 9.1224 \\
\hline Identifiable? & Yes & Yes & Weakly & Weakly & Weakly & Yes \\
\hline
\end{tabular}
\end{table}

\begin{table}[h!]
\centering
\caption{MC approach results for recognized prevalence data with no spillover ($s=0$). The values in the table represent the average relative error (ARE) for each parameter with respect to the noise level added to the data. }\label{tab:Practical_Perceived_Prevalence_None}
\begin{tabular}{|c||c|c|c|c|c|c|c|c||}
\hline
\multicolumn{1}{|c||}{} & \multicolumn{8}{|c||}{\textbf{No spillover ($s=0$)}}  \\ 
\hline $\sigma$  & $\beta_A$	& $\beta_B$ & $k$ & $\tau_R$ & $\tau_I$ & $\tau_P$	 & $K_A$ & 	$K_B$  \\
\hline 0\% & 0 & 0 & 0 & 0 & 0 & 0 & 0 & 0    \\
\hline 1\% & 0.5573 & 0.9723 & 2.3101 & 2.4096 & 2.4066 & 0.5471 & 3.2341 & 3.5235  \\
\hline 5\% & 2.8243 & 4.7754 & 7.8151 & 14.1952 & 10.459 & 2.6675 & 10.1982 & 11.3486  \\
\hline 10\% & 6.8513 & 10.9262 & 15.5526 & 65.5675 & 20.5047 & 5.8883 & 17.7665 & 20.12  \\
\hline 20\% & 19.3635 & 28.8108 & 29.8194 & 321.0108 & 40.2069 & 14.7349 & 29.6125 & 33.7876  \\
\hline 30\% & 29.23 & 44.8878 & 38.6998 & 374.87 & 48.8166 & 24.1475 & 33.3833 & 38.1048  \\
\hline Identifiable? & Yes & Weakly & Weakly & Weakly & Weakly & Yes & Weakly & Weakly  \\
\hline
\end{tabular} 
\end{table}

\begin{table}[h!]
\centering
\caption{MC approach results for recognized prevalence data with imperfect spillover where $s=0.1$. The values in the table represent the average relative error (ARE) for each parameter with respect to the noise level added to the data. }\label{tab:Practical_Perceived_Prevalence_Imperfect_0.1}
\begin{tabular}{|c||c|c|c|c|c|c|c|c||}
\hline
\multicolumn{1}{|c||}{} & \multicolumn{8}{|c||}{\textbf{Imperfect Spillover $s=0.1$}} \\ 
\hline $\sigma$  & $\beta_A$	& $\beta_B$ & $k$ & $\tau_R$	& $\tau_I$ & $\tau_P$	 & $K_A$ & 	$K_B$ \\
\hline 0\% & 0 & 0 & 0 & 0 & 0 & 0 & 0 & 0 \\
\hline 1\% & 0.5742 & 0.9721 & 2.374 & 2.5283 & 2.3035 & 0.5562 & 2.8904 & 3.1886 \\
\hline 5\% & 2.7926 & 4.7816 & 7.8388 & 15.849 & 10.3548 & 2.7315 & 9.2464 & 10.583 \\
\hline 10\% & 6.5938 & 10.6576 & 14.2412 & 76.3384 & 20.1283 & 5.9298 & 14.2889 & 16.6075 \\
\hline 20\% & 17.6301 & 26.5648 & 28.4594 & 323.3021 & 40.0006 & 14.4302 & 26.4062 & 30.68 \\
\hline 30\% & 26.6671 & 40.9656 & 36.5904 & 399.4829 & 46.0969 & 22.8819 & 30.2587 & 35.1408 \\
\hline Identifiable? & Yes & Weakly & Weakly & No & Weakly & Yes & Weakly & Weakly \\
\hline
\end{tabular}
\end{table}

\begin{table}[h!]
\centering
\caption{MC approach results for recognized prevalence data with imperfect spillover where $s=0.5$. The values in the table represent the average relative error (ARE) for each parameter with respect to the noise level added to the data. }\label{tab:Practical_Perceived_Prevalence_Imperfect_0.5}
\begin{tabular}{|c||c|c|c|c|c|c|c|c||}
\hline
\multicolumn{1}{|c||}{} & \multicolumn{8}{|c||}{\textbf{Imperfect Spillover $s=0.5$}} \\ 
\hline $\sigma$  &  $\beta_A$	& $\beta_B$ & $k$ & $\tau_R$	& $\tau_I$ & $\tau_P$	 & $K_A$ & 	$K_B$ \\
\hline 0\% & 0 & 0 & 0 & 0 & 0 & 0 & 0 & 0 \\
\hline 1\% & 0.5356 & 0.9212 & 2.4137 & 3.1561 & 1.7311 & 0.4895 & 1.8389 & 2.4389 \\
\hline 5\% & 2.4712 & 4.262 & 8.2662 & 31.8741 & 8.1389 & 2.7837 & 6.2542 & 8.3069 \\
\hline 10\% & 5.6866 & 9.3309 & 14.1045 & 108.2123 & 15.5261 & 6.0023 & 11.6017 & 15.3663 \\
\hline 20\% & 16.3214 & 25.2412 & 24.8117 & 269.8939 & 31.812 & 14.3797 & 22.9048 & 32.395 \\
\hline 30\% & 23.4349 & 35.1604 & 32.9777 & 388.8011 & 40.1668 & 23.4831 & 29.4802 & 41.1973 \\
\hline Identifiable? & Yes & Weakly & Weakly & No & Weakly & Yes & Weakly & Weakly \\
\hline
\end{tabular}
\end{table}

\begin{table}[h!]
\centering
\caption{MC approach results for recognized prevalence data with imperfect spillover where $s=0.9$. The values in the table represent the average relative error (ARE) for each parameter with respect to the noise level added to the data. }\label{tab:Practical_Perceived_Prevalence_Imperfect_0.9}
\begin{tabular}{|c||c|c|c|c|c|c|c|c||}
\hline
\multicolumn{1}{|c||}{} & \multicolumn{8}{|c||}{\textbf{Imperfect Spillover $s=0.9$}} \\ 
\hline $\sigma$  & $\beta_A$	& $\beta_B$ & $k$ & $\tau_R$	& $\tau_I$ & $\tau_P$	 & $K_A$ & 	$K_B$ \\
\hline 0\% & 0 & 0 & 0 & 0 & 0 & 0 & 0 & 0 \\
\hline 1\% & 0.6688 & 1.6347 & 3.5109 & 4.6411 & 2.9101 & 0.6689 & 4.6772 & 3.5536 \\
\hline 5\% & 3.4076 & 6.5134 & 12.3635 & 41.3099 & 12.104 & 3.6436 & 14.3734 & 27.1759 \\
\hline 10\% & 7.8052 & 12.3904 & 20.0273 & 129.5587 & 21.2405 & 7.9192 & 20.7891 & 41.2142 \\
\hline 20\% & 17.4125 & 27.2516 & 30.9329 & 346.1011 & 38.6858 & 15.4679 & 34.1449 & 51.6439 \\
\hline 30\% & 25.4263 & 36.785 & 37.5017 & 420.8189 & 45.7195 & 24.1191 & 37.775 & 59.4871 \\
\hline Identifiable? & Yes & Weakly & Weakly & No & Weakly & Yes & Weakly & Weakly \\
\hline
\end{tabular}
\end{table}

\begin{table}[h!]
\centering
\caption{MC approach results for recognized prevalence data with perfect spillover ($s=1$). The values in the table represent the average relative error (ARE) for each parameter with respect to the noise level added to the data. }\label{tab:Practical_Perceived_Prevalence_Perfect}
\begin{tabular}{|c||c|c|c|c|c|c|c|c||}
\hline
\multicolumn{1}{|c||}{} & \multicolumn{8}{|c||}{\textbf{Perfect Spillover ($s=1$})}  \\ 
\hline $\sigma$  & $\beta_A$	& $\beta_B$ & $k$ & $\tau_R$ & $\tau_I$ & $\tau_P$	 & $K_A$ & 	$K_B$ \\
\hline 0\% & 0 & 0 & 0 & 0 & 0 & 0 & 0 & 0  \\
\hline 1\% &  0.725 & 2.1585 & 4.4542 & 5.8929 & 3.941 & 0.9023 & 6.9895 & 12.8302 \\
\hline 5\% &  4.1666 & 10.2032 & 15.4967 & 44.8981 & 15.7003 & 4.0507 & 22.0269 & 49.5525 \\
\hline 10\% &  9.2879 & 17.5328 & 23.7134 & 109.7555 & 26.8304 & 7.7913 & 30.2851 & 59.8575 \\
\hline 20\% & 18.8488 & 30.7405 & 34.9604 & 265.5261 & 46.4454 & 16.0036 & 39.68 & 63.0998 \\
\hline 30\% & 25.989 & 40.3422 & 41.3798 & 399.2245 & 54.0062 & 23.6666 & 41.1436 & 66.5048 \\
\hline Identifiable? &  Yes & Weakly & Weakly & No & Weakly & Yes & Weakly & Weakly \\
\hline
\end{tabular} 
\end{table}

\clearpage
\section{Numerical results with fractional formulation of behavioral feedback}
\label{sec:app:frac}

To show the robustness of our results for behavioral spillover, we repeat our numerical simulations with another formulation for behavioral feedback.
Rather than the exponential formulation in Eq. \eqref{eq:e}, we use a fractional formulation as in \cite{rahmandad2022enhancing}, given by
\begin{align}
\label{eq:e_frac}
& m_i=\frac{1} {\left(1+\alpha \widetilde{I}_i\right)^{\gamma}}, 
\end{align}
where $\alpha = 100$ and $\gamma = 2$, consistent with analyses in \cite{rahmandad2022enhancing,lejeune2024mathematical}. We reproduce what is shown in Figures \ref{fig:dyanmics_spillover} and \ref{fig:persistence_exceed} using this formulation. Overall, the results are qualitatively identical, while the quantitative results differ slightly. In particular, there are two possible outcomes: co-existence of both diseases or exclusion of disease-$B$, which can happen when the basic reproductive number of disease-$B$ is above 1 as long as spillover is high enough. Quantitatively, the epidemic peaks do not reach as high a level (Figs. \ref{fig:dyanmics_spillover} and \ref{fig:dyanmics_spillover_frac}), which reduces the amount of difference between disease $A$ and disease $B$ (Figs. \ref{fig:persistence_exceed} and \ref{fig:persistence_exceed_frac}, panels c,d). Thus, the features we discuss in detail in the main text are the result of the behavioral spillover, rather than the precise formulation of behavioral feedback. We do not repeat the analytical work for the fractional formulation so all results presented are the outcome of numerical simulations.

\begin{figure}[h]
\centering

\begin{subfigure}[t]{.32\textwidth}
\centering    \includegraphics[width=\linewidth]{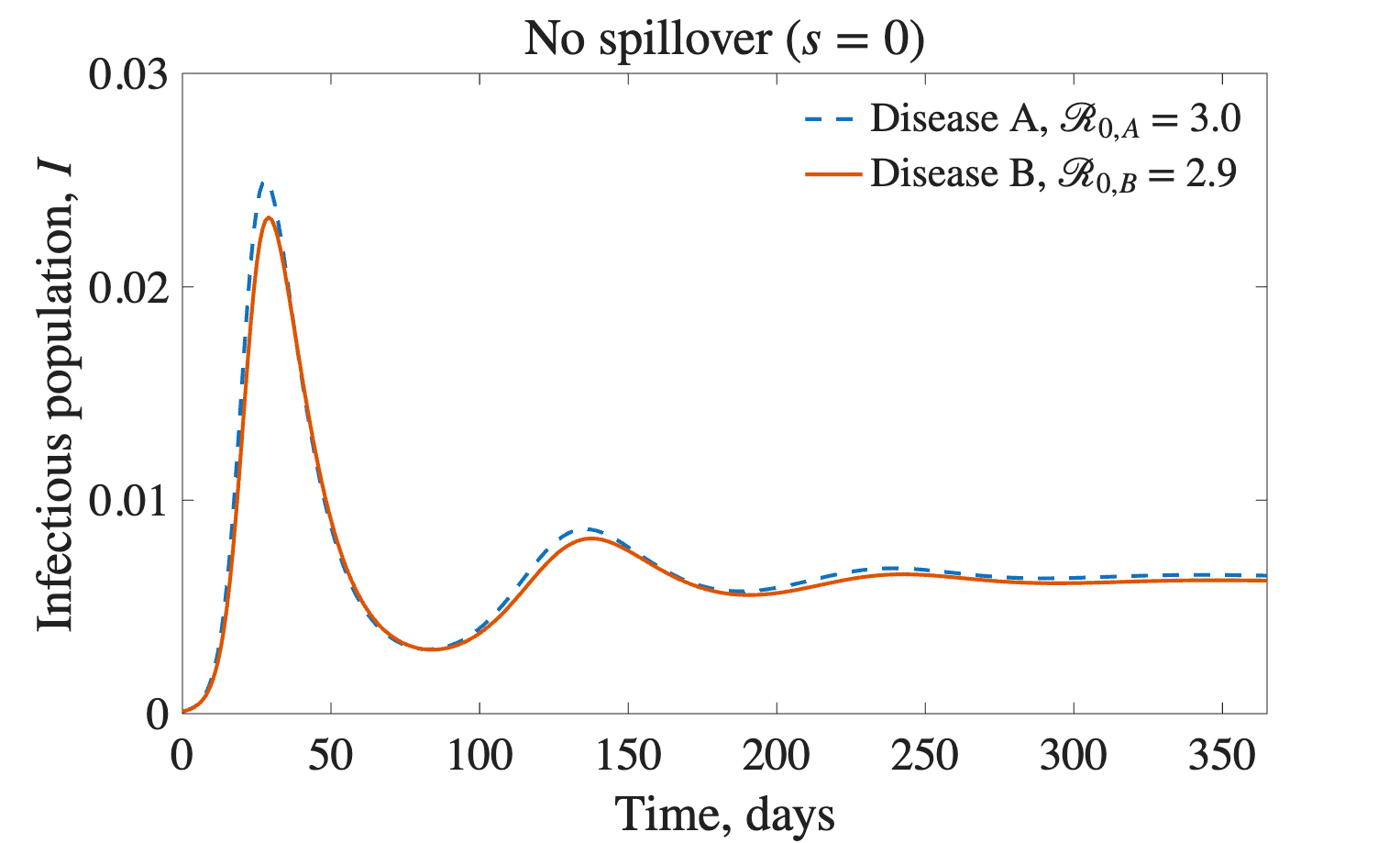}
\caption{}
\end{subfigure} 
\begin{subfigure}[t]{.32\textwidth}
\centering  
\includegraphics[width=\linewidth]{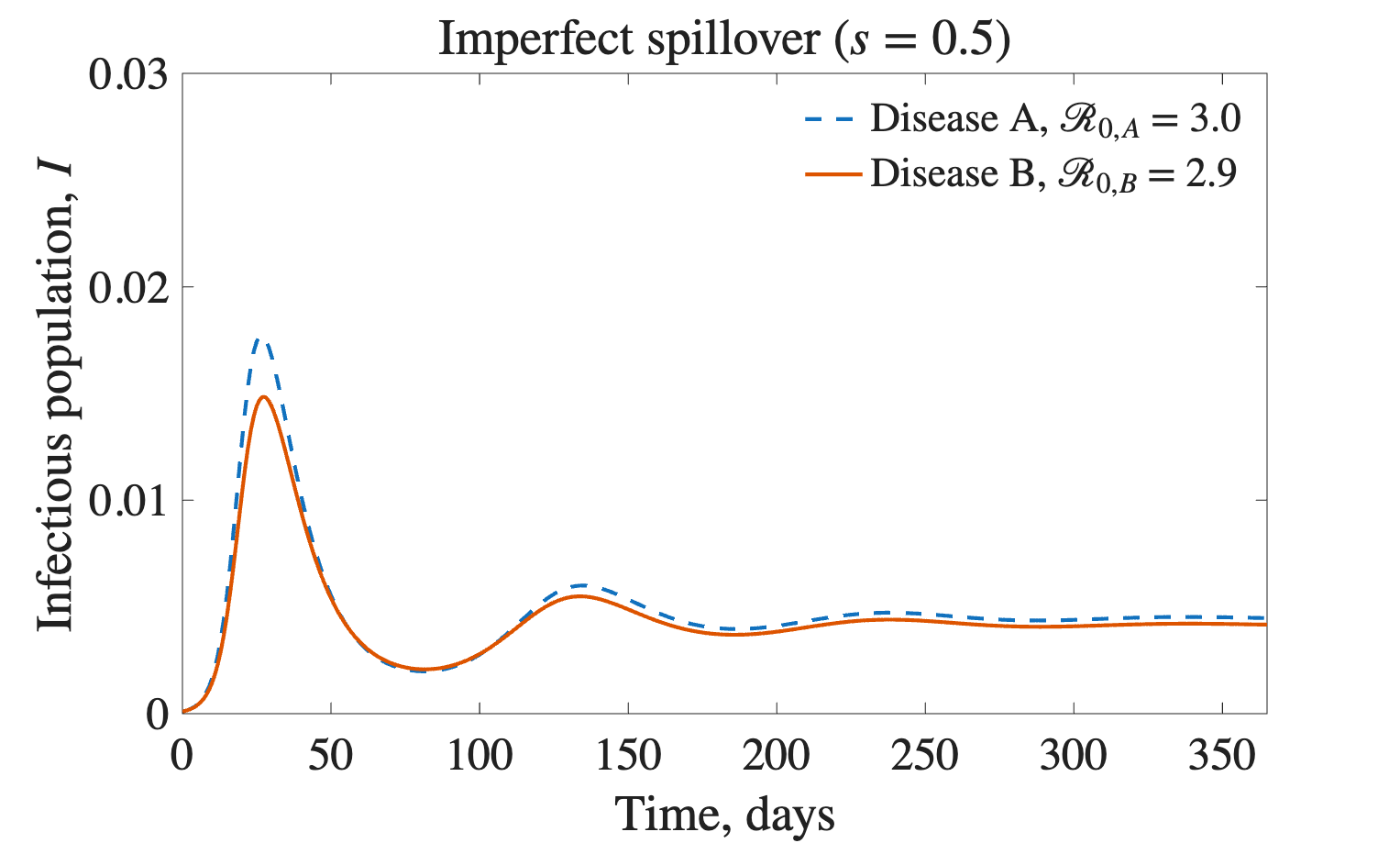}
\caption{}
\end{subfigure}
\begin{subfigure}[t]{.32\textwidth}
\centering \includegraphics[width=\linewidth]{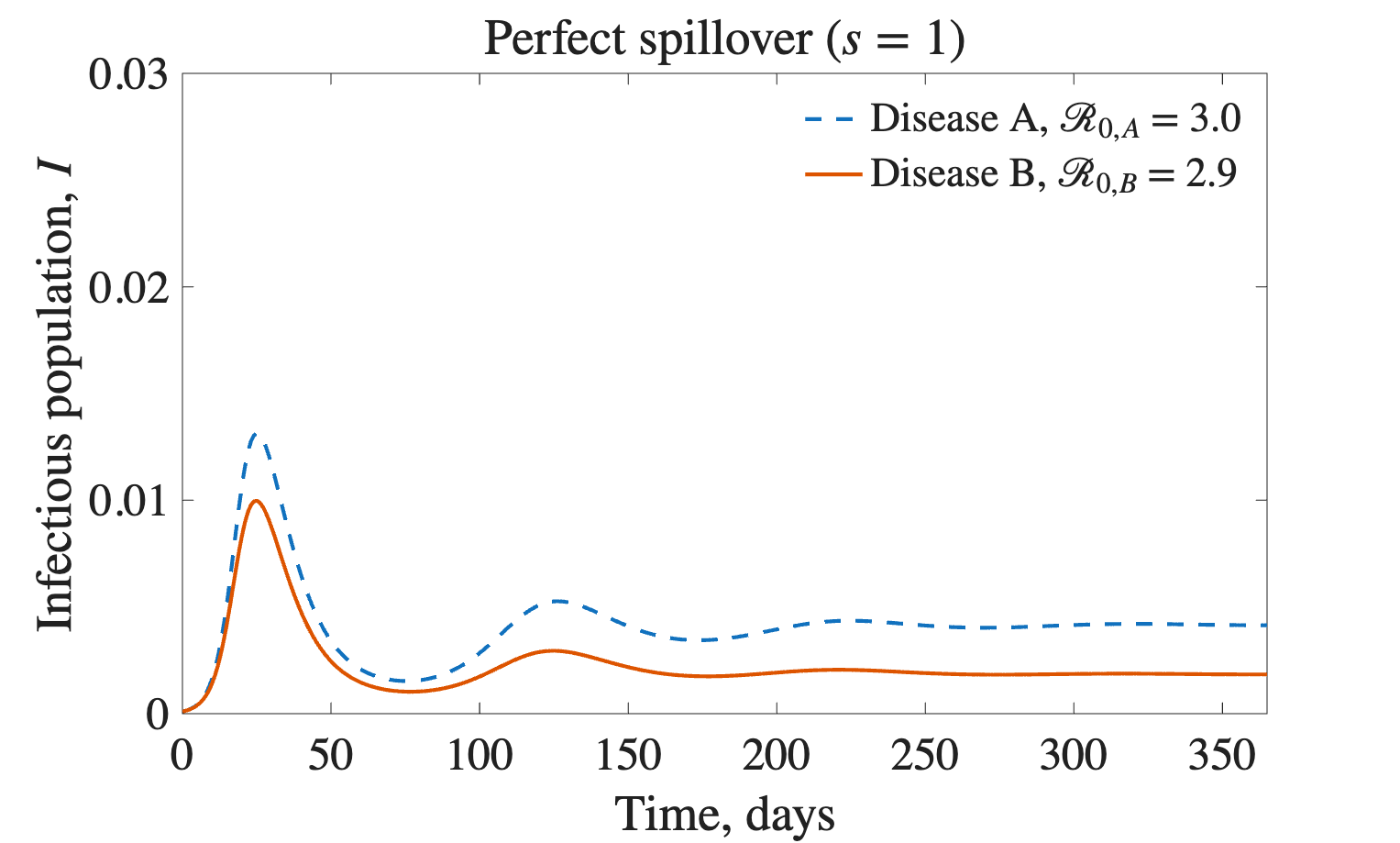}
\caption{}
\end{subfigure}\\

\begin{subfigure}[t]{.32\textwidth}
\centering \includegraphics[width=\linewidth]  {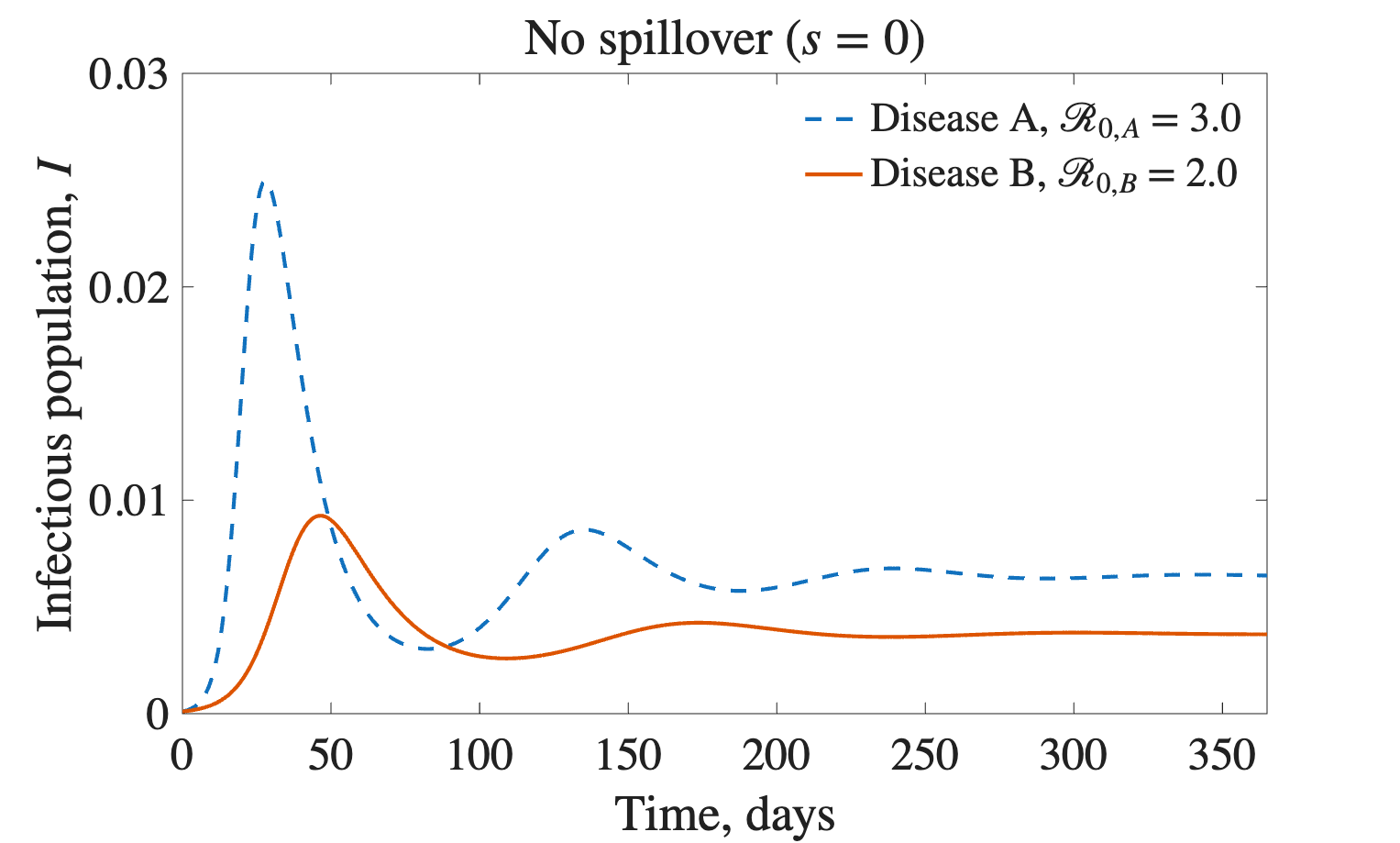}
\caption{}
\end{subfigure}
\begin{subfigure}[t]{.32\textwidth}
\centering \includegraphics[width=\linewidth]{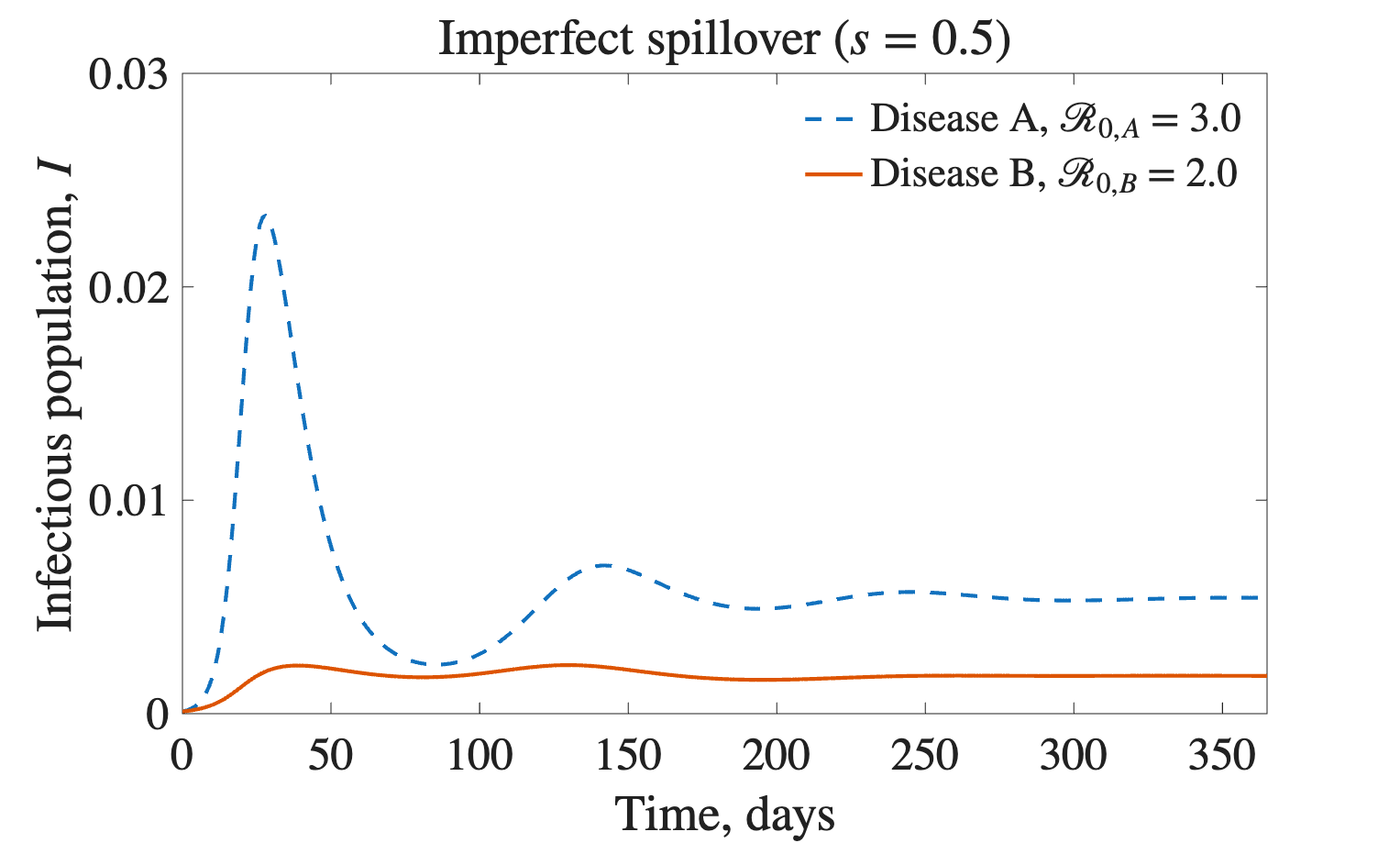}
\caption{}
\end{subfigure}
\begin{subfigure}[t]{.32\textwidth}
\centering \includegraphics[width=\linewidth]{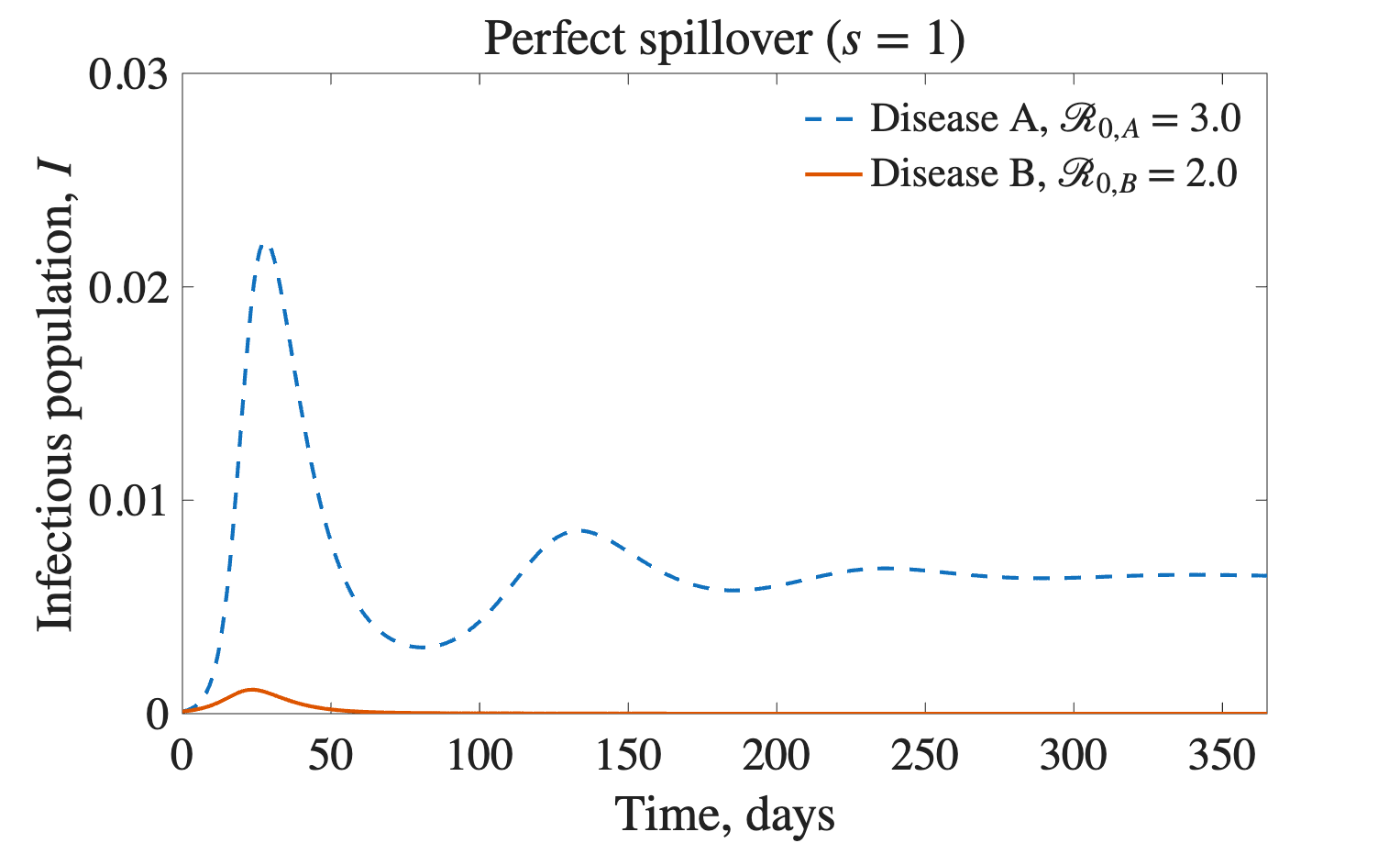}
\caption{}
\end{subfigure}\\

\begin{subfigure}[t]{.32\textwidth}
\centering \includegraphics[width=\linewidth]{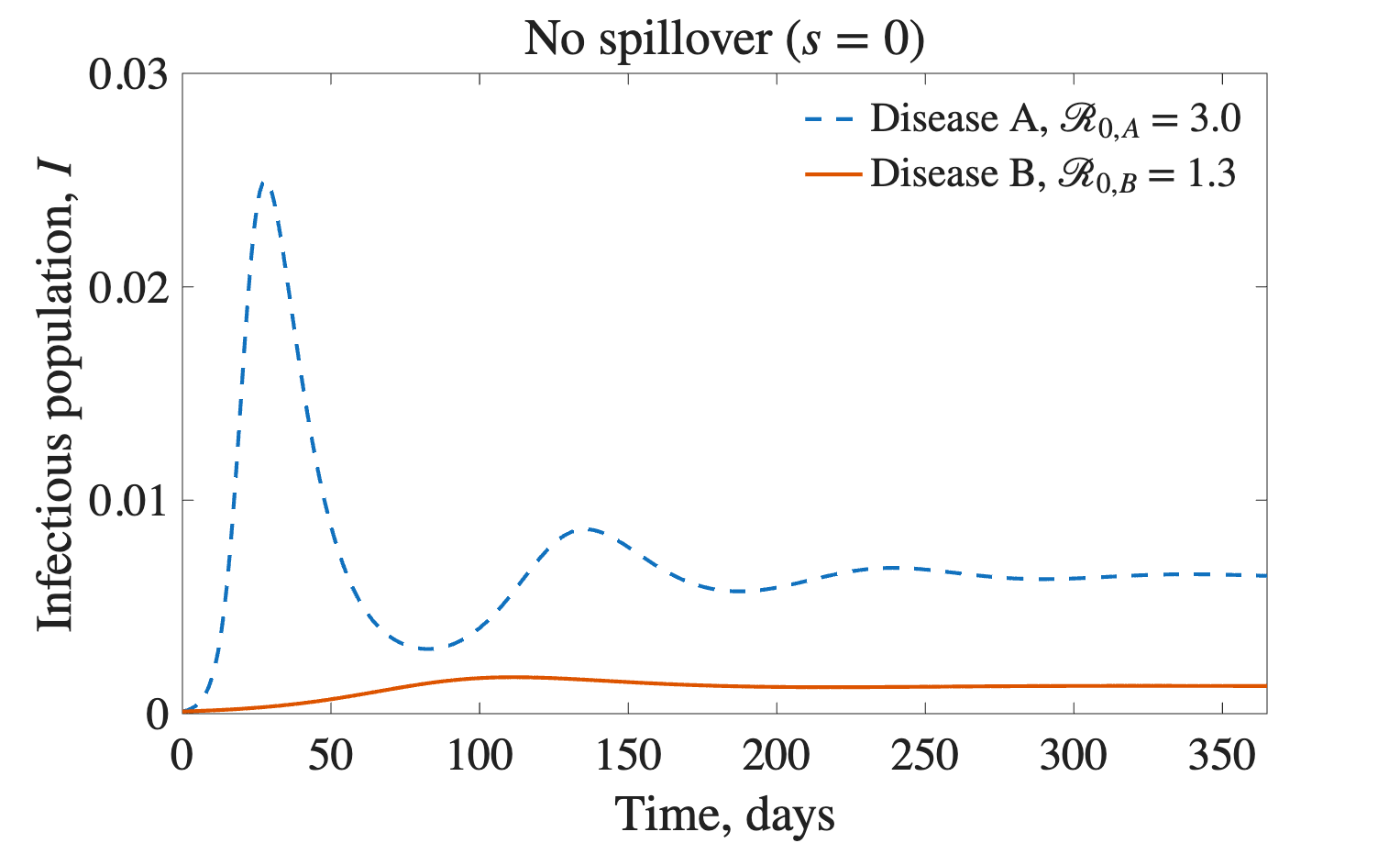}
\caption{}
\end{subfigure}
\begin{subfigure}[t]{.32\textwidth}
\centering \includegraphics[width=\linewidth]{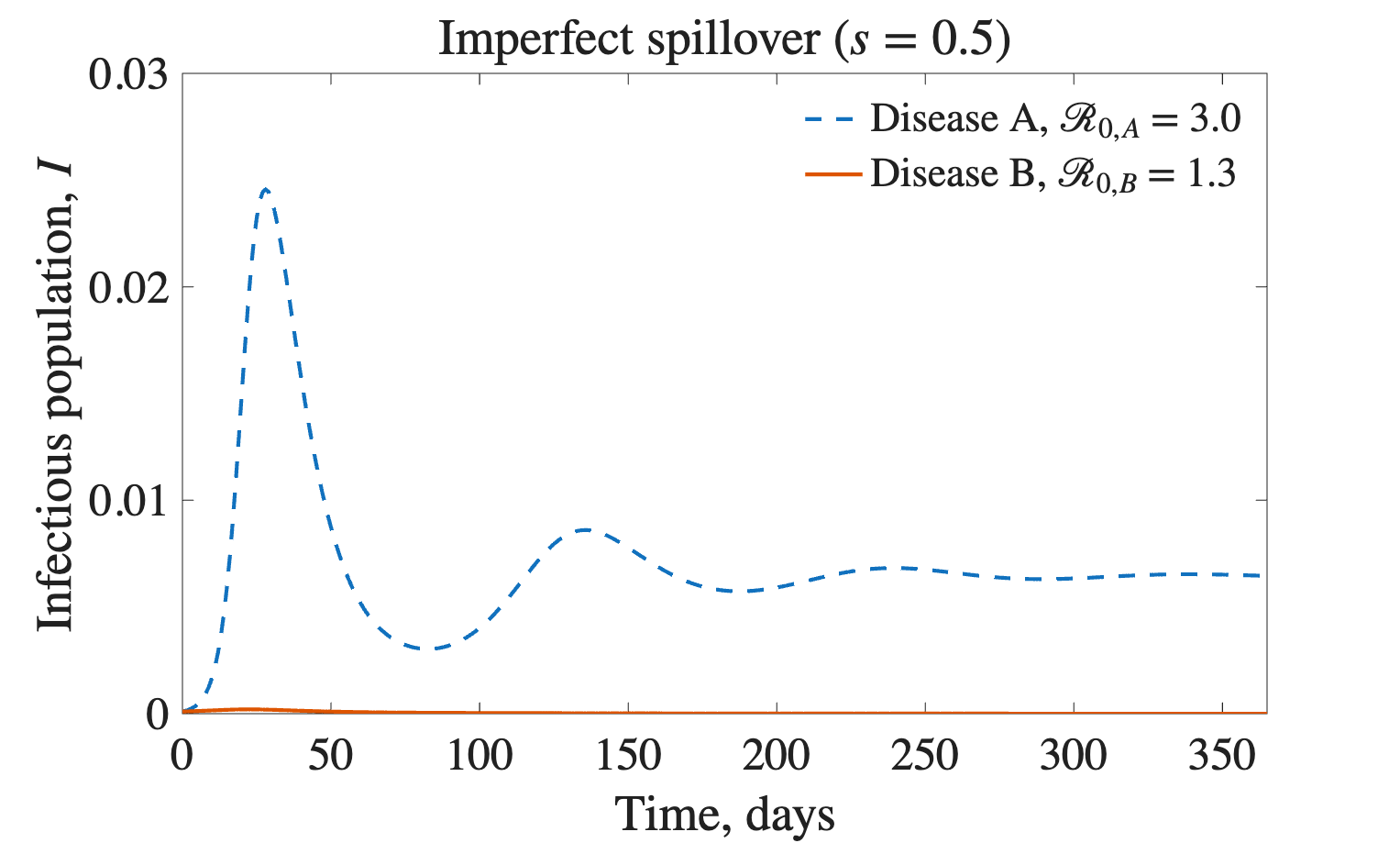}
\caption{}
\end{subfigure}
\begin{subfigure}[t]{.32\textwidth}
\centering \includegraphics[width=\linewidth]{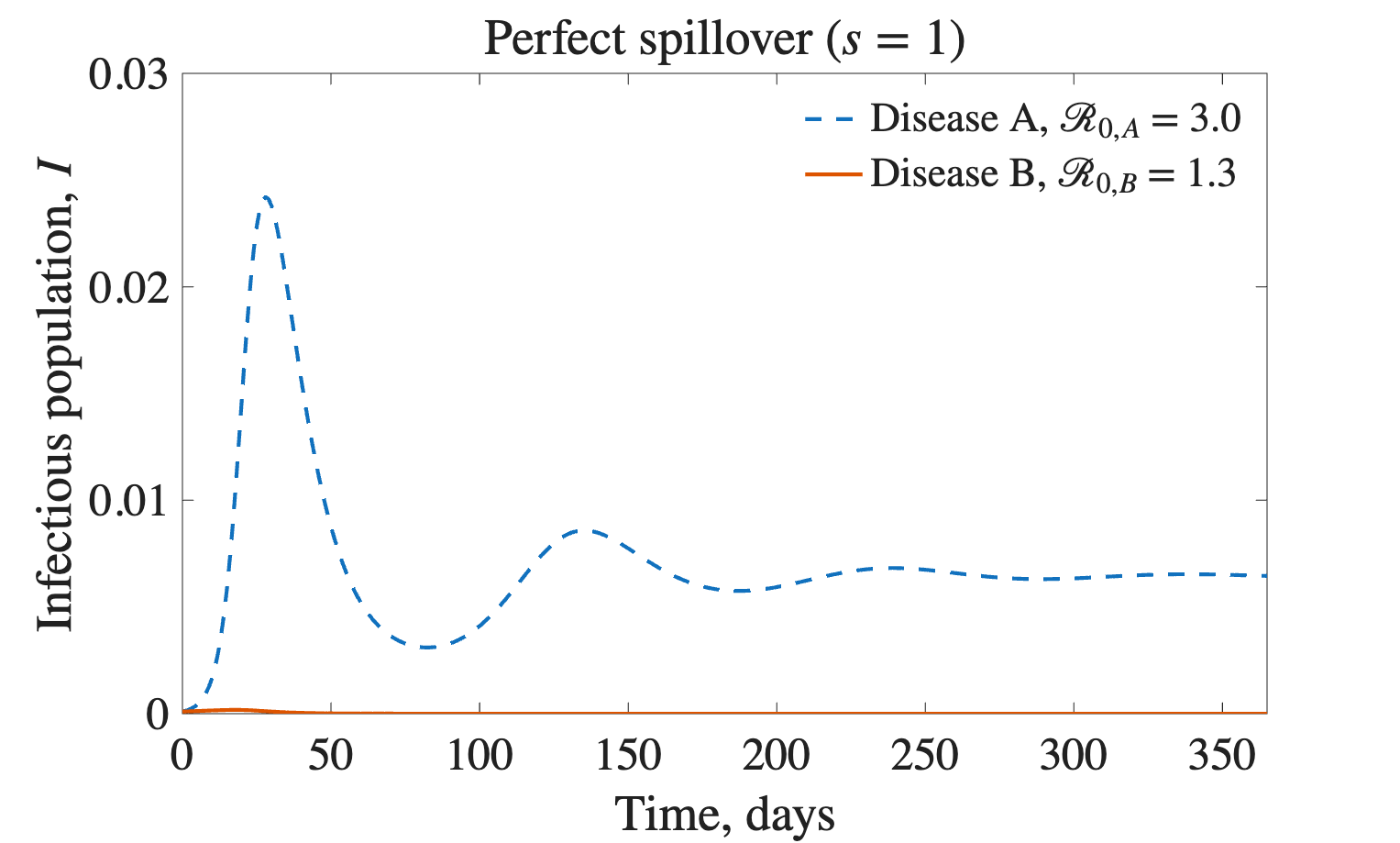}
\caption{}
\end{subfigure}\\

\caption{Dynamics with three levels of spillover and three values of the basic reproduction number of disease $B$ with a fractional formulation for behavioral feedback (as opposed to exponential formulation in Figure \ref{fig:dyanmics_spillover_frac}). For (a), (d) and (g), there is no spillover, $s=0$, which corresponds to the two diseases spreading independently (i.e., the curve is identical to that of the single disease spreading in the population); for (b), (e), (h), there is  imperfect spillover, $s=0.5$; and for (c), (f), (i), there is perfect spillover, $s=1$. For (a)-(c), $\mathcal{R}_{0,B}=2.9$; for (d)-(f), $\mathcal{R}_{0,B}=2$; and for (g)-(i) $\mathcal{R}_{0,B}=1.3$). Other parameters as in Table \ref{table:params} with $\mathcal{R}_{0,A}=3$ and $\alpha=100$ and $\gamma=2$, as discussed in Section \ref{sec:app:frac}.}
    \label{fig:dyanmics_spillover_frac}
\end{figure}

\begin{figure}[h]
    \centering
    \begin{subfigure}[t]{.49\textwidth}
\centering
    \includegraphics[width=\linewidth]{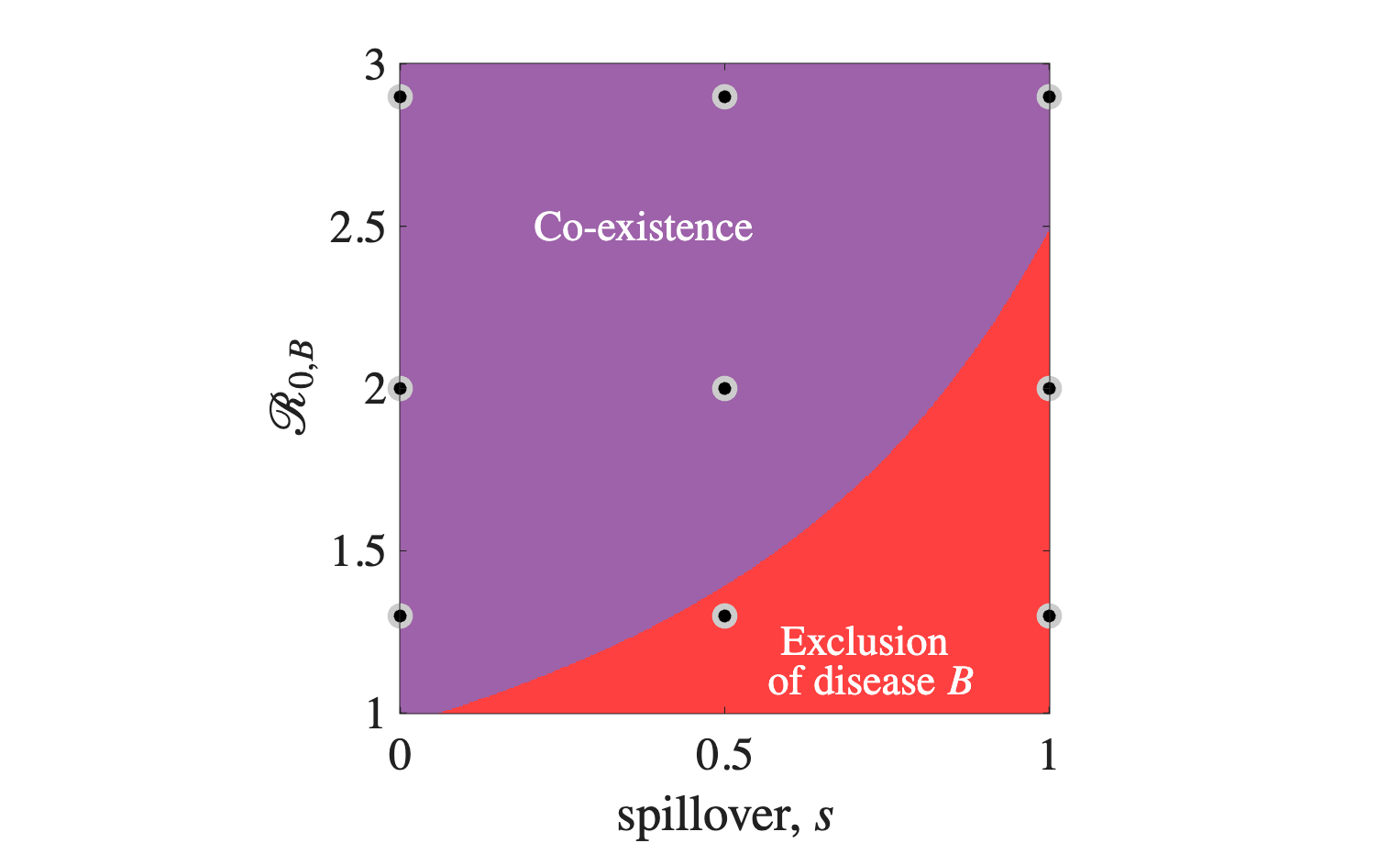}
    \caption{Co-existence vs. exclusion (simulation)}\label{fig:persistence_exceed_A_frac}
\end{subfigure}%
    \begin{subfigure}[t]{.49\textwidth}
\centering
    \includegraphics[width=\linewidth]{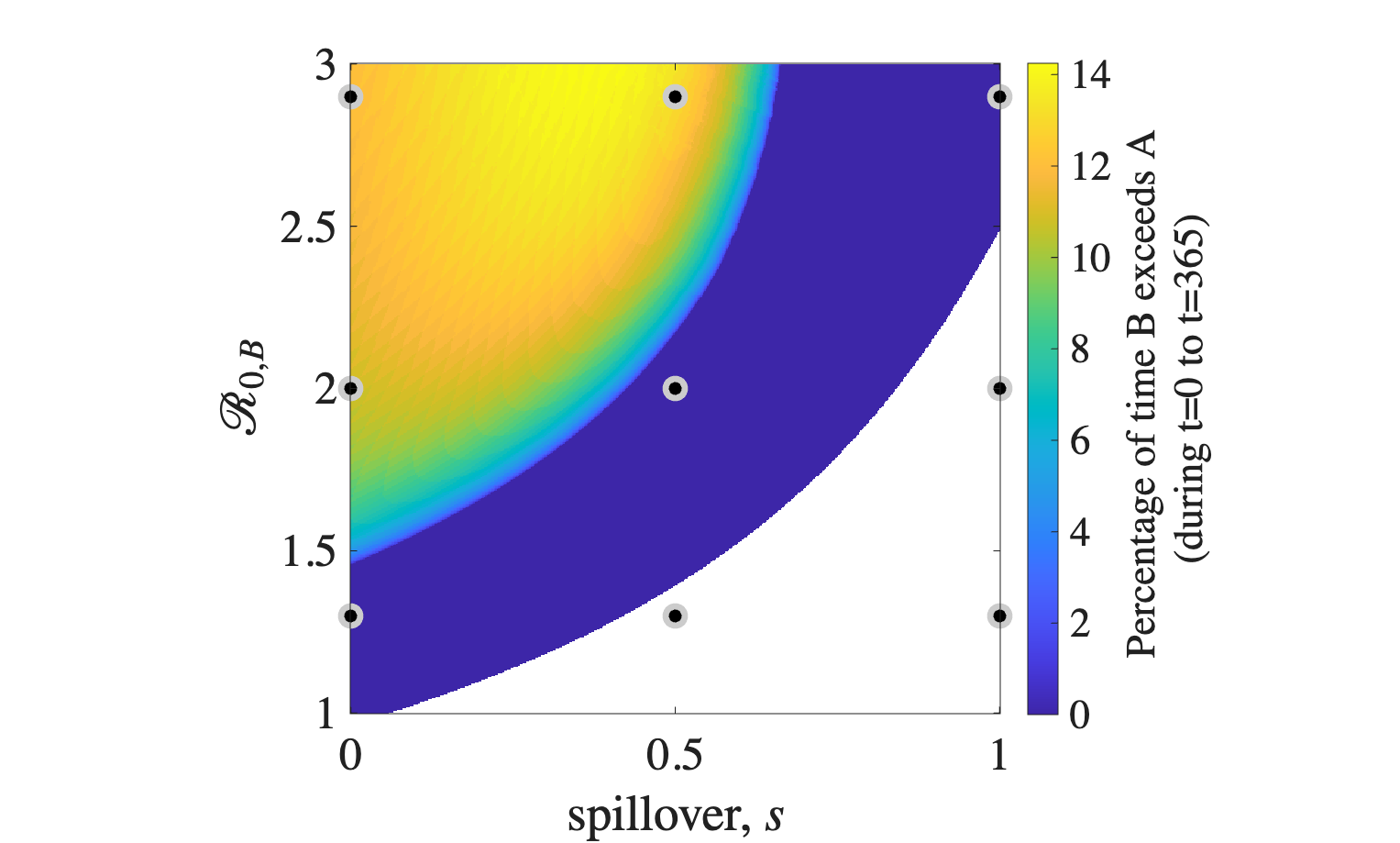}
    \caption{Time disease $B$ exceeds disease  $A$}\label{fig:persistence_exceed_B_frac}
\end{subfigure}
    \begin{subfigure}[t]{.49\textwidth}
\centering
    \includegraphics[width=\linewidth]{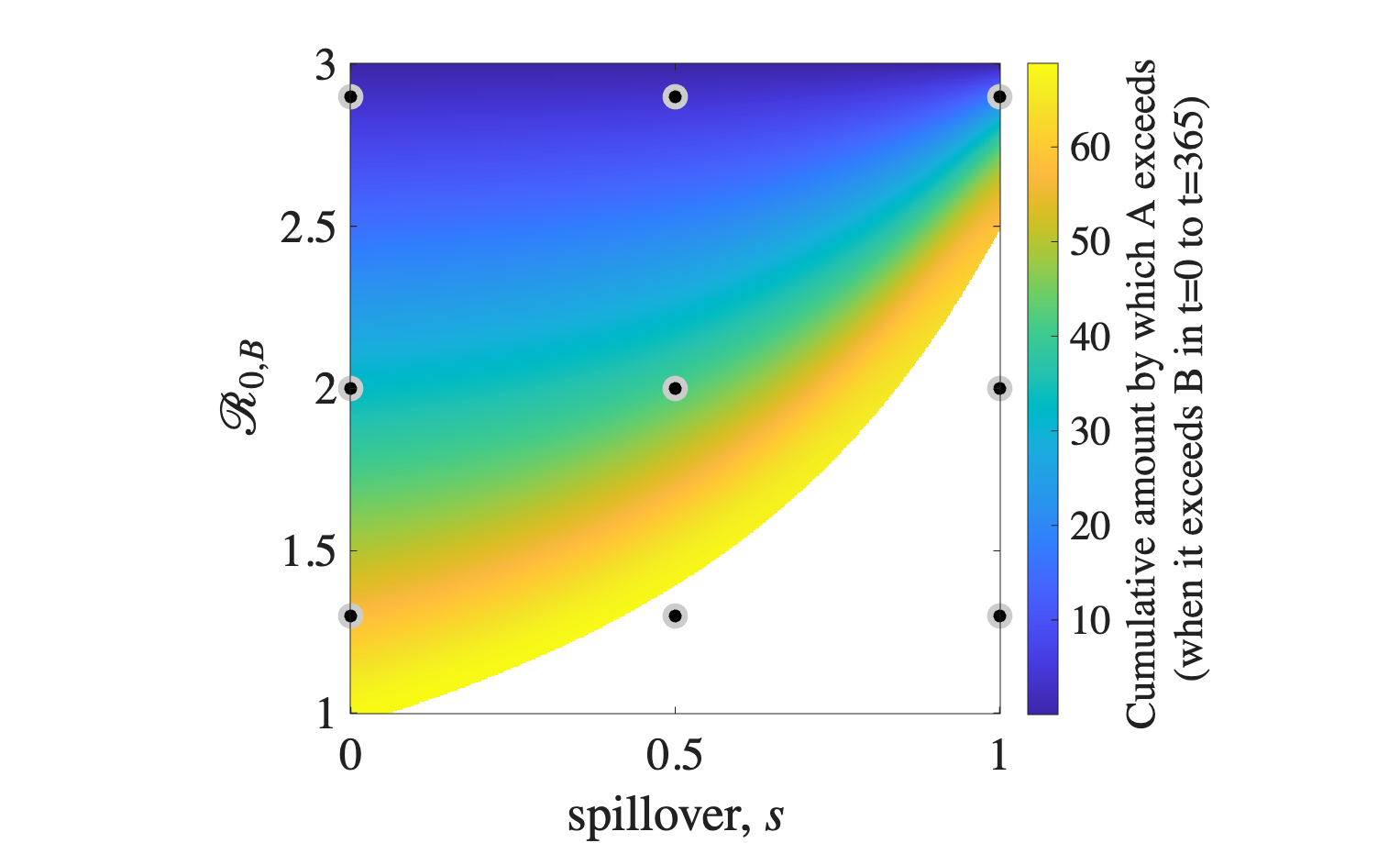}
    \caption{Amount disease $A$ exceeds disease  $B$}\label{fig:persistence_A_exceed_B_frac}
\end{subfigure}% 
    \begin{subfigure}[t]{.49\textwidth}
\centering
    \includegraphics[width=\linewidth]{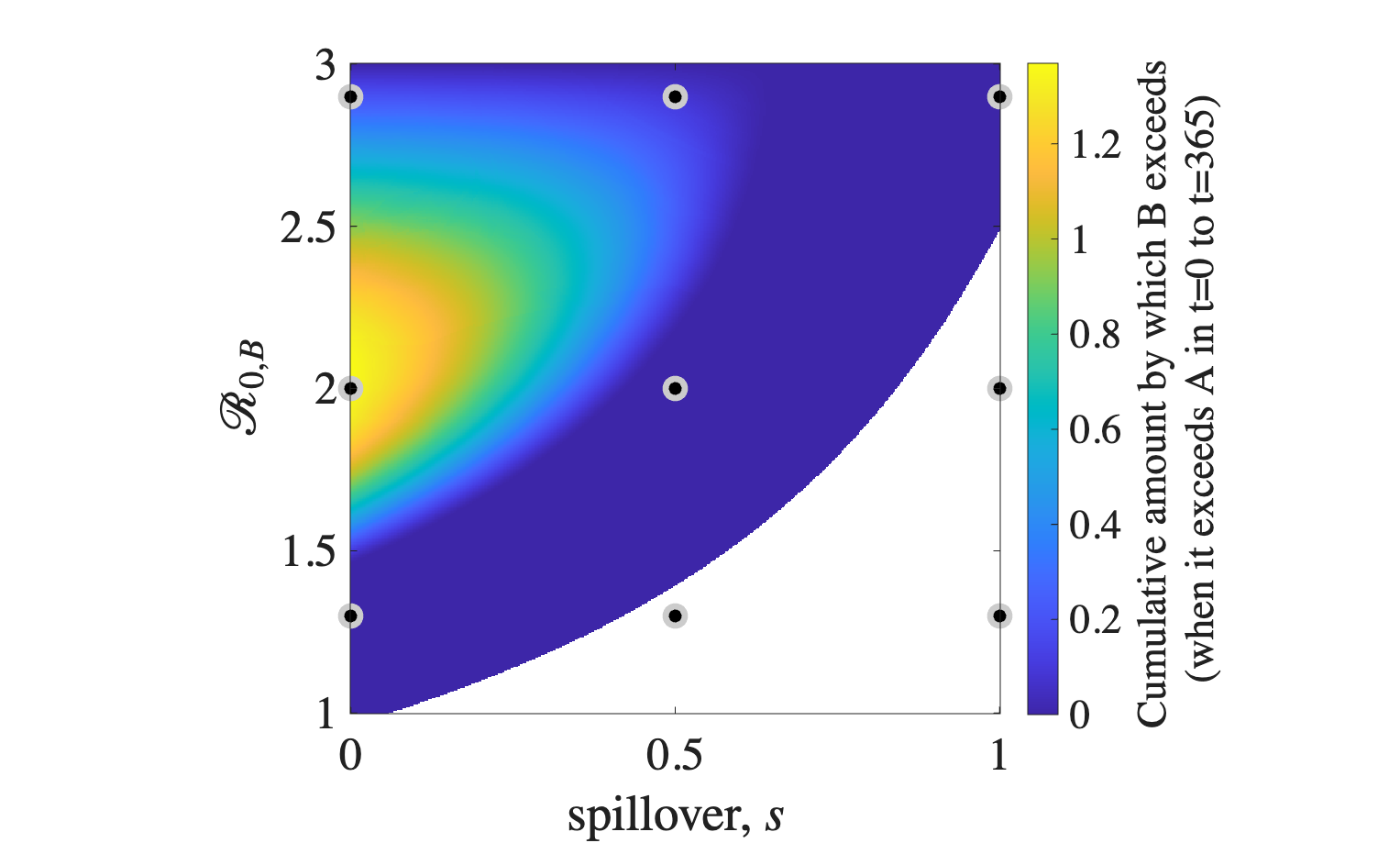}
    \caption{Amount disease $B$ exceeds disease  $A$}\label{fig:persistence_B_exceed_A_frac}
\end{subfigure}% 
    \caption{Persistence and superiority of diseases across one year as computed through numerical simulation using a fractional formulation for behavioral feedback (as opposed to exponential formulation in Figure \ref{fig:persistence_exceed}). (a) Persistence of both diseases (purple) or only disease $A$ (red), which has the higher basic reproduction number. (b) Percentage of the first year that disease $B$ prevalence is above disease $A$ prevalence. For $\mathcal{R}_{0,A}=\mathcal{R}_{0,B}=3$ the diseases produce independent, identical outbreaks so neither disease exceeds the other one (white line at the tope). White space at the bottom right corresponds to when there is persistence of only disease $A$.  (c) Cumulative amount during the first year that disease $A$ exceeds disease $B$. (d) Cumulative amount during the first year that disease $B$ exceeds disease $A$. Note the different color scaling in panels (b)-(d). Dots correspond to combination of spillover ($s$) and basic reproduction number of disease $B$ ($\mathcal{R}_{0,B}$) used for plots in Figure \ref{fig:dyanmics_spillover}. %Note that for $\mathcal{R}_{0,A}=\mathcal{R}_{0,B}=3$, the strains are independent and identical, regardless of the level of spillover, so neither are considered dominant.
    }
    \label{fig:persistence_exceed_frac}
\end{figure}

\end{document}